\documentclass[11pt,a4paper]{amsart}

\usepackage{amsmath}
\usepackage{amssymb}
\usepackage{latexsym}
\usepackage{xypic}
\usepackage{amscd,amsthm}

\usepackage{epsfig,epic}

\input{epsf}

\newcommand{\card}[1]{{\mid\! #1 \!\mid}}
\newcommand{\pro}[2]{\langle #1, #2 \rangle}

\def\C{{\mathbb C}}
\def\N{{\mathbb N}}

\def\P{{\mathbb P}}

\def\R{{\mathbb R}}
\def\Z{{\mathbb Z}}

\def\F{{\mathcal F}}

\def\V{{\mathcal V}}

\def\conv{{\rm conv}}

\def\Vol{{\rm Vol}}

\newtheorem*{theorem*}{Theorem}
\newtheorem*{proposition*}{Proposition}
\newtheorem*{corollary*}{Corollary}
\newtheorem{theorem}{Theorem}[section]
\newtheorem{definition}[theorem]{Definition}
\newtheorem{remark}[theorem]{Remark}

\newtheorem{corollary}[theorem]{Corollary}
\newtheorem{proposition}[theorem]{Proposition}
\newtheorem{algorithm}[theorem]{Algorithm}

\include{pst-plot}      \psset{unit=3pt}               \unitlength=3pt
\def\BP{\begin{picture}} \def\EP{\end{picture}}         

\def\subdef#1{\gdef\globalColor##1{##1}}     
\def\TC#1{\subdef{cmyk #1}}      \subdef{Black}

      \def\blue       {\TC{1. 1. 0. 0.}}
\def\black      {\TC{0. 0. 0. 1.}}      

\def\lab#1)#2#3{\put#1){\makebox(0,0)[#2]{\small #3}}}
\def\putlin#1,#2,#3,#4,#5){\put#1,#2){\line(#3,#4){#5}}} 
\def\putvec#1,#2,#3,#4,#5){\put#1,#2){\vector(#3,#4){#5}}}
        \def\Xlab#1 {\put(#1,#1){\drawline(-.5,.5)(.5,-.5)}}
        \def\Ylab#1 {\put(-#1,#1){\drawline(.5,.5)(-.5,-.5)}}

\newcount\hsum \newcount\hdif

\begin{document}
\title{Classification of toric Fano $5$-folds}

\author[Maximilian Kreuzer]{Maximilian Kreuzer}
\address{Institute for Theoretical Physics, Vienna University of Technology, Wiedner Hauptstrasse 8-10, 1040 Vienna, Austria}
\email{maximilian.kreuzer@tuwien.ac.at}

\author[Benjamin Nill]{Benjamin Nill}
\address{Research Group Lattice Polytopes, Freie Universit\"at Berlin, Arnimallee 3, 14195 Berlin, Germany}
\email{nill@math.fu-berlin.de}

\begin{abstract}
We obtain $866$ isomorphism classes of five-dimensional nonsingular toric Fano varieties using a computer program and 
the database of four-dimensional reflexive polytopes. The algorithm is based on the 
existence of facets of Fano polytopes having small integral distance from any vertex.
\end{abstract}

\maketitle

\section*{Introduction}

Nonsingular projective varieties over $\C$ are called Fano, if their anticanonical sheaf is ample. 
Many efforts were taken to classify nonsingular Fano varieties up to deformation, and there is now a complete list 
in dimensions two and three \cite{MM04}. In the toric case $d$-dimensional nonsingular toric Fano varieties 
were completely classified for $d \leq 4$ up to biregular isomorphisms, for $d\leq 3$ by Batyrev \cite{Bat82} and Watanabe and Watanabe 
\cite{WW82}, and for $d=4$ by Batyrev \cite{Bat99} and Sato \cite{Sat00}. They found $1$, $5$, $18$, $124$ isomorphism classes 
for $d=1,2,3,4$.

The combinatorial objects corresponding to $d$-dimensional nonsingular toric Fano varieties are $d$-dimensional {\em Fano polytopes}. 
These are lattice polytopes having the origin in its interior such that any facet is a simplex, and moreover, 
the vertices of any facet form a lattice basis. The fan of the corresponding toric variety is given by the cones over 
the faces of the polytope. Biregular isomorphism classes of nonsingular toric Fano varieties are in one-to-one 
correspondence to lattice isomorphism classes of Fano polytopes. 

In \cite{Bat99} Batyrev investigated the projection of $d$-dimensional Fano polytopes 
along a vertex, and showed that the image is a $(d-1)$-dimensional reflexive polytope. A {\em reflexive polytope} 
is a lattice polytope having the origin in its interior such that the origin has integral distance one from any facet. 
Reflexive polytopes correspond to Gorenstein toric Fano varieties and were used to construct mirror pairs of Calabi-Yau varieties \cite{Bat94}. 
Skarke and the first author succeeded in classifying all $d$-dimensional reflexive polytopes up to lattice isomorphisms 
for $d \leq 4$ using the computer program PALP \cite{KS98, KS00, KS02, KS04}.
They found $1$, $16$, $4319$, $473 800 776$ 
isomorphism classes for $d=1,2,3,4$.

In this article, we describe how to recover all $d$-dimensional Fano polytopes from the classification of $(d-1)$-dimensional reflexive polytopes.
The effectiveness of this procedure is based on the following key observation 
(for the precise and slightly stronger statement see Theorem \ref{facet}, Cor. \ref{curves}):

\begin{proposition*}
Let $X$ be a $d$-dimensional nonsingular toric Fano variety, and $X \not\cong \P^d$. Then 
there exists a facet $F$ of the associated $d$-dimensional Fano polytope $P$ such that any vertex of $P$ has integral distance 
at most $d$ from $F$. In particular, there exist $d$ torus-invariant 
curves on $X$, intersecting in a common torus fixpoint, all having anticanonical degree at most $d$.
\end{proposition*}

The idea of the classification is to take the vertices of the facet 
$F$ as a lattice basis $e_1, \ldots, e_d$ of $\Z^d$, and to 
project $P$ along $e_d$. Now, for any vertex of $P$ its first $d-1$ coordinates are determined by the 
projection, and the assumption on $F$ yields an effective bound on its last coordinate. 
We implemented the algorithm in PALP and obtained the following main result:

\begin{theorem*}
There are $866$ isomorphism classes of five-dimensional nonsingular toric Fano varieties.
\label{theo}
\end{theorem*}

The list of the five-dimensional Fano polytopes and their facet inequalities, 
in form of the vertices of the dual polytope, is available on the webpage \cite{KN07}.

\smallskip

The paper is organized in the following way:

In the first section we prove the existence of a facet that has small integral distance from any vertex, extending an argument 
due to Casagrande in \cite{Cas06}. In the second section we recall the necessary results about projections of Fano polytopes. 
The third section contains the description of the classification algorithm. 
In the last section we list some observations on the properties of the classified nonsingular toric Fano varieties, and we also 
give the Hodge numbers of associated complete intersection Calabi-Yau $3$-folds.

\smallskip

\textit{Remark:} After the publication of this work {\O}bro presented in \cite{Oebro07b} a direct algorithm for the classification of 
$d$-dimensional Fano polytopes, which he successfully applied for $d \leq 7$.

\vspace{-0.6ex}

\subsection*{Acknowledgments:}
The work of the first author was supported in part by the Austrian Research
Funds FWF grant Nr. P18679-N16. The second author is a member of the research group Lattice Polytopes, 
supported by Emmy Noether fellowship HA 4383/1 of the German Research Foundation. He also 
thanks the Vienna University of Technology for hospitality and support. 
We are grateful to the anonymous referees for valuable suggestions.

\vspace{-0.4ex}

\section{Facets having small integral distance from any vertex}

In this section we show that any Fano polytope $P$ has a facet $F$ 
such that any vertex of $P$ has small integral distance from $F$.

First let us recall the basic combinatorial notions.

Let $P \subseteq \R^d$ be a lattice polytope, i.e., its vertices are contained in the lattice $\Z^d$. 
Throughout we assume that $P$ is $d$-dimensional and contains the origin of the lattice in its interior. 
We also use the notion of a $d$-polytope to denote a $d$-dimensional polytope.

A facet of $P$ is called {\em simplex}, if the vertices of the facet are affinely independent, 
and the facet is called {\em unimodular}, if the vertices even form an affine lattice basis. 

\begin{definition}\ {\rm 
\begin{enumerate}
\item $P$ is called {\em reflexive}, if the {\em dual polytope}
\[P^* := \{y \in (\R^d)^* \;:\; \pro{y}{x} \geq -1 \;\forall\, x \in P\}\]
is a lattice polytope with respect to the dual lattice $(\Z^d)^*$.
\item $P$ is called {\em simplicial}, if any facet of $P$ is a simplex.
\item $P$ is called {\em Fano polytope}, if the vertices of any facet of $P$ form a lattice basis.
\end{enumerate}
}
\end{definition}

Hence, a Fano polytope is precisely a simplicial reflexive polytope, where all facets are unimodular.

\begin{remark}{\rm 
One should note that in the literature Fano polytopes are sometimes also called smooth Fano polytopes, e.g., 
\cite{Deb03}, \cite{Nil05}.
}
\end{remark}

Now, the starting point is the following bound on the number of vertices of a simplicial reflexive polytope, 
\cite{Cas06}:

\begin{theorem}[Casagrande]
A simplicial reflexive $d$-polytope $P$ has at most $3d$ vertices for $d$ even and 
$3d-1$ for $d$ odd. If $P$ has $3d$ vertices, then the associated nonsingular toric Fano variety is 
a $d/2$-product of $\P^2$ blown-up in three torus-invariant points.
\label{casa}
\end{theorem}

For our purposes we need a simplified and generalized version of the previous result, 
based on the idea of the proof of Casagrande. For this we first introduce some 
notation:

\begin{definition} {\rm We define:
\begin{enumerate}
\item $\V(P)$ is the set of vertices, $\F(P)$ the set of facets of $P$.
\item $\eta_F$ is the unique inner normal of $F \in \F(P)$ 
with $\pro{\eta_F}{F} = -1$. Hence, $\V(P^*) = \{\eta_F \,:\, F \in \F(P)\}$.
\item For $F \in \F(P)$ we define 
$a_F := \sum_{v \in \V(P)} (1 - \pro{\eta_F}{v})$.
\end{enumerate}
}
\end{definition}

Here is the main observation:

\begin{proposition}
Let $P$ be a lattice $d$-polytope containing the origin in its interior. 
Then $\card{\V(P)}$ is a proper convex combination of $\{a_F\}_{F \in \F(P)}$.
\label{casaprop}
\end{proposition}

\begin{proof}

Since also $P^*$ contains the origin in its interior, we can represent the origin as 
a proper convex combination of the vertices $\{\eta_F\}_{F \in \F(P)}$:
\[\sum_{F \in \F(P)} \lambda_F \eta_F = 0,\] 
where $\lambda_F > 0$ for any $F \in \F(P)$, and $\sum_{F \in \F(P)} \lambda_F = 1$. Hence 
\[\sum_{F \in \F(P)} \lambda_F a_F = 
\sum_{v \in \V(P)} \sum_{F \in \F(P)} \lambda_F (1 - \pro{\eta_F}{v}) = 
\sum_{v \in \V(P)} (1 - 0) = \card{\V(P)}.\]
\end{proof}

Recall that, if $P$ is a reflexive polytope and $F \in \F(P)$, then $\eta_F \in \V(P^*)$ is a primitive lattice point, hence 
$v \in \V(P)$ has {\em integral distance} $\pro{\eta_F}{v} + 1$ from the facet $F$. 

Now, we can derive a first result on the existence of a facet having small integral distance from any vertex. This 
result is essential for the effectiveness of the classification algorithm:

\begin{corollary}
Let $P$ be a simplicial reflexive $d$-polytope. Then there exists a facet $F$ such that 
\[\sum_{v \in \V(P) \,:\, \pro{\eta_F}{v} \geq 1} 
\pro{\eta_F}{v} \leq d.\]
In particular, any vertex of $P$ has integral distance at most $d+1$ from $F$.
\label{simplizial}
\end{corollary}

\begin{proof}

By Prop. \ref{casaprop} there is a facet $F \in \F(P)$ such that 
\[\card{\V(P)} \leq a_F = \card{\V(P)} - 
\sum_{v \in \V(P)} \pro{\eta_F}{v}.\]
 
Hence, since $F$ has $d$ vertices, 
\[0 \geq \sum_{v \in \V(P)} \pro{\eta_F}{v} = 
-d + \sum_{v \in \V(P) \,:\, \pro{\eta_F}{v} \geq 1} \pro{\eta_F}{v}.\]
\end{proof}

While the previous result gives a strong restriction on the number of vertices having integral distance $> 1$ from such a facet, 
the following result yields a bound on the remaining vertices, \cite[Remark 5(2)]{Deb03} and \cite[Lemma 5.5]{Nil05}:

\begin{proposition}
Let $P$ be a simplicial reflexive $d$-polytope. Then for any facet $F \in \F(P)$ there are at most $d$ vertices $v \in \V(P)$ with 
$\pro{\eta_F}{v} = 0$.
\label{zero}
\end{proposition}

In particular, we see that, if for any facet there is a vertex far away, then the overall number of vertices cannot be too large:

\begin{corollary}
Let $P$ be a simplicial reflexive $d$-polytope. Let $w$ be the maximal integral distance a vertex of $P$ has from the facet $F$ 
of Cor. \ref{simplizial}. Then $P$ has at most $3d + 2 - w$ vertices.
\end{corollary}

\begin{proof}

We have $\pro{\eta_F}{v} \in \{-1,0,1, \ldots, w-1\}$ for $v \in \V(P)$. Let $k_i := \card{\{v \in \V(P) \,:\, \pro{\eta_F}{v} =i\}}$ for 
$i = -1,0,1, \ldots, w-1$. By Prop. \ref{zero} $k_0 \leq d$. Hence, since $k_{w-1} > 0$, we get $\card{\V(P)} = \sum_{i=-1}^{w-1} k_i = d + k_0 + 
\sum_{i=1}^{w-1} k_i \leq 2 d + (\sum_{i=1}^{w-1} k_i i) - k_{w-1} (w-2) \leq 3 d - (w-2)$.
\end{proof}

In the case of a Fano polytope we can sharpen Cor. \ref{simplizial} even further (note that a $d$-dimensional Fano polytope is a simplex if 
and only if the associated nonsingular toric Fano variety is isomorphic to $\P^d$):

\begin{theorem}
Let $P$ be a Fano $d$-polytope. If $P$ is not a simplex, then 
there exists a facet $F$ such that 
\[\sum_{v \in \V(P) \,:\, \pro{\eta_F}{v} \geq 1} \pro{\eta_F}{v} \leq d\]
and any vertex of $P$ has integral distance at most $d$ from $F$.
\label{facet}
\end{theorem}

\begin{proof}

Let $F \in \F(P)$ be chosen as in Cor. \ref{simplizial}. The vertices of $F$ are denoted by $e_1, \ldots, e_d$. 
Assume that there is a vertex $v \in \V(P)$ with $\pro{\eta_F}{v} \geq d$. This implies 
$\pro{\eta_F}{v} = d$, and moreover, any other vertex $v' \in \V(P)$ with $v' \not= v$ and $v' \not\in F$ 
satisfies $\pro{\eta_F}{v'} = 0$. Now, since $P$ is not a simplex, there exists a subset $I \subsetneq \{1, \ldots, d\}$ 
such that $V := \{v, e_i\,:\, i \in I\}$ is not the vertex set of a face of $P$, however any proper subset of $V$ is. $V$ is called 
{\em primitive collection} in the language of Batyrev, \cite[2.6]{Bat91}. This implies that there exists a so-called {\em primitive relation} 
\cite[2.8,3.1]{Bat91}
\[v + \sum_{i \in I} e_i = \sum_{j=1}^r k_j v_j,\]
where $r \in \N_{\geq 0}$, $k_j \in \N_{>0}$ and $v_j \in \V(P)\backslash V$ for $j = 1, \ldots, r$. 
Hence, applying $\eta_F$ yields $\pro{\eta_F}{v_j} \in \{-1,0\}$ for $j = 1, \ldots, r$, thus
\[0 < d-\card{I} = \sum_{j=1}^r k_j \pro{\eta_F}{v_j} \leq 0,\]
a contradiction.
\end{proof}

In particular, from this result one may immediately derive the classification of toric del Pezzo surfaces.

\begin{remark}{\rm
Looking at the reflexive hexagon $P$, a Fano $2$-polytope, we see that one cannot in general expect the existence of a facet $F$
with 
\[\sum_{v \in \V(P) \,:\, \pro{\eta_F}{v} \geq 1} \pro{\eta_F}{v} \leq d-1.\]
However, such a facet does not exist if and only if 
$\sum_{v \in \V(P)} \pro{\eta_F}{v} \geq 0$ for all $F \in \F(P)$, which is equivalent to
\[\sum_{v \in \V(P)} v = 0.\]
}
\end{remark}

We finish with a "dual" application of Theorem \ref{facet} (compare with \cite[Prop. 2.1]{Lat96}):

\begin{corollary}
Let $X$ be a $d$-dimensional nonsingular toric Fano variety, and $X \not\cong \P^d$. Then there exist $d$ torus-invariant 
curves on $X$, intersecting in a common torus fixpoint, all having anticanonical degree at most $d$.
\label{curves}
\end{corollary}

\begin{proof}

Let $P$ be the associated Fano polytope. Let $F \in \F(P)$ be chosen as in Theorem \ref{facet}. 
Let $u \in \V(P^*)$ be one of the $d$ vertices joined with $\eta_F \in \V(P^*)$ by an edge, say $e$. 
Let $G$ be a facet of $P^*$ containing $u$ but not $\eta_F$. By duality, $\eta_G \in \V(P)$, and the assumption on $F$ 
yields that $\eta_F$ has integral distance $\pro{\eta_G}{\eta_F}+1 \leq d$ from $G$. 
Therefore, there can be at most $d+1$ lattice points on $e$. 
Hence, as is well-known (e.g., see the proof of Cor. 3.6 in \cite{Lat96}), 
the torus-invariant curve on $X$ corresponding to $e$ has anticanonical degree at most $d$.
\end{proof}

\section{Projecting Fano polytopes}

In this section we recall the basic notions and results about 
projecting Fano polytopes along vertices.

Here is the main observation, \cite[Prop. 2.4.4, Prop. 2.4.10]{Bat99}:

\begin{proposition}[Batyrev]
Let $P$ be a Fano $d$-polytope. Then the projection $P'$ of $P$ along 
a vertex $v$ is a reflexive $(d-1)$-polytope.

For any lattice point $x'$ in $P'$ there exists a vertex $x$ 
in $P$ mapping to $x'$. 

If $x'=0$, then $x = v$ or $x = -v$.

If $x' \not= 0$, then there are at most two 
vertices $x_1, x_2$ of $P$ mapping to $x'$. 
In this case $x_1 - x_2 = v$ or $x_2 - x_1 = v$.
\label{double}
\end{proposition}

This motivates the following definition, \cite[Def. 2.4.7]{Bat99}:

\begin{definition}{\rm
Let $P'$ be the projection of a Fano polytope along a vertex. 
A lattice point $x'$ in $P'$ with two vertices in $P$ mapping 
to $x'$ is called a {\em double point}.
}
\end{definition}

The next result \cite[Thm. 3.4]{Cas03} gives a strong obstruction on the existence of double points, and 
can also be generalized to simplicial reflexive polytopes, \cite{Nil05}:

\begin{theorem}[Casagrande]
Let $P'$ be the projection of a Fano polytope $P$ along a vertex $v$. 

Then there are at most two double points $\not= 0$. 

Equality can only hold, if these two double points $x',y'$ are centrally-symmetric. 
In this case, let us denote by $x$ the unique vertex of $P$ mapping to $x'$ such that $x$ is contained in a common 
face with $v$; in the same way we define $y$. Then $x+y=v$.
\label{casadouble}
\end{theorem}

Now, the following result shows that any facet of a projection of a Fano polytope is either unimodular or a circuit with unimodular triangulation:

\begin{proposition}[Batyrev]
Let $P'$ be the projection of a Fano $d$-polytope $P$ along a vertex $v$. 

If $C$ is a face of $P$ contained in a facet $F$ of $P$ such that $v \in F$ but $v \not\in C$, then 
the projection of $C$ is contained in a facet of $P'$.

On the other hand, let $G'$ be a facet of $P'$, and $G$ the face of $P$ mapping to $G'$, thus $v \not\in G$. There are two cases:

\begin{enumerate}
\item $G'$ has $d-1$ lattice points.

Then $G'$ is an unimodular $(d-2)$-simplex. 
There exists at most one double point on $G'$. 

If there is no double point on $G'$, then $G$ is $(d-2)$-dimensional. If there is one double point on $G'$, then 
$G$ is a facet of $P$.
\item $G'$ has $d$ lattice points $v'_1, \ldots, v'_d$.

Then there exists a relation among $v'_1, \ldots, v'_d$ of the following type:
\begin{equation}\sum_{i \in I} v'_i = \sum_{j \in J} k_j v'_j \;\; \text{ with }\;\; \sum_{j \in J} k_j = \card{I} = \Vol(G'),
\tag{$\star$}
\label{eq1}
\end{equation}
where $I, J$ are non-empty disjoint subsets of $\{1, \ldots, d\}$, and $k_j$ are positive natural numbers.

For any $i \in I$ the set 
$\{v'_1, \ldots, v'_d\}\backslash\{v'_i\}$ is a lattice basis. There are no double points on $G'$. 

Here, $G$ is a facet of $P$. Moreover, there is a {\em unique} relation of type (\ref{eq1}) such that 
projecting yields a one-to-one correspondence between the following two sets of 
cardinality $\card{I}$:
\begin{itemize}
\item[(a)] $(d-2)$-dimensional faces $C$ of $G$ such that the convex hull of $C$ and $v$ is a facet of $P$
\item[(b)] $(d-2)$-simplices with vertices $\{v'_1, \ldots, v'_d\}\backslash\{v'_i\}$ (for $i \in I$) 
\end{itemize}
\end{enumerate}
\label{proj}
\end{proposition}

\begin{proof}

The main results can be found in \cite[Cor. 2.4.6, Prop. 2.4.8]{Bat99}. 
Only the last statement is not explicitly stated in \cite{Bat99}. Here, the relation (\ref{eq1}) is given as the projection of 
a so-called primitive relation, see \cite[Def. 2.1.4, Cor. 2.4.6]{Bat99}. Namely, 
$v + \sum_{i \in I} v_i = \sum_{j \in J} k_j v_j$, with $\V(G) = \{v_1, \ldots, v_d\}$, where by definition 
$\{v\} \cup \{v_i \,:\, i \in I\}$ is not contained in a face of $P$. Hence, the 
projection of $\V(C)$, for $C$ as in (a), does not contain $\{v'_i \,:\, i \in I\}$. 
On the other hand, it follows from \cite[Thm. 2.3.3]{Bat99} and a well-known theorem of Reid, see \cite[Thm. 2.3]{Cas03}, 
that for any $i \in I$ the convex hull of $v$ and $\{v_1, \ldots, v_d\}\backslash\{v_i\}$  is a facet of $P$.
\end{proof}

\begin{remark}{\em The relation (\ref{eq1}) is called {\em circuit relation}, since $d$ lattice points, whose 
affine hull is $(d-2)$-dimensional, form a so-called {\em circuit}. This circuit relation is unique, if there is some 
$j \in J$ with $k_j> 1$. Otherwise, interchanging $I$ and $J$ yields another circuit relation of type (\ref{eq1}). 
}
\end{remark}

In a special case there is yet another restriction, \cite[Prop. 2.4.9]{Bat99}:

\begin{proposition}[Batyrev]
If in the situation of Prop. \ref{proj} the origin of $P'$ is a double point, then 
in the case of (2) we have $k_j=1$ for all $j \in J$.
\label{sym}
\end{proposition}

\section{The classification algorithm}

In this section we describe the algorithm 
which we used to classify all five-dimensional Fano polytopes up to lattice isomorphisms:

\begin{algorithm}\ 
\label{algo}

{\rm 
{\em Input:} A finite list of all isomorphism classes of reflexive $(d-1)$-polytopes.

{\em Output:} A finite list of all isomorphism classes of Fano $d$-polytopes.\\

\pagebreak

{\em Procedure:}
\begin{enumerate}
\item Output the Fano simplex associated to $\P^d$. If $d$ is even, then output the Fano 
polytope associated to the $d/2$-product of $\P^2$ blown-up in three torus-invariant points.
\item Determine all $(d-1)$-dimensional reflexive polytopes $P'$ with $\leq 3 d-1$ lattice points such that 
any facet of $P'$ is either unimodular or a circuit with a circuit relation of type (\ref{eq1}) in Prop. \ref{proj}(2). 
\item For any such chosen $P'$ proceed through all the facets $G'$ of $P'$. If $G'$ is unimodular, then define $C' := G'$. 
If $G' = \conv(v'_1, \ldots, v'_d)$ is a circuit, then for any (at most two) circuit relations of type (\ref{eq1}) and any  $i \in I$ 
(see (\ref{eq1})) define $C' := \conv(\{v'_1, \ldots, v'_d\}\backslash\{v'_i\})$.
\item For any such $C'$ perform a lattice transformation mapping the vertices of $C'$ onto the first $d-1$ basis vectors 
$e_1, \ldots, e_{d-1}$ of $\Z^{d-1}$.
\item Define $k(0)=k(e_1)= \cdots = k(e_{d-1}) = -1$. 
Now, proceed through all possible choices of natural numbers $k(w') \in \{0, \ldots, d-1\}$ 
assigned to the lattice points $w'$ of $P'$ with $w' \not\in \{0,e_1, \ldots, e_{d-1}\}$ satisfying
\[\sum_{w' \in (P' \cap \Z^{d-1})\backslash\{0, e_1, \ldots, e_{d-1}\}} k(w') \leq d.\] 
Hereby, one is allowed to choose $k(w') = 0$ only for at most $d$ lattice points $w'$.
\item For any such choice of $\{k(w')\}_{w' \in P' \cap \Z^{d-1}}$ 
define a list $V$ of lattice points in $\Z^d$ by assigning to any lattice point $w'$ of $P'$ 
a lattice point $w \in V$ having the same first $d-1$ coordinates as $w'$, and the last coordinate 
$w_d := -w_1 - \cdots - w_{d-1} - k(w')$. For instance, to the lattice point $0 \in P'$ this associates $e_d \in V$.
\item Finally, choose among the lattice points of $P'$ at most three so-called double points, where 
the number of lattice points of $P'$ plus the number of chosen double points must not exceed $3d-1$. 
Here, any non-zero double point must not be contained in a circuit facet of $P'$, 
and any unimodular facet of $P'$ may contain at most one double point. 
The origin $0$ may only be chosen as a double point, if for 
any circuit facet of $P'$ all coefficients in a circuit relation of type (\ref{eq1}) are equal to $1$.  
Moreover, two non-zero double points $w',u'$ may only be chosen, if $u'=-w'$ and the associated lattice points $w,u \in V$ 
satisfy $w+u=e_d$.
\item For any chosen non-zero double point $w'$ add the lattice point $w-e_d$ to the set $V$. If the origin is chosen as a double point, 
add $-e_d$ to $V$.
\item Check whether $\conv(V)$ is a Fano polytope. If yes, output $\conv(V)$.
\end{enumerate}
}
\end{algorithm}

\begin{proposition}
Algorithm \ref{algo} is correct.
\end{proposition}

\begin{proof}

Let $P$ be a Fano $d$-polytope. By (1) and Theorem \ref{casa} we may assume that $P$ is not a simplex and has $\leq 3d-1$ vertices. 
Let $F$ be a facet as in Theorem \ref{facet}. The vertices $e_1, \ldots, e_d$ of $F$ form a lattice basis of $\Z^d$. 
By Prop. \ref{double} projecting $P$ along $e_d$ yields a reflexive $(d-1)$-polytope $P'$, that has at most $3d-1$ lattice points by Prop. \ref{double}. 
Therefore the isomorphism type of $P'$ appears in the list of (2). 

Let $C := \conv(e_2, \ldots, e_d)$. Then by Prop. \ref{proj} projecting $C$ yields a $(d-2)$-dimensional simplex $C'$ contained in a facet $G'$ of $P'$. 
In the situation of Prop. \ref{proj} either $C' = G'$ or $C' = \conv(\{v'_1, \ldots, v'_d\}\backslash\{v'_i\})$ for some $i \in I$ in a circuit 
relation of type (\ref{eq1}). 
Now, let $w$ be a vertex of $P$, with $w+e_d \not\in P$, projecting to a lattice point $w'$. Then 
\[k(w') := \pro{\eta_F}{w} = - w'_1 - \cdots - w'_{d-1} - w_d.\]
Moreover, by the assumption on $F$ in Theorem \ref{facet} $k(w') \leq d-1$ for all vertices $w$, 
the sum of all $k(w') \geq 0$ does not exceed $d$, and by Prop. \ref{zero} we have $k(w') = 0$ for at most $d$ vertices. 

Hence, for these choices of $P',G',C', \{k(w')\}$ in steps (3)--(5) the set $V$ in (6) contains precisely all vertices $w$ of $P$ with $w+e_d \not\in P$. 
Now, by Prop. \ref{double} any other vertex of $P$ projects onto a double point. 
The restrictions on the double points of $P'$ in (7) 
follow from Theorem \ref{casadouble} and Prop. \ref{sym}. Hence, choosing these double points in step (8) yields $V = \V(P)$, thus 
the Algorithm \ref{algo} outputs $P$ in step (9).
\end{proof}

For $d=5$, we first determined all $4605$ isomorphism classes of reflexive $4$-polytopes satisfying the conditions in (2). 
Then the remaining computation steps of the algorithm took less than an hour on a standard PC. Note, that the program PALP, \cite{KS04}, 
can determine for any lattice polytope a normal form, which makes it easy to check whether two lattice polytopes 
are lattice isomorphic. Using this feature we determined the list of all $866$ mutually non-isomorphic Fano $5$-polytopes. 

\begin{remark}{\em 
Here are three comments in favour of the validity of our computer calculations: 

First, the algorithm was implemented for arbitrary $d$. For $d \leq 4$ we regained all $d$-dimensional nonsingular toric Fano varieties, which were 
previously classified using completely different methods.

Second, in \cite{Sat06} Sato classified all five-dimensional nonsingular toric Fano varieties with index $2$. He found $10$ isomorphism 
classes. In our list of isomorphism classes one variety had index six (namely, $\P^5$), one had index three (namely, $\P_{\P^3}(O_{\P^3} \oplus 
O_{\P^3} \oplus O_{\P^3}(1))$), 
and $11$ varieties had index two. It turned out that Sato had missed one five-dimensional nonsingular toric Fano variety with index $2$. 
In the notation of \cite[Thm. 3.6]{Sat06}, this variety is the $\P^1$-bundle over $\P_{\P^3}(O_{\P^3} \oplus O_{\P^3}(1))$ whose primitive relations are 
$x_1+x_2+x_3+x_4=x_5+x_7$, $x_5+x_6=0$ and $x_7+x_8=0$, where $G(\Sigma)=\{x_1, \ldots, x_8\}$.

Third, 
the classification of $866$ five-dimensional Fano polytopes was recently also independently achieved by {\O}bro using a completely different 
algorithm, see \cite{Oebro07b}.
}
\end{remark}

\section{Classification results}

In this section we list some of the properties of $d$-dimensional nonsingular toric Fano varieties for $d \leq 5$. Many of these observations give 
rise to conjectures for general $d$. For the toric dictionary see \cite{Ful93}.

\subsection{Maximal Picard number}

The Picard number of a $d$-dimensional nonsingular toric Fano variety equals the number of vertices of the associated Fano polytope minus $d$. 
Let $S_2$ be $\P^2$ blown-up in two torus-invariant fixpoints, and $S_3$ be $\P^2$ blown-up in three torus-invariant fixpoints. 

\begin{proposition}
Let $d \leq 5$, and $X$ be a $d$-dimensional nonsingular toric Fano variety with Picard number $\rho_X$.

If $d$ is even, then $\rho_X \leq 2d$, and there is up to isomorphisms only one $X$ with 
$\rho_X = 2 d$, namely $(S_3)^{d/2}$, and one with $\rho_X = 2d-1$, namely $S_2 \times (S_3)^{(d-2)/2}$.

If $d$ is odd, then $\rho_X \leq 2d-1$, and there are up to isomorphisms precisely two $X$ with 
$\rho_X = 2 d-1$, namely $\P^1 \times (S_3)^{(d-1)/2}$ or a uniquely determined toric $(S_3)^{(d-1)/2}$-fiber bundle over $\P^1$.
\label{vertices}
\end{proposition}

This result has a classification-free proof for $d=5$ due to Casagrande, see \cite[Thm. 4.2]{Cas03}.

\begin{remark}{\em 
A complete proof of Prop. \ref{vertices} for general $d$ was announced by the second author in the preprint of this article, {\texttt math.AG/0702890}. 
In the meantime a proof has been published by {\O}bro, \cite{Oebro07a}.
}
\end{remark}

\subsection{Maximal Euler characteristic}

The topological Euler characteristic of a nonsingular toric Fano variety equals the number of facets of the associated Fano polytope.

\begin{proposition}
Let $d \leq 5$, and $X$ be a $d$-dimensional nonsingular toric Fano variety. 
Then the topological Euler characteristic $\chi(X)$ of $X$ is maximal if and only if the Picard number $\rho_X$ is maximal. 

If $d$ is even, then $\chi(X) \leq 6^{d/2}$.

If $d$ is odd, then $\chi(X) \leq 2 \,\cdot\, 6^{(d-1)/2}$.
\end{proposition}

\begin{remark}{\em Since the dual of a Fano polytope is a reflexive polytope, 
the previous result should be seen in view of the open conjecture 
that any $d$-dimensional reflexive polytope has at most $6^{d/2}$ vertices, with equality 
only if $d$ is even and only for the $d/2$-product of the reflexive hexagon, \cite{Nil05}. 
}
\end{remark}

\subsection{Maximal degree}

The anticanonical degree of a $d$-dimensional nonsingular toric Fano variety equals the lattice volume of the dual of the associated Fano polytope.

\begin{proposition}
Let $d \leq 5$, and $X$ be a $d$-dimensional nonsingular toric Fano variety. 
Then there is precisely one $X$ with maximal anticanonical degree $(-K_X)^d$:\\

\begin{tabular}{l||l|l|l|l}
Dimension $d$ & $2$ & $3$ & $4$ & $5$\\
\cline{1-5}Maximal $(-K_X)^d$ & $9$ & $64$ & $800$ & $14762$\\
Unique $X$ & $\P^2$ & $\P^3$ & $\P_{\P^3}(O_{\P^3} \oplus O_{\P^3}(3))$ & $\P_{\P^4}(O_{\P^4} \oplus O_{\P^4}(4))$
\end{tabular}
\end{proposition}

\begin{remark}{\rm It is an open problem to determine the maximal degree of a nonsingular toric Fano variety. Projective bundles with large degree 
are constructed in \cite[Prop. 5.22]{Deb01}. An upper bound is proven in \cite[Thm. 9]{Deb03}.
}
\end{remark}

\subsection{Maximal dimension of the space of global sections of the anticanonical sheaf}

The dimension of the space of global sections of the anticanonical sheaf of a 
nonsingular toric Fano variety equals the number of lattice points of the dual of the associated Fano polytope.

\begin{proposition}
Let $d \leq 5$, and $X$ be a $d$-dimensional nonsingular toric Fano variety. 
Then $h^0(X,O_X(-K_X))$ is maximal if and only if the anticanonical degree $(-K_X)^d$ is maximal.\\

\begin{tabular}{l||l|l|l|l}
Dimension $d$ & $2$ & $3$ & $4$ & $5$\\
\cline{1-5}Maximal $h^0(X,O_X(-K_X))$ & $10$ & $35$ & $159$ & $846$
\end{tabular}
\end{proposition}

\subsection{Maximal degree of torus-invariant curves}

There is a one-to-one correspondence between torus-invariant curves on a nonsingular toric Fano variety and 
the edges of the dual of the associated Fano polytope. Thereby, the anticanonical degree of a torus-invariant curve equals the edge length, 
i.e., the number of lattice points minus one, of the corresponding edge.

\begin{proposition}
Let $w(d)$ be the maximal anticanonical degree $-K_X . C$ of a torus-invariant curve on a $d$-dimensional nonsingular toric Fano variety. 
Then $w(d)=3,5,7,11$ for $d=2,3,4,5$.
\end{proposition}

One can construct a nonsingular toric Fano variety that is a 
$\P^2$-bundle over a $\P^2$-bundle over $\P^2$ such that $w(6) \geq 15$. However, the authors are not 
aware of any reasonable upper bound on $w(d)$ giving the presumably correct asymptotics.

\subsection{Embedding of Fano polytopes}

There is a tempting question about Fano polytopes due to Ewald \cite{Ewa88} in 1988. 
Can one always embed any $d$-dimensional Fano polytope as a lattice polytope 
into $[-1,1]^d$? There is still no counterexample known.

\begin{proposition}
Let $d \leq 5$, and $P \subseteq \R^d$ a $d$-dimensional Fano polytope. Then there exists a lattice basis of $\Z^d$ such that 
any vertex of $P$ has coordinates in $\{-1,0,1\}$.
\end{proposition}

This proposition was checked by finding in the dual polytope $P^*$ 
lattice points $\pm e_1, \ldots, \pm e_d$ such that $e_1, \ldots, e_d$ is a lattice basis.

\subsection{Hodge numbers of complete intersection Calabi-Yau $3$-folds}

Recall that a projective variety with certain mild singularities is called a Calabi-Yau variety, 
if it has trivial canonical bundle, see \cite[Def. 4.1.8]{Bat94}. 
By a construction of Batyrev and Borisov \cite{BB96a} one can define complete 
intersection 
Calabi-Yau varieties (CICY) of codimension $d-r$ in a $d$-dimensional Gorenstein toric Fano variety 
using special Minkowski summands of length $r$, so called nef-partitions, of the dual of the associated reflexive polytope. 
Moreover, they gave combinatorial formulas for determining the stringy 
Hodge numbers of the CICYs \cite{BB96b,KRS}.

The five-dimensional Fano polytopes have between $3$ and $159$ nef-partitions of length two per dual Fano polytope, 
which can be used to construct $36233$ CICY 3-folds with $421$ different pairs of Hodge numbers. 
Here is a plot of $(h_{1 1}, h_{2 1})$, where $59$ of these pairs do not 
come from anticanonical hypersurfaces in four-dimensional Gorenstein toric 
Fano varieties:

\def\FanoCICYs{
\h1.73.\h1.79.\h1.89.\h19.19.\h2.100.\h2.102.\h2.112.\h2.116.\h2.122.
\h2.128.\h2.56.\h2.58.\h2.59.\h2.60.\h2.62.\h2.64.\h2.66.\h2.68.\h2.70.\h2.72.
\h2.74.\h2.76.\h2.77.\h2.78.\h2.80.\h2.82.\h2.83.\h2.86.\h2.90.\h2.92.\h2.95.
\h3.101.\h3.103.\h3.105.\h3.107.\h3.111.\h3.113.\h3.123.\h3.39.\h3.41.\h3.43.
\h3.44.\h3.45.\h3.47.\h3.48.\h3.49.\h3.50.\h3.51.\h3.52.\h3.53.\h3.54.\h3.55.
\h3.56.\h3.57.\h3.58.\h3.59.\h3.60.\h3.61.\h3.62.\h3.63.\h3.64.\h3.65.\h3.66.
\h3.67.\h3.68.\h3.69.\h3.70.\h3.71.\h3.72.\h3.73.\h3.75.\h3.76.\h3.77.\h3.78.
\h3.79.\h3.80.\h3.81.\h3.83.\h3.84.\h3.85.\h3.87.\h3.89.\h3.91.\h3.93.\h3.95.
\h3.99.\h4.100.\h4.101.\h4.102.\h4.104.\h4.112.\h4.121.\h4.34.\h4.36.\h4.37.
\h4.38.\h4.39.\h4.40.\h4.41.\h4.42.\h4.43.\h4.44.\h4.45.\h4.46.\h4.47.\h4.48.
\h4.49.\h4.50.\h4.51.\h4.52.\h4.53.\h4.54.\h4.55.\h4.56.\h4.57.\h4.58.\h4.59.
\h4.60.\h4.61.\h4.62.\h4.63.\h4.64.\h4.65.\h4.66.\h4.67.\h4.68.\h4.69.\h4.70.
\h4.71.\h4.72.\h4.73.\h4.74.\h4.75.\h4.76.\h4.77.\h4.78.\h4.79.\h4.80.\h4.81.
\h4.82.\h4.83.\h4.84.\h4.85.\h4.86.\h4.88.\h4.89.\h4.90.\h4.91.\h4.92.\h4.93.
\h4.94.\h4.96.\h4.97.\h5.101.\h5.29.\h5.31.\h5.32.\h5.33.\h5.34.\h5.35.\h5.36.
\h5.37.\h5.38.\h5.39.\h5.40.\h5.41.\h5.42.\h5.43.\h5.44.\h5.45.\h5.46.\h5.47.
\h5.48.\h5.49.\h5.50.\h5.51.\h5.52.\h5.53.\h5.54.\h5.55.\h5.56.\h5.57.\h5.58.
\h5.59.\h5.60.\h5.61.\h5.62.\h5.63.\h5.64.\h5.65.\h5.66.\h5.67.\h5.68.\h5.69.
\h5.70.\h5.71.\h5.72.\h5.73.\h5.74.\h5.75.\h5.76.\h5.77.\h5.78.\h5.79.\h5.80.
\h5.81.\h5.82.\h5.83.\h5.85.\h5.87.\h5.89.\h5.91.\h5.93.\h6.26.\h6.28.\h6.29.
\h6.30.\h6.31.\h6.32.\h6.33.\h6.34.\h6.35.\h6.36.\h6.37.\h6.38.\h6.39.\h6.40.
\h6.41.\h6.42.\h6.43.\h6.44.\h6.45.\h6.46.\h6.47.\h6.48.\h6.49.\h6.50.\h6.51.
\h6.52.\h6.53.\h6.54.\h6.55.\h6.56.\h6.57.\h6.58.\h6.59.\h6.60.\h6.61.\h6.62.
\h6.63.\h6.64.\h6.65.\h6.66.\h6.67.\h6.68.\h6.69.\h6.70.\h6.71.\h6.72.\h6.74.
\h6.75.\h6.76.\h6.78.\h6.80.\h6.81.\h6.84.\h6.90.\h7.23.\h7.25.\h7.27.\h7.28.
\h7.29.\h7.30.\h7.31.\h7.32.\h7.33.\h7.34.\h7.35.\h7.36.\h7.37.\h7.38.\h7.39.
\h7.40.\h7.41.\h7.42.\h7.43.\h7.44.\h7.45.\h7.46.\h7.47.\h7.48.\h7.49.\h7.50.
\h7.51.\h7.52.\h7.53.\h7.54.\h7.55.\h7.57.\h7.58.\h7.59.\h7.60.\h7.61.\h7.62.
\h7.63.\h7.65.\h7.67.\h7.69.\h7.71.\h7.75.\h7.79.\h8.22.\h8.24.\h8.25.\h8.26.
\h8.27.\h8.28.\h8.29.\h8.30.\h8.31.\h8.32.\h8.33.\h8.34.\h8.35.\h8.36.\h8.37.
\h8.38.\h8.39.\h8.40.\h8.41.\h8.42.\h8.43.\h8.44.\h8.45.\h8.46.\h8.47.\h8.48.
\h8.49.\h8.50.\h8.51.\h8.52.\h8.54.\h8.56.\h8.57.\h8.58.\h8.59.\h8.60.\h8.62.
\h8.64.\h8.68.\h8.74.\h9.21.\h9.23.\h9.24.\h9.25.\h9.27.\h9.29.\h9.30.\h9.31.
\h9.32.\h9.33.\h9.35.\h9.36.\h9.37.\h9.38.\h9.39.\h9.41.\h9.42.\h9.43.\h9.45.
\h9.47.\h9.49.\h9.51.\h9.53.\h9.57.
\h10.24.\h10.26.\h10.27.\h10.28.\h10.29.\h10.30.\h10.31.\h10.32.\h10.34.
\h10.36.\h10.37.\h10.38.\h10.40.\h10.42.\h10.43.\h10.46.\h10.52.\h1.101.
\h11.23.\h11.25.\h11.27.\h11.28.\h1.129.\h11.29.\h11.31.\h11.33.\h11.35.
\h11.37.\h11.39.\h11.41.\h1.149.\h12.20.\h12.22.\h12.24.\h12.26.\h12.27.
\h12.28.\h12.29.\h12.30.\h12.32.\h13.21.\h13.23.\h13.24.\h13.25.\h14.20.
\h21.21.
}

\def\HypersurfCY{
\h14.14.\h15.15.\h16.12.\h16.13.\h16.14.\h16.15.\h16.16.\h17.12.
\h17.13.\h17.14.\h17.15.\h17.16.\h17.17.\h18.10.\h18.11.\h18.12.
\h18.13.\h18.14.\h18.15.\h18.16.\h18.17.\h18.18.\h19.7.\h19.9.
\h19.11.\h19.13.\h19.14.\h19.15.\h19.16.\h19.17.\h19.18.\h19.19.
\h20.5.\h20.10.\h20.11.\h20.12.\h20.13.\h20.14.\h20.15.\h20.16.
\h20.17.\h20.18.\h20.19.\h20.20.\h21.1.\h21.9.\h21.11.\h21.12.
\h21.13.\h21.14.\h21.15.\h21.16.\h21.17.\h21.18.\h21.19.\h21.20.
\h21.21.\h22.8.\h22.9.\h22.10.\h22.11.\h22.12.\h22.13.\h22.14.
\h22.15.\h22.16.\h22.17.\h22.18.\h22.19.\h22.20.\h22.21.\h22.22.
\h23.7.\h23.9.\h23.10.\h23.11.\h23.12.\h23.13.\h23.14.\h23.15.
\h23.16.\h23.17.\h23.18.\h23.19.\h23.20.\h23.21.\h23.22.\h23.23.
\h24.8.\h24.9.\h24.10.\h24.11.\h24.12.\h24.13.\h24.14.\h24.15.
\h24.16.\h24.17.\h24.18.\h24.19.\h24.20.\h24.21.\h24.22.\h24.23.
\h24.24.\h25.7.\h25.8.\h25.9.\h25.10.\h25.11.\h25.12.\h25.13.
\h25.14.\h25.15.\h25.16.\h25.17.\h25.18.\h25.19.\h25.20.\h25.21.
\h25.22.\h25.23.\h25.24.\h25.25.\h26.6.\h26.8.\h26.9.\h26.10.
\h26.11.\h26.12.\h26.13.\h26.14.\h26.15.\h26.16.\h26.17.\h26.18.
\h26.19.\h26.20.\h26.21.\h26.22.\h26.23.\h26.24.\h26.25.\h26.26.
\h27.6.\h27.7.\h27.8.\h27.9.\h27.10.\h27.11.\h27.12.\h27.13.
\h27.14.\h27.15.\h27.16.\h27.17.\h27.18.\h27.19.\h27.20.\h27.21.
\h27.22.\h27.23.\h27.24.\h27.25.\h27.26.\h27.27.\h28.4.\h28.6.
\h28.7.\h28.8.\h28.9.\h28.10.\h28.11.\h28.12.\h28.13.\h28.14.
\h28.15.\h28.16.\h28.17.\h28.18.\h28.19.\h28.20.\h28.21.\h28.22.
\h28.23.\h28.24.\h28.25.\h28.26.\h28.27.\h28.28.\h29.2.\h29.5.
\h29.7.\h29.8.\h29.9.\h29.10.\h29.11.\h29.12.\h29.13.\h29.14.
\h29.15.\h29.16.\h29.17.\h29.18.\h29.19.\h29.20.\h29.21.\h29.22.
\h29.23.\h29.24.\h29.25.\h29.26.\h29.27.\h29.28.\h29.29.\h30.6.
\h30.7.\h30.8.\h30.9.\h30.10.\h30.11.\h30.12.\h30.13.\h30.14.
\h30.15.\h30.16.\h30.17.\h30.18.\h30.19.\h30.20.\h30.21.\h30.22.
\h30.23.\h30.24.\h30.25.\h30.26.\h30.27.\h30.28.\h30.29.\h30.30.
\h31.5.\h31.7.\h31.8.\h31.9.\h31.10.\h31.11.\h31.12.\h31.13.
\h31.14.\h31.15.\h31.16.\h31.17.\h31.18.\h31.19.\h31.20.\h31.21.
\h31.22.\h31.23.\h31.24.\h31.25.\h31.26.\h31.27.\h31.28.\h31.29.
\h31.30.\h31.31.\h32.6.\h32.7.\h32.8.\h32.9.\h32.10.\h32.11.
\h32.12.\h32.13.\h32.14.\h32.15.\h32.16.\h32.17.\h32.18.\h32.19.
\h32.20.\h32.21.\h32.22.\h32.23.\h32.24.\h32.25.\h32.26.\h32.27.
\h32.28.\h32.29.\h32.30.\h32.31.\h32.32.\h33.5.\h33.6.\h33.7.
\h33.8.\h33.9.\h33.10.\h33.11.\h33.12.\h33.13.\h33.14.\h33.15.
\h33.16.\h33.17.\h33.18.\h33.19.\h33.20.\h33.21.\h33.22.\h33.23.
\h33.24.\h33.25.\h33.26.\h33.27.\h33.28.\h33.29.\h33.30.\h33.31.
\h33.32.\h33.33.\h34.4.\h34.6.\h34.7.\h34.8.\h34.9.\h34.10.
\h34.11.\h34.12.\h34.13.\h34.14.\h34.15.\h34.16.\h34.17.\h34.18.
\h34.19.\h34.20.\h34.21.\h34.22.\h34.23.\h34.24.\h34.25.\h34.26.
\h34.27.\h34.28.\h34.29.\h34.30.\h34.31.\h34.32.\h34.33.\h34.34.
\h35.5.\h35.7.\h35.8.\h35.9.\h35.10.\h35.11.\h35.12.\h35.13.
\h35.14.\h35.15.\h35.16.\h35.17.\h35.18.\h35.19.\h35.20.\h35.21.
\h35.22.\h35.23.\h35.24.\h35.25.\h35.26.\h35.27.\h35.28.\h35.29.
\h35.30.\h35.31.\h35.32.\h35.33.\h35.34.\h35.35.\h36.4.\h36.6.
\h36.7.\h36.8.\h36.9.\h36.10.\h36.11.\h36.12.\h36.13.\h36.14.
\h36.15.\h36.16.\h36.17.\h36.18.\h36.19.\h36.20.\h36.21.\h36.22.
\h36.23.\h36.24.\h36.25.\h36.26.\h36.27.\h36.28.\h36.29.\h36.30.
\h36.31.\h36.32.\h36.33.\h36.34.\h36.35.\h36.36.\h37.4.\h37.5.
\h37.6.\h37.7.\h37.8.\h37.9.\h37.10.\h37.11.\h37.12.\h37.13.
\h37.14.\h37.15.\h37.16.\h37.17.\h37.18.\h37.19.\h37.20.\h37.21.
\h37.22.\h37.23.\h37.24.\h37.25.\h37.26.\h37.27.\h37.28.\h37.29.
\h37.30.\h37.31.\h37.32.\h37.33.\h37.34.\h37.35.\h37.36.\h37.37.
\h38.2.\h38.5.\h38.6.\h38.7.\h38.8.\h38.9.\h38.10.\h38.11.
\h38.12.\h38.13.\h38.14.\h38.15.\h38.16.\h38.17.\h38.18.\h38.19.
\h38.20.\h38.21.\h38.22.\h38.23.\h38.24.\h38.25.\h38.26.\h38.27.
\h38.28.\h38.29.\h38.30.\h38.31.\h38.32.\h38.33.\h38.34.\h38.35.
\h38.36.\h38.37.\h38.38.\h39.5.\h39.6.\h39.7.\h39.8.\h39.9.
\h39.10.\h39.11.\h39.12.\h39.13.\h39.14.\h39.15.\h39.16.\h39.17.
\h39.18.\h39.19.\h39.20.\h39.21.\h39.22.\h39.23.\h39.24.\h39.25.
\h39.26.\h39.27.\h39.28.\h39.29.\h39.30.\h39.31.\h39.32.\h39.33.
\h39.34.\h39.35.\h39.36.\h39.37.\h39.38.\h39.39.\h40.4.\h40.5.
\h40.6.\h40.7.\h40.8.\h40.9.\h40.10.\h40.11.\h40.12.\h40.13.
\h40.14.\h40.15.\h40.16.\h40.17.\h40.18.\h40.19.\h40.20.\h40.21.
\h40.22.\h40.23.\h40.24.\h40.25.\h40.26.\h40.27.\h40.28.\h40.29.
\h40.30.\h40.31.\h40.32.\h40.33.\h40.34.\h40.35.\h40.36.\h40.37.
\h40.38.\h40.39.\h40.40.\h41.5.\h41.6.\h41.7.\h41.8.\h41.9.
\h41.10.\h41.11.\h41.12.\h41.13.\h41.14.\h41.15.\h41.16.\h41.17.
\h41.18.\h41.19.\h41.20.\h41.21.\h41.22.\h41.23.\h41.24.\h41.25.
\h41.26.\h41.27.\h41.28.\h41.29.\h41.30.\h41.31.\h41.32.\h41.33.
\h41.34.\h41.35.\h41.36.\h41.37.\h41.38.\h41.39.\h41.40.\h41.41.
\h42.4.\h42.5.\h42.6.\h42.7.\h42.8.\h42.9.\h42.10.\h42.11.
\h42.12.\h42.13.\h42.14.\h42.15.\h42.16.\h42.17.\h42.18.\h42.19.
\h42.20.\h42.21.\h42.22.\h42.23.\h42.24.\h42.25.\h42.26.\h42.27.
\h42.28.\h42.29.\h42.30.\h42.31.\h42.32.\h42.33.\h42.34.\h42.35.
\h42.36.\h42.37.\h42.38.\h42.39.\h42.40.\h42.41.\h42.42.\h43.3.
\h43.5.\h43.6.\h43.7.\h43.8.\h43.9.\h43.10.\h43.11.\h43.12.
\h43.13.\h43.14.\h43.15.\h43.16.\h43.17.\h43.18.\h43.19.\h43.20.
\h43.21.\h43.22.\h43.23.\h43.24.\h43.25.\h43.26.\h43.27.\h43.28.
\h43.29.\h43.30.\h43.31.\h43.32.\h43.33.\h43.34.\h43.35.\h43.36.
\h43.37.\h43.38.\h43.39.\h43.40.\h43.41.\h43.42.\h43.43.\h44.4.
\h44.5.\h44.6.\h44.7.\h44.8.\h44.9.\h44.10.\h44.11.\h44.12.
\h44.13.\h44.14.\h44.15.\h44.16.\h44.17.\h44.18.\h44.19.\h44.20.
\h44.21.\h44.22.\h44.23.\h44.24.\h44.25.\h44.26.\h44.27.\h44.28.
\h44.29.\h44.30.\h44.31.\h44.32.\h44.33.\h44.34.\h44.35.\h44.36.
\h44.37.\h44.38.\h44.39.\h44.40.\h44.41.\h44.42.\h44.43.\h44.44.
\h45.3.\h45.5.\h45.6.\h45.7.\h45.8.\h45.9.\h45.10.\h45.11.
\h45.12.\h45.13.\h45.14.\h45.15.\h45.16.\h45.17.\h45.18.\h45.19.
\h45.20.\h45.21.\h45.22.\h45.23.\h45.24.\h45.25.\h45.26.\h45.27.
\h45.28.\h45.29.\h45.30.\h45.31.\h45.32.\h45.33.\h45.34.\h45.35.
\h45.36.\h45.37.\h45.38.\h45.39.\h45.40.\h45.41.\h45.42.\h45.43.
\h45.44.\h45.45.\h46.4.\h46.5.\h46.6.\h46.7.\h46.8.\h46.9.
\h46.10.\h46.11.\h46.12.\h46.13.\h46.14.\h46.15.\h46.16.\h46.17.
\h46.18.\h46.19.\h46.20.\h46.21.\h46.22.\h46.23.\h46.24.\h46.25.
\h46.26.\h46.27.\h46.28.\h46.29.\h46.30.\h46.31.\h46.32.\h46.33.
\h46.34.\h46.35.\h46.36.\h46.37.\h46.38.\h46.39.\h46.40.\h46.41.
\h46.42.\h46.43.\h46.44.\h46.45.\h46.46.\h47.5.\h47.6.\h47.7.
\h47.8.\h47.9.\h47.10.\h47.11.\h47.12.\h47.13.\h47.14.\h47.15.
\h47.16.\h47.17.\h47.18.\h47.19.\h47.20.\h47.21.\h47.22.\h47.23.
\h47.24.\h47.25.\h47.26.\h47.27.\h47.28.\h47.29.\h47.30.\h47.31.
\h47.32.\h47.33.\h47.34.\h47.35.\h47.36.\h47.37.\h47.38.\h47.39.
\h47.40.\h47.41.\h47.42.\h47.43.\h47.44.\h47.45.\h47.46.\h47.47.
\h48.4.\h48.5.\h48.6.\h48.7.\h48.8.\h48.9.\h48.10.\h48.11.
\h48.12.\h48.13.\h48.14.\h48.15.\h48.16.\h48.17.\h48.18.\h48.19.
\h48.20.\h48.21.\h48.22.\h48.23.\h48.24.\h48.25.\h48.26.\h48.27.
\h48.28.\h48.29.\h48.30.\h48.31.\h48.32.\h48.33.\h48.34.\h48.35.
\h48.36.\h48.37.\h48.38.\h48.39.\h48.40.\h48.41.\h48.42.\h48.43.
\h48.44.\h48.45.\h48.46.\h48.47.\h48.48.\h49.4.\h49.5.\h49.6.
\h49.7.\h49.8.\h49.9.\h49.10.\h49.11.\h49.12.\h49.13.\h49.14.
\h49.15.\h49.16.\h49.17.\h49.18.\h49.19.\h49.20.\h49.21.\h49.22.
\h49.23.\h49.24.\h49.25.\h49.26.\h49.27.\h49.28.\h49.29.\h49.30.
\h49.31.\h49.32.\h49.33.\h49.34.\h49.35.\h49.36.\h49.37.\h49.38.
\h49.39.\h49.40.\h49.41.\h49.42.\h49.43.\h49.44.\h49.45.\h49.46.
\h49.47.\h49.48.\h49.49.\h50.4.\h50.5.\h50.6.\h50.7.\h50.8.
\h50.9.\h50.10.\h50.11.\h50.12.\h50.13.\h50.14.\h50.15.\h50.16.
\h50.17.\h50.18.\h50.19.\h50.20.\h50.21.\h50.22.\h50.23.\h50.24.
\h50.25.\h50.26.\h50.27.\h50.28.\h50.29.\h50.30.\h50.31.\h50.32.
\h50.33.\h50.34.\h50.35.\h50.36.\h50.37.\h50.38.\h50.39.\h50.40.
\h50.41.\h50.42.\h50.43.\h50.44.\h50.45.\h50.46.\h50.47.\h50.48.
\h50.49.\h50.50.\h51.3.\h51.5.\h51.6.\h51.7.\h51.8.\h51.9.
\h51.10.\h51.11.\h51.12.\h51.13.\h51.14.\h51.15.\h51.16.\h51.17.
\h51.18.\h51.19.\h51.20.\h51.21.\h51.22.\h51.23.\h51.24.\h51.25.
\h51.26.\h51.27.\h51.28.\h51.29.\h51.30.\h51.31.\h51.32.\h51.33.
\h51.34.\h51.35.\h51.36.\h51.37.\h51.38.\h51.39.\h51.40.\h51.41.
\h51.42.\h51.43.\h51.44.\h51.45.\h51.46.\h51.47.\h51.48.\h51.49.
\h51.50.\h51.51.\h52.4.\h52.5.\h52.6.\h52.7.\h52.8.\h52.9.
\h52.10.\h52.11.\h52.12.\h52.13.\h52.14.\h52.15.\h52.16.\h52.17.
\h52.18.\h52.19.\h52.20.\h52.21.\h52.22.\h52.23.\h52.24.\h52.25.
\h52.26.\h52.27.\h52.28.\h52.29.\h52.30.\h52.31.\h52.32.\h52.33.
\h52.34.\h52.35.\h52.36.\h52.37.\h52.38.\h52.39.\h52.40.\h52.41.
\h52.42.\h52.43.\h52.44.\h52.45.\h52.46.\h52.47.\h52.48.\h52.49.
\h52.50.\h52.51.\h52.52.\h53.5.\h53.6.\h53.7.\h53.8.\h53.9.
\h53.10.\h53.11.\h53.12.\h53.13.\h53.14.\h53.15.\h53.16.\h53.17.
\h53.18.\h53.19.\h53.20.\h53.21.\h53.22.\h53.23.\h53.24.\h53.25.
\h53.26.\h53.27.\h53.28.\h53.29.\h53.30.\h53.31.\h53.32.\h53.33.
\h53.34.\h53.35.\h53.36.\h53.37.\h53.38.\h53.39.\h53.40.\h53.41.
\h53.42.\h53.43.\h53.44.\h53.45.\h53.46.\h53.47.\h53.48.\h53.49.
\h53.50.\h53.51.\h53.52.\h53.53.\h54.4.\h54.5.\h54.6.\h54.7.
\h54.8.\h54.9.\h54.10.\h54.11.\h54.12.\h54.13.\h54.14.\h54.15.
\h54.16.\h54.17.\h54.18.\h54.19.\h54.20.\h54.21.\h54.22.\h54.23.
\h54.24.\h54.25.\h54.26.\h54.27.\h54.28.\h54.29.\h54.30.\h54.31.
\h54.32.\h54.33.\h54.34.\h54.35.\h54.36.\h54.37.\h54.38.\h54.39.
\h54.40.\h54.41.\h54.42.\h54.43.\h54.44.\h54.45.\h54.46.\h54.47.
\h54.48.\h54.49.\h54.50.\h54.51.\h54.52.\h54.53.\h54.54.\h55.4.
\h55.5.\h55.6.\h55.7.\h55.8.\h55.9.\h55.10.\h55.11.\h55.12.
\h55.13.\h55.14.\h55.15.\h55.16.\h55.17.\h55.18.\h55.19.\h55.20.
\h55.21.\h55.22.\h55.23.\h55.24.\h55.25.\h55.26.\h55.27.\h55.28.
\h55.29.\h55.30.\h55.31.\h55.32.\h55.33.\h55.34.\h55.35.\h55.36.
\h55.37.\h55.38.\h55.39.\h55.40.\h55.41.\h55.42.\h55.43.\h55.44.
\h55.45.\h55.46.\h55.47.\h55.48.\h55.49.\h55.50.\h55.51.\h55.52.
\h55.53.\h55.54.\h55.55.\h56.4.\h56.5.\h56.6.\h56.7.\h56.8.
\h56.9.\h56.10.\h56.11.\h56.12.\h56.13.\h56.14.\h56.15.\h56.16.
\h56.17.\h56.18.\h56.19.\h56.20.\h56.21.\h56.22.\h56.23.\h56.24.
\h56.25.\h56.26.\h56.27.\h56.28.\h56.29.\h56.30.\h56.31.\h56.32.
\h56.33.\h56.34.\h56.35.\h56.36.\h56.37.\h56.38.\h56.39.\h56.40.
\h56.41.\h56.42.\h56.43.\h56.44.\h56.45.\h56.46.\h56.47.\h56.48.
\h56.49.\h56.50.\h56.51.\h56.52.\h56.53.\h56.54.\h56.55.\h56.56.
\h57.3.\h57.4.\h57.5.\h57.6.\h57.7.\h57.8.\h57.9.\h57.10.
\h57.11.\h57.12.\h57.13.\h57.14.\h57.15.\h57.16.\h57.17.\h57.18.
\h57.19.\h57.20.\h57.21.\h57.22.\h57.23.\h57.24.\h57.25.\h57.26.
\h57.27.\h57.28.\h57.29.\h57.30.\h57.31.\h57.32.\h57.33.\h57.34.
\h57.35.\h57.36.\h57.37.\h57.38.\h57.39.\h57.40.\h57.41.\h57.42.
\h57.43.\h57.44.\h57.45.\h57.46.\h57.47.\h57.48.\h57.49.\h57.50.
\h57.51.\h57.52.\h57.53.\h57.54.\h57.55.\h57.56.\h57.57.\h58.4.
\h58.5.\h58.6.\h58.7.\h58.8.\h58.9.\h58.10.\h58.11.\h58.12.
\h58.13.\h58.14.\h58.15.\h58.16.\h58.17.\h58.18.\h58.19.\h58.20.
\h58.21.\h58.22.\h58.23.\h58.24.\h58.25.\h58.26.\h58.27.\h58.28.
\h58.29.\h58.30.\h58.31.\h58.32.\h58.33.\h58.34.\h58.35.\h58.36.
\h58.37.\h58.38.\h58.39.\h58.40.\h58.41.\h58.42.\h58.43.\h58.44.
\h58.45.\h58.46.\h58.47.\h58.48.\h58.49.\h58.50.\h58.51.\h58.52.
\h58.53.\h58.54.\h58.55.\h58.56.\h58.57.\h58.58.\h59.3.\h59.4.
\h59.5.\h59.6.\h59.7.\h59.8.\h59.9.\h59.10.\h59.11.\h59.12.
\h59.13.\h59.14.\h59.15.\h59.16.\h59.17.\h59.18.\h59.19.\h59.20.
\h59.21.\h59.22.\h59.23.\h59.24.\h59.25.\h59.26.\h59.27.\h59.28.
\h59.29.\h59.30.\h59.31.\h59.32.\h59.33.\h59.34.\h59.35.\h59.36.
\h59.37.\h59.38.\h59.39.\h59.40.\h59.41.\h59.42.\h59.43.\h59.44.
\h59.45.\h59.46.\h59.47.\h59.48.\h59.49.\h59.50.\h59.51.\h59.52.
\h59.53.\h59.54.\h59.55.\h59.56.\h59.57.\h59.58.\h59.59.\h60.4.
\h60.5.\h60.6.\h60.7.\h60.8.\h60.9.\h60.10.\h60.11.\h60.12.
\h60.13.\h60.14.\h60.15.\h60.16.\h60.17.\h60.18.\h60.19.\h60.20.
\h60.21.\h60.22.\h60.23.\h60.24.\h60.25.\h60.26.\h60.27.\h60.28.
\h60.29.\h60.30.\h60.31.\h60.32.\h60.33.\h60.34.\h60.35.\h60.36.
\h60.37.\h60.38.\h60.39.\h60.40.\h60.41.\h60.42.\h60.43.\h60.44.
\h60.45.\h60.46.\h60.47.\h60.48.\h60.49.\h60.50.\h60.51.\h60.52.
\h60.53.\h60.54.\h60.55.\h60.56.\h60.57.\h60.58.\h60.59.\h60.60.
\h61.4.\h61.5.\h61.6.\h61.7.\h61.8.\h61.9.\h61.10.\h61.11.
\h61.12.\h61.13.\h61.14.\h61.15.\h61.16.\h61.17.\h61.18.\h61.19.
\h61.20.\h61.21.\h61.22.\h61.23.\h61.24.\h61.25.\h61.26.\h61.27.
\h61.28.\h61.29.\h61.30.\h61.31.\h61.32.\h61.33.\h61.34.\h61.35.
\h61.36.\h61.37.\h61.38.\h61.39.\h61.40.\h61.41.\h61.42.\h61.43.
\h61.44.\h61.45.\h61.46.\h61.47.\h61.48.\h61.49.\h61.50.\h61.51.
\h61.52.\h61.53.\h61.54.\h61.55.\h61.56.\h61.57.\h61.58.\h61.59.
\h61.60.\h61.61.\h62.4.\h62.5.\h62.6.\h62.7.\h62.8.\h62.9.
\h62.10.\h62.11.\h62.12.\h62.13.\h62.14.\h62.15.\h62.16.\h62.17.
\h62.18.\h62.19.\h62.20.\h62.21.\h62.22.\h62.23.\h62.24.\h62.25.
\h62.26.\h62.27.\h62.28.\h62.29.\h62.30.\h62.31.\h62.32.\h62.33.
\h62.34.\h62.35.\h62.36.\h62.37.\h62.38.\h62.39.\h62.40.\h62.41.
\h62.42.\h62.43.\h62.44.\h62.45.\h62.46.\h62.47.\h62.48.\h62.49.
\h62.50.\h62.51.\h62.52.\h62.53.\h62.54.\h62.55.\h62.56.\h62.57.
\h62.58.\h62.59.\h62.60.\h62.61.\h62.62.\h63.3.\h63.4.\h63.5.
\h63.6.\h63.7.\h63.8.\h63.9.\h63.10.\h63.11.\h63.12.\h63.13.
\h63.14.\h63.15.\h63.16.\h63.17.\h63.18.\h63.19.\h63.20.\h63.21.
\h63.22.\h63.23.\h63.24.\h63.25.\h63.26.\h63.27.\h63.28.\h63.29.
\h63.30.\h63.31.\h63.32.\h63.33.\h63.34.\h63.35.\h63.36.\h63.37.
\h63.38.\h63.39.\h63.40.\h63.41.\h63.42.\h63.43.\h63.44.\h63.45.
\h63.46.\h63.47.\h63.48.\h63.49.\h63.50.\h63.51.\h63.52.\h63.53.
\h63.54.\h63.55.\h63.56.\h63.57.\h63.58.\h63.59.\h63.60.\h63.61.
\h63.62.\h63.63.\h64.4.\h64.5.\h64.6.\h64.7.\h64.8.\h64.9.
\h64.10.\h64.11.\h64.12.\h64.13.\h64.14.\h64.15.\h64.16.\h64.17.
\h64.18.\h64.19.\h64.20.\h64.21.\h64.22.\h64.23.\h64.24.\h64.25.
\h64.26.\h64.27.\h64.28.\h64.29.\h64.30.\h64.31.\h64.32.\h64.33.
\h64.34.\h64.35.\h64.36.\h64.37.\h64.38.\h64.39.\h64.40.\h64.41.
\h64.42.\h64.43.\h64.44.\h64.45.\h64.46.\h64.47.\h64.48.\h64.49.
\h64.50.\h64.51.\h64.52.\h64.53.\h64.54.\h64.55.\h64.56.\h64.57.
\h64.58.\h64.59.\h64.60.\h64.61.\h64.62.\h64.63.\h64.64.\h65.3.
\h65.4.\h65.5.\h65.6.\h65.7.\h65.8.\h65.9.\h65.10.\h65.11.
\h65.12.\h65.13.\h65.14.\h65.15.\h65.16.\h65.17.\h65.18.\h65.19.
\h65.20.\h65.21.\h65.22.\h65.23.\h65.24.\h65.25.\h65.26.\h65.27.
\h65.28.\h65.29.\h65.30.\h65.31.\h65.32.\h65.33.\h65.34.\h65.35.
\h65.36.\h65.37.\h65.38.\h65.39.\h65.40.\h65.41.\h65.42.\h65.43.
\h65.44.\h65.45.\h65.46.\h65.47.\h65.48.\h65.49.\h65.50.\h65.51.
\h65.52.\h65.53.\h65.54.\h65.55.\h65.56.\h65.57.\h65.58.\h65.59.
\h65.60.\h65.61.\h65.62.\h65.63.\h65.64.\h65.65.\h66.3.\h66.4.
\h66.5.\h66.6.\h66.7.\h66.8.\h66.9.\h66.10.\h66.11.\h66.12.
\h66.13.\h66.14.\h66.15.\h66.16.\h66.17.\h66.18.\h66.19.\h66.20.
\h66.21.\h66.22.\h66.23.\h66.24.\h66.25.\h66.26.\h66.27.\h66.28.
\h66.29.\h66.30.\h66.31.\h66.32.\h66.33.\h66.34.\h66.35.\h66.36.
\h66.37.\h66.38.\h66.39.\h66.40.\h66.41.\h66.42.\h66.43.\h66.44.
\h66.45.\h66.46.\h66.47.\h66.48.\h66.49.\h66.50.\h66.51.\h66.52.
\h66.53.\h66.54.\h66.55.\h66.56.\h66.57.\h66.58.\h66.59.\h66.60.
\h66.61.\h66.62.\h66.63.\h66.64.\h66.65.\h66.66.\h67.3.\h67.4.
\h67.5.\h67.6.\h67.7.\h67.8.\h67.9.\h67.10.\h67.11.\h67.12.
\h67.13.\h67.14.\h67.15.\h67.16.\h67.17.\h67.18.\h67.19.\h67.20.
\h67.21.\h67.22.\h67.23.\h67.24.\h67.25.\h67.26.\h67.27.\h67.28.
\h67.29.\h67.30.\h67.31.\h67.32.\h67.33.\h67.34.\h67.35.\h67.36.
\h67.37.\h67.38.\h67.39.\h67.40.\h67.41.\h67.42.\h67.43.\h67.44.
\h67.45.\h67.46.\h67.47.\h67.48.\h67.49.\h67.50.\h67.51.\h67.52.
\h67.53.\h67.54.\h67.55.\h67.56.\h67.57.\h67.58.\h67.59.\h67.60.
\h67.61.\h67.62.\h67.63.\h67.64.\h67.65.\h67.66.\h67.67.\h68.4.
\h68.5.\h68.6.\h68.7.\h68.8.\h68.9.\h68.10.\h68.11.\h68.12.
\h68.13.\h68.14.\h68.15.\h68.16.\h68.17.\h68.18.\h68.19.\h68.20.
\h68.21.\h68.22.\h68.23.\h68.24.\h68.25.\h68.26.\h68.27.\h68.28.
\h68.29.\h68.30.\h68.31.\h68.32.\h68.33.\h68.34.\h68.35.\h68.36.
\h68.37.\h68.38.\h68.39.\h68.40.\h68.41.\h68.42.\h68.43.\h68.44.
\h68.45.\h68.46.\h68.47.\h68.48.\h68.49.\h68.50.\h68.51.\h68.52.
\h68.53.\h68.54.\h68.55.\h68.56.\h68.57.\h68.58.\h68.59.\h68.60.
\h68.61.\h68.62.\h68.63.\h68.64.\h68.65.\h68.66.\h68.67.\h68.68.
\h69.3.\h69.4.\h69.5.\h69.6.\h69.7.\h69.8.\h69.9.\h69.10.
\h69.11.\h69.12.\h69.13.\h69.14.\h69.15.\h69.16.\h69.17.\h69.18.
\h69.19.\h69.20.\h69.21.\h69.22.\h69.23.\h69.24.\h69.25.\h69.26.
\h69.27.\h69.28.\h69.29.\h69.30.\h69.31.\h69.32.\h69.33.\h69.34.
\h69.35.\h69.36.\h69.37.\h69.38.\h69.39.\h69.40.\h69.41.\h69.42.
\h69.43.\h69.44.\h69.45.\h69.46.\h69.47.\h69.48.\h69.49.\h69.50.
\h69.51.\h69.52.\h69.53.\h69.54.\h69.55.\h69.56.\h69.57.\h69.58.
\h69.59.\h69.60.\h69.61.\h69.62.\h69.63.\h69.64.\h69.65.\h69.66.
\h69.67.\h69.68.\h69.69.\h70.4.\h70.5.\h70.6.\h70.7.\h70.8.
\h70.9.\h70.10.\h70.11.\h70.12.\h70.13.\h70.14.\h70.15.\h70.16.
\h70.17.\h70.18.\h70.19.\h70.20.\h70.21.\h70.22.\h70.23.\h70.24.
\h70.25.\h70.26.\h70.27.\h70.28.\h70.29.\h70.30.\h70.31.\h70.32.
\h70.33.\h70.34.\h70.35.\h70.36.\h70.37.\h70.38.\h70.39.\h70.40.
\h70.41.\h70.42.\h70.43.\h70.44.\h70.45.\h70.46.\h70.47.\h70.48.
\h70.49.\h70.50.\h70.51.\h70.52.\h70.53.\h70.54.\h70.55.\h70.56.
\h70.57.\h70.58.\h70.59.\h70.60.\h70.61.\h70.62.\h70.63.\h70.64.
\h70.65.\h70.66.\h70.67.\h70.68.\h70.69.\h70.70.\h71.3.\h71.4.
\h71.5.\h71.6.\h71.7.\h71.8.\h71.9.\h71.10.\h71.11.\h71.12.
\h71.13.\h71.14.\h71.15.\h71.16.\h71.17.\h71.18.\h71.19.\h71.20.
\h71.21.\h71.22.\h71.23.\h71.24.\h71.25.\h71.26.\h71.27.\h71.28.
\h71.29.\h71.30.\h71.31.\h71.32.\h71.33.\h71.34.\h71.35.\h71.36.
\h71.37.\h71.38.\h71.39.\h71.40.\h71.41.\h71.42.\h71.43.\h71.44.
\h71.45.\h71.46.\h71.47.\h71.48.\h71.49.\h71.50.\h71.51.\h71.52.
\h71.53.\h71.54.\h71.55.\h71.56.\h71.57.\h71.58.\h71.59.\h71.60.
\h71.61.\h71.62.\h71.63.\h71.64.\h71.65.\h71.66.\h71.67.\h71.68.
\h71.69.\h71.70.\h71.71.\h72.3.\h72.4.\h72.5.\h72.6.\h72.7.
\h72.8.\h72.9.\h72.10.\h72.11.\h72.12.\h72.13.\h72.14.\h72.15.
\h72.16.\h72.17.\h72.18.\h72.19.\h72.20.\h72.21.\h72.22.\h72.23.
\h72.24.\h72.25.\h72.26.\h72.27.\h72.28.\h72.29.\h72.30.\h72.31.
\h72.32.\h72.33.\h72.34.\h72.35.\h72.36.\h72.37.\h72.38.\h72.39.
\h72.40.\h72.41.\h72.42.\h72.43.\h72.44.\h72.45.\h72.46.\h72.47.
\h72.48.\h72.49.\h72.50.\h72.51.\h72.52.\h72.53.\h72.54.\h72.55.
\h72.56.\h72.57.\h72.58.\h72.59.\h72.60.\h72.61.\h72.62.\h72.63.
\h72.64.\h72.65.\h72.66.\h72.67.\h72.68.\h72.69.\h72.70.\h72.71.
\h72.72.\h73.3.\h73.4.\h73.5.\h73.6.\h73.7.\h73.8.\h73.9.
\h73.10.\h73.11.\h73.12.\h73.13.\h73.14.\h73.15.\h73.16.\h73.17.
\h73.18.\h73.19.\h73.20.\h73.21.\h73.22.\h73.23.\h73.24.\h73.25.
\h73.26.\h73.27.\h73.28.\h73.29.\h73.30.\h73.31.\h73.32.\h73.33.
\h73.34.\h73.35.\h73.36.\h73.37.\h73.38.\h73.39.\h73.40.\h73.41.
\h73.42.\h73.43.\h73.44.\h73.45.\h73.46.\h73.47.\h73.48.\h73.49.
\h73.50.\h73.51.\h73.52.\h73.53.\h73.54.\h73.55.\h73.56.\h73.57.
\h73.58.\h73.59.\h73.60.\h73.61.\h73.62.\h73.63.\h73.64.\h73.65.
\h73.66.\h73.67.\h73.68.\h73.69.\h73.70.\h73.71.\h73.72.\h73.73.
\h74.2.\h74.4.\h74.5.\h74.6.\h74.7.\h74.8.\h74.9.\h74.10.
\h74.11.\h74.12.\h74.13.\h74.14.\h74.15.\h74.16.\h74.17.\h74.18.
\h74.19.\h74.20.\h74.21.\h74.22.\h74.23.\h74.24.\h74.25.\h74.26.
\h74.27.\h74.28.\h74.29.\h74.30.\h74.31.\h74.32.\h74.33.\h74.34.
\h74.35.\h74.36.\h74.37.\h74.38.\h74.39.\h74.40.\h74.41.\h74.42.
\h74.43.\h74.44.\h74.45.\h74.46.\h74.47.\h74.48.\h74.49.\h74.50.
\h74.51.\h74.52.\h74.53.\h74.54.\h74.55.\h74.56.\h74.57.\h74.58.
\h74.59.\h74.60.\h74.61.\h74.62.\h74.63.\h74.64.\h74.65.\h74.66.
\h74.67.\h74.68.\h74.69.\h74.70.\h74.71.\h74.72.\h74.73.\h74.74.
\h75.3.\h75.5.\h75.6.\h75.7.\h75.8.\h75.9.\h75.10.\h75.11.
\h75.12.\h75.13.\h75.14.\h75.15.\h75.16.\h75.17.\h75.18.\h75.19.
\h75.20.\h75.21.\h75.22.\h75.23.\h75.24.\h75.25.\h75.26.\h75.27.
\h75.28.\h75.29.\h75.30.\h75.31.\h75.32.\h75.33.\h75.34.\h75.35.
\h75.36.\h75.37.\h75.38.\h75.39.\h75.40.\h75.41.\h75.42.\h75.43.
\h75.44.\h75.45.\h75.46.\h75.47.\h75.48.\h75.49.\h75.50.\h75.51.
\h75.52.\h75.53.\h75.54.\h75.55.\h75.56.\h75.57.\h75.58.\h75.59.
\h75.60.\h75.61.\h75.62.\h75.63.\h75.64.\h75.65.\h75.66.\h75.67.
\h75.68.\h75.69.\h75.70.\h75.71.\h75.72.\h75.73.\h75.74.\h75.75.
\h76.3.\h76.4.\h76.5.\h76.6.\h76.7.\h76.8.\h76.9.\h76.10.
\h76.11.\h76.12.\h76.13.\h76.14.\h76.15.\h76.16.\h76.17.\h76.18.
\h76.19.\h76.20.\h76.21.\h76.22.\h76.23.\h76.24.\h76.25.\h76.26.
\h76.27.\h76.28.\h76.29.\h76.30.\h76.31.\h76.32.\h76.33.\h76.34.
\h76.35.\h76.36.\h76.37.\h76.38.\h76.39.\h76.40.\h76.41.\h76.42.
\h76.43.\h76.44.\h76.45.\h76.46.\h76.47.\h76.48.\h76.49.\h76.50.
\h76.51.\h76.52.\h76.53.\h76.54.\h76.55.\h76.56.\h76.57.\h76.58.
\h76.59.\h76.60.\h76.61.\h76.62.\h76.63.\h76.64.\h76.65.\h76.66.
\h76.67.\h76.68.\h76.69.\h76.70.\h76.71.\h76.72.\h76.73.\h76.74.
\h76.75.\h76.76.\h77.3.\h77.4.\h77.5.\h77.6.\h77.7.\h77.8.
\h77.9.\h77.10.\h77.11.\h77.12.\h77.13.\h77.14.\h77.15.\h77.16.
\h77.17.\h77.18.\h77.19.\h77.20.\h77.21.\h77.22.\h77.23.\h77.24.
\h77.25.\h77.26.\h77.27.\h77.28.\h77.29.\h77.30.\h77.31.\h77.32.
\h77.33.\h77.34.\h77.35.\h77.36.\h77.37.\h77.38.\h77.39.\h77.40.
\h77.41.\h77.42.\h77.43.\h77.44.\h77.45.\h77.46.\h77.47.\h77.48.
\h77.49.\h77.50.\h77.51.\h77.52.\h77.53.\h77.54.\h77.55.\h77.56.
\h77.57.\h77.58.\h77.59.\h77.60.\h77.61.\h77.62.\h77.63.\h77.64.
\h77.65.\h77.66.\h77.67.\h77.68.\h77.69.\h77.70.\h77.71.\h77.72.
\h77.73.\h77.74.\h77.75.\h77.76.\h77.77.\h78.3.\h78.4.\h78.5.
\h78.6.\h78.7.\h78.8.\h78.9.\h78.10.\h78.11.\h78.12.\h78.13.
\h78.14.\h78.15.\h78.16.\h78.17.\h78.18.\h78.19.\h78.20.\h78.21.
\h78.22.\h78.23.\h78.24.\h78.25.\h78.26.\h78.27.\h78.28.\h78.29.
\h78.30.\h78.31.\h78.32.\h78.33.\h78.34.\h78.35.\h78.36.\h78.37.
\h78.38.\h78.39.\h78.40.\h78.41.\h78.42.\h78.43.\h78.44.\h78.45.
\h78.46.\h78.47.\h78.48.\h78.49.\h78.50.\h78.51.\h78.52.\h78.53.
\h78.54.\h78.55.\h78.56.\h78.57.\h78.58.\h78.59.\h78.60.\h78.61.
\h78.62.\h78.63.\h78.64.\h78.65.\h78.66.\h78.67.\h78.68.\h78.69.
\h78.70.\h78.71.\h78.72.\h78.73.\h78.74.\h78.75.\h78.76.\h78.77.
\h78.78.\h79.3.\h79.4.\h79.5.\h79.6.\h79.7.\h79.8.\h79.9.
\h79.10.\h79.11.\h79.12.\h79.13.\h79.14.\h79.15.\h79.16.\h79.17.
\h79.18.\h79.19.\h79.20.\h79.21.\h79.22.\h79.23.\h79.24.\h79.25.
\h79.26.\h79.27.\h79.28.\h79.29.\h79.30.\h79.31.\h79.32.\h79.33.
\h79.34.\h79.35.\h79.36.\h79.37.\h79.38.\h79.39.\h79.40.\h79.41.
\h79.42.\h79.43.\h79.44.\h79.45.\h79.46.\h79.47.\h79.48.\h79.49.
\h79.50.\h79.51.\h79.52.\h79.53.\h79.54.\h79.55.\h79.56.\h79.57.
\h79.58.\h79.59.\h79.60.\h79.61.\h79.62.\h79.63.\h79.64.\h79.65.
\h79.66.\h79.67.\h79.68.\h79.69.\h79.70.\h79.71.\h79.72.\h79.73.
\h79.74.\h79.75.\h79.76.\h79.77.\h79.78.\h79.79.\h80.4.\h80.5.
\h80.6.\h80.7.\h80.8.\h80.9.\h80.10.\h80.11.\h80.12.\h80.13.
\h80.14.\h80.15.\h80.16.\h80.17.\h80.18.\h80.19.\h80.20.\h80.21.
\h80.22.\h80.23.\h80.24.\h80.25.\h80.26.\h80.27.\h80.28.\h80.29.
\h80.30.\h80.31.\h80.32.\h80.33.\h80.34.\h80.35.\h80.36.\h80.37.
\h80.38.\h80.39.\h80.40.\h80.41.\h80.42.\h80.43.\h80.44.\h80.45.
\h80.46.\h80.47.\h80.48.\h80.49.\h80.50.\h80.51.\h80.52.\h80.53.
\h80.54.\h80.55.\h80.56.\h80.57.\h80.58.\h80.59.\h80.60.\h80.61.
\h80.62.\h80.63.\h80.64.\h80.65.\h80.66.\h80.67.\h80.68.\h80.69.
\h80.70.\h80.71.\h80.72.\h80.73.\h80.74.\h80.75.\h80.76.\h80.77.
\h80.78.\h80.79.\h80.80.\h81.3.\h81.4.\h81.5.\h81.6.\h81.7.
\h81.8.\h81.9.\h81.10.\h81.11.\h81.12.\h81.13.\h81.14.\h81.15.
\h81.16.\h81.17.\h81.18.\h81.19.\h81.20.\h81.21.\h81.22.\h81.23.
\h81.24.\h81.25.\h81.26.\h81.27.\h81.28.\h81.29.\h81.30.\h81.31.
\h81.32.\h81.33.\h81.34.\h81.35.\h81.36.\h81.37.\h81.38.\h81.39.
\h81.40.\h81.41.\h81.42.\h81.43.\h81.44.\h81.45.\h81.46.\h81.47.
\h81.48.\h81.49.\h81.50.\h81.51.\h81.52.\h81.53.\h81.54.\h81.55.
\h81.56.\h81.57.\h81.58.\h81.59.\h81.60.\h81.61.\h81.62.\h81.63.
\h81.64.\h81.65.\h81.66.\h81.67.\h81.68.\h81.69.\h81.70.\h81.71.
\h81.72.\h81.73.\h81.74.\h81.75.\h81.76.\h81.77.\h81.78.\h81.79.
\h81.80.\h81.81.\h82.4.\h82.5.\h82.6.\h82.7.\h82.8.\h82.9.
\h82.10.\h82.11.\h82.12.\h82.13.\h82.14.\h82.15.\h82.16.\h82.17.
\h82.18.\h82.19.\h82.20.\h82.21.\h82.22.\h82.23.\h82.24.\h82.25.
\h82.26.\h82.27.\h82.28.\h82.29.\h82.30.\h82.31.\h82.32.\h82.33.
\h82.34.\h82.35.\h82.36.\h82.37.\h82.38.\h82.39.\h82.40.\h82.41.
\h82.42.\h82.43.\h82.44.\h82.45.\h82.46.\h82.47.\h82.48.\h82.49.
\h82.50.\h82.51.\h82.52.\h82.53.\h82.54.\h82.55.\h82.56.\h82.57.
\h82.58.\h82.59.\h82.60.\h82.61.\h82.62.\h82.63.\h82.64.\h82.65.
\h82.66.\h82.67.\h82.68.\h82.69.\h82.70.\h82.71.\h82.72.\h82.73.
\h82.74.\h82.75.\h82.76.\h82.77.\h82.78.\h82.79.\h82.80.\h82.81.
\h82.82.\h83.2.\h83.3.\h83.5.\h83.6.\h83.7.\h83.8.\h83.9.
\h83.10.\h83.11.\h83.12.\h83.13.\h83.14.\h83.15.\h83.16.\h83.17.
\h83.18.\h83.19.\h83.20.\h83.21.\h83.22.\h83.23.\h83.24.\h83.25.
\h83.26.\h83.27.\h83.28.\h83.29.\h83.30.\h83.31.\h83.32.\h83.33.
\h83.34.\h83.35.\h83.36.\h83.37.\h83.38.\h83.39.\h83.40.\h83.41.
\h83.42.\h83.43.\h83.44.\h83.45.\h83.46.\h83.47.\h83.48.\h83.49.
\h83.50.\h83.51.\h83.52.\h83.53.\h83.54.\h83.55.\h83.56.\h83.57.
\h83.58.\h83.59.\h83.60.\h83.61.\h83.62.\h83.63.\h83.64.\h83.65.
\h83.66.\h83.67.\h83.68.\h83.69.\h83.70.\h83.71.\h83.72.\h83.73.
\h83.74.\h83.75.\h83.76.\h83.77.\h83.78.\h83.79.\h83.80.\h83.81.
\h83.82.\h83.83.\h84.2.\h84.3.\h84.4.\h84.5.\h84.6.\h84.7.
\h84.8.\h84.9.\h84.10.\h84.11.\h84.12.\h84.13.\h84.14.\h84.15.
\h84.16.\h84.17.\h84.18.\h84.19.\h84.20.\h84.21.\h84.22.\h84.23.
\h84.24.\h84.25.\h84.26.\h84.27.\h84.28.\h84.29.\h84.30.\h84.31.
\h84.32.\h84.33.\h84.34.\h84.35.\h84.36.\h84.37.\h84.38.\h84.39.
\h84.40.\h84.41.\h84.42.\h84.43.\h84.44.\h84.45.\h84.46.\h84.47.
\h84.48.\h84.49.\h84.50.\h84.51.\h84.52.\h84.53.\h84.54.\h84.55.
\h84.56.\h84.57.\h84.58.\h84.59.\h84.60.\h84.61.\h84.62.\h84.63.
\h84.64.\h84.65.\h84.66.\h84.67.\h84.68.\h84.69.\h84.70.\h84.71.
\h84.72.\h84.73.\h84.74.\h84.75.\h84.76.\h84.77.\h84.78.\h84.79.
\h84.80.\h84.81.\h84.82.\h84.83.\h84.84.\h85.3.\h85.4.\h85.5.
\h85.6.\h85.7.\h85.8.\h85.9.\h85.10.\h85.11.\h85.12.\h85.13.
\h85.14.\h85.15.\h85.16.\h85.17.\h85.18.\h85.19.\h85.20.\h85.21.
\h85.22.\h85.23.\h85.24.\h85.25.\h85.26.\h85.27.\h85.28.\h85.29.
\h85.30.\h85.31.\h85.32.\h85.33.\h85.34.\h85.35.\h85.36.\h85.37.
\h85.38.\h85.39.\h85.40.\h85.41.\h85.42.\h85.43.\h85.44.\h85.45.
\h85.46.\h85.47.\h85.48.\h85.49.\h85.50.\h85.51.\h85.52.\h85.53.
\h85.54.\h85.55.\h85.56.\h85.57.\h85.58.\h85.59.\h85.60.\h85.61.
\h85.62.\h85.63.\h85.64.\h85.65.\h85.66.\h85.67.\h85.68.\h85.69.
\h85.70.\h85.71.\h85.72.\h85.73.\h85.74.\h85.75.\h85.76.\h85.77.
\h85.78.\h85.79.\h85.80.\h85.81.\h85.82.\h85.83.\h85.84.\h85.85.
\h86.2.\h86.4.\h86.5.\h86.6.\h86.7.\h86.8.\h86.9.\h86.10.
\h86.11.\h86.12.\h86.13.\h86.14.\h86.15.\h86.16.\h86.17.\h86.18.
\h86.19.\h86.20.\h86.21.\h86.22.\h86.23.\h86.24.\h86.25.\h86.26.
\h86.27.\h86.28.\h86.29.\h86.30.\h86.31.\h86.32.\h86.33.\h86.34.
\h86.35.\h86.36.\h86.37.\h86.38.\h86.39.\h86.40.\h86.41.\h86.42.
\h86.43.\h86.44.\h86.45.\h86.46.\h86.47.\h86.48.\h86.49.\h86.50.
\h86.51.\h86.52.\h86.53.\h86.54.\h86.55.\h86.56.\h86.57.\h86.58.
\h86.59.\h86.60.\h86.61.\h86.62.\h86.63.\h86.64.\h86.65.\h86.66.
\h86.67.\h86.68.\h86.69.\h86.70.\h86.71.\h86.72.\h86.73.\h86.74.
\h86.75.\h86.76.\h86.77.\h86.78.\h86.79.\h86.80.\h86.81.\h86.82.
\h86.83.\h86.84.\h86.85.\h86.86.\h87.3.\h87.5.\h87.6.\h87.7.
\h87.8.\h87.9.\h87.10.\h87.11.\h87.12.\h87.13.\h87.14.\h87.15.
\h87.16.\h87.17.\h87.18.\h87.19.\h87.20.\h87.21.\h87.22.\h87.23.
\h87.24.\h87.25.\h87.26.\h87.27.\h87.28.\h87.29.\h87.30.\h87.31.
\h87.32.\h87.33.\h87.34.\h87.35.\h87.36.\h87.37.\h87.38.\h87.39.
\h87.40.\h87.41.\h87.42.\h87.43.\h87.44.\h87.45.\h87.46.\h87.47.
\h87.48.\h87.49.\h87.50.\h87.51.\h87.52.\h87.53.\h87.54.\h87.55.
\h87.56.\h87.57.\h87.58.\h87.59.\h87.60.\h87.61.\h87.62.\h87.63.
\h87.64.\h87.65.\h87.66.\h87.67.\h87.68.\h87.69.\h87.70.\h87.71.
\h87.72.\h87.73.\h87.74.\h87.75.\h87.76.\h87.77.\h87.78.\h87.79.
\h87.80.\h87.81.\h87.82.\h87.83.\h87.84.\h87.85.\h87.86.\h87.87.
\h88.4.\h88.5.\h88.6.\h88.7.\h88.8.\h88.9.\h88.10.\h88.11.
\h88.12.\h88.13.\h88.14.\h88.15.\h88.16.\h88.17.\h88.18.\h88.19.
\h88.20.\h88.21.\h88.22.\h88.23.\h88.24.\h88.25.\h88.26.\h88.27.
\h88.28.\h88.29.\h88.30.\h88.31.\h88.32.\h88.33.\h88.34.\h88.35.
\h88.36.\h88.37.\h88.38.\h88.39.\h88.40.\h88.41.\h88.42.\h88.43.
\h88.44.\h88.45.\h88.46.\h88.47.\h88.48.\h88.49.\h88.50.\h88.51.
\h88.52.\h88.53.\h88.54.\h88.55.\h88.56.\h88.57.\h88.58.\h88.59.
\h88.60.\h88.61.\h88.62.\h88.63.\h88.64.\h88.65.\h88.66.\h88.67.
\h88.68.\h88.69.\h88.70.\h88.71.\h88.72.\h88.73.\h88.74.\h88.75.
\h88.76.\h88.77.\h88.78.\h88.79.\h88.80.\h88.81.\h88.82.\h88.83.
\h88.84.\h88.85.\h88.86.\h88.87.\h88.88.\h89.3.\h89.4.\h89.5.
\h89.6.\h89.7.\h89.8.\h89.9.\h89.10.\h89.11.\h89.12.\h89.13.
\h89.14.\h89.15.\h89.16.\h89.17.\h89.18.\h89.19.\h89.20.\h89.21.
\h89.22.\h89.23.\h89.24.\h89.25.\h89.26.\h89.27.\h89.28.\h89.29.
\h89.30.\h89.31.\h89.32.\h89.33.\h89.34.\h89.35.\h89.36.\h89.37.
\h89.38.\h89.39.\h89.40.\h89.41.\h89.42.\h89.43.\h89.44.\h89.45.
\h89.46.\h89.47.\h89.48.\h89.49.\h89.50.\h89.51.\h89.52.\h89.53.
\h89.54.\h89.55.\h89.56.\h89.57.\h89.58.\h89.59.\h89.60.\h89.61.
\h89.62.\h89.63.\h89.64.\h89.65.\h89.66.\h89.67.\h89.68.\h89.69.
\h89.70.\h89.71.\h89.72.\h89.73.\h89.74.\h89.75.\h89.76.\h89.77.
\h89.78.\h89.79.\h89.80.\h89.81.\h89.82.\h89.83.\h89.84.\h89.85.
\h89.86.\h89.87.\h89.89.\h90.2.\h90.4.\h90.5.\h90.6.\h90.7.
\h90.8.\h90.9.\h90.10.\h90.11.\h90.12.\h90.13.\h90.14.\h90.15.
\h90.16.\h90.17.\h90.18.\h90.19.\h90.20.\h90.21.\h90.22.\h90.23.
\h90.24.\h90.25.\h90.26.\h90.27.\h90.28.\h90.29.\h90.30.\h90.31.
\h90.32.\h90.33.\h90.34.\h90.35.\h90.36.\h90.37.\h90.38.\h90.39.
\h90.40.\h90.41.\h90.42.\h90.43.\h90.44.\h90.45.\h90.46.\h90.47.
\h90.48.\h90.49.\h90.50.\h90.51.\h90.52.\h90.53.\h90.54.\h90.55.
\h90.56.\h90.57.\h90.58.\h90.59.\h90.60.\h90.61.\h90.62.\h90.63.
\h90.64.\h90.65.\h90.66.\h90.67.\h90.68.\h90.69.\h90.70.\h90.71.
\h90.72.\h90.73.\h90.74.\h90.75.\h90.76.\h90.77.\h90.78.\h90.79.
\h90.80.\h90.81.\h90.82.\h90.83.\h90.84.\h90.85.\h90.86.\h90.87.
\h90.88.\h90.89.\h90.90.\h91.3.\h91.4.\h91.5.\h91.6.\h91.7.
\h91.8.\h91.9.\h91.10.\h91.11.\h91.12.\h91.13.\h91.14.\h91.15.
\h91.16.\h91.17.\h91.18.\h91.19.\h91.20.\h91.21.\h91.22.\h91.23.
\h91.24.\h91.25.\h91.26.\h91.27.\h91.28.\h91.29.\h91.30.\h91.31.
\h91.32.\h91.33.\h91.34.\h91.35.\h91.36.\h91.37.\h91.38.\h91.39.
\h91.40.\h91.41.\h91.42.\h91.43.\h91.44.\h91.45.\h91.46.\h91.47.
\h91.48.\h91.49.\h91.50.\h91.51.\h91.52.\h91.53.\h91.54.\h91.55.
\h91.56.\h91.57.\h91.58.\h91.59.\h91.60.\h91.61.\h91.62.\h91.63.
\h91.64.\h91.65.\h91.66.\h91.67.\h91.68.\h91.69.\h91.70.\h91.71.
\h91.72.\h91.73.\h91.74.\h91.75.\h91.76.\h91.77.\h91.78.\h91.79.
\h91.80.\h91.81.\h91.82.\h91.83.\h91.84.\h91.85.\h91.86.\h91.87.
\h91.88.\h91.89.\h91.91.\h92.2.\h92.4.\h92.5.\h92.6.\h92.7.
\h92.8.\h92.9.\h92.10.\h92.11.\h92.12.\h92.13.\h92.14.\h92.15.
\h92.16.\h92.17.\h92.18.\h92.19.\h92.20.\h92.21.\h92.22.\h92.23.
\h92.24.\h92.25.\h92.26.\h92.27.\h92.28.\h92.29.\h92.30.\h92.31.
\h92.32.\h92.33.\h92.34.\h92.35.\h92.36.\h92.37.\h92.38.\h92.39.
\h92.40.\h92.41.\h92.42.\h92.43.\h92.44.\h92.45.\h92.46.\h92.47.
\h92.48.\h92.49.\h92.50.\h92.51.\h92.52.\h92.53.\h92.54.\h92.55.
\h92.56.\h92.57.\h92.58.\h92.59.\h92.60.\h92.61.\h92.62.\h92.63.
\h92.64.\h92.65.\h92.66.\h92.67.\h92.68.\h92.69.\h92.70.\h92.71.
\h92.72.\h92.73.\h92.74.\h92.75.\h92.76.\h92.77.\h92.78.\h92.79.
\h92.80.\h92.81.\h92.82.\h92.83.\h92.84.\h92.85.\h92.86.\h92.88.
\h92.89.\h92.90.\h92.91.\h92.92.\h93.3.\h93.4.\h93.5.\h93.6.
\h93.7.\h93.8.\h93.9.\h93.10.\h93.11.\h93.12.\h93.13.\h93.14.
\h93.15.\h93.16.\h93.17.\h93.18.\h93.19.\h93.20.\h93.21.\h93.22.
\h93.23.\h93.24.\h93.25.\h93.26.\h93.27.\h93.28.\h93.29.\h93.30.
\h93.31.\h93.32.\h93.33.\h93.34.\h93.35.\h93.36.\h93.37.\h93.38.
\h93.39.\h93.40.\h93.41.\h93.42.\h93.43.\h93.44.\h93.45.\h93.46.
\h93.47.\h93.48.\h93.49.\h93.50.\h93.51.\h93.52.\h93.53.\h93.54.
\h93.55.\h93.56.\h93.57.\h93.58.\h93.59.\h93.60.\h93.61.\h93.62.
\h93.63.\h93.64.\h93.65.\h93.66.\h93.67.\h93.68.\h93.69.\h93.70.
\h93.71.\h93.72.\h93.73.\h93.74.\h93.75.\h93.76.\h93.77.\h93.78.
\h93.79.\h93.80.\h93.81.\h93.82.\h93.83.\h93.84.\h93.85.\h93.86.
\h93.87.\h93.88.\h93.89.\h93.90.\h93.91.\h93.93.\h94.4.\h94.5.
\h94.6.\h94.7.\h94.8.\h94.9.\h94.10.\h94.11.\h94.12.\h94.13.
\h94.14.\h94.15.\h94.16.\h94.17.\h94.18.\h94.19.\h94.20.\h94.21.
\h94.22.\h94.23.\h94.24.\h94.25.\h94.26.\h94.27.\h94.28.\h94.29.
\h94.30.\h94.31.\h94.32.\h94.33.\h94.34.\h94.35.\h94.36.\h94.37.
\h94.38.\h94.39.\h94.40.\h94.41.\h94.42.\h94.43.\h94.44.\h94.45.
\h94.46.\h94.47.\h94.48.\h94.49.\h94.50.\h94.51.\h94.52.\h94.53.
\h94.54.\h94.55.\h94.56.\h94.57.\h94.58.\h94.59.\h94.60.\h94.61.
\h94.62.\h94.63.\h94.64.\h94.65.\h94.66.\h94.67.\h94.68.\h94.69.
\h94.70.\h94.71.\h94.72.\h94.73.\h94.74.\h94.75.\h94.76.\h94.77.
\h94.78.\h94.79.\h94.80.\h94.81.\h94.82.\h94.84.\h94.85.\h94.86.
\h94.87.\h94.88.\h94.89.\h94.90.\h94.91.\h94.92.\h94.93.\h94.94.
\h95.2.\h95.3.\h95.5.\h95.6.\h95.7.\h95.8.\h95.9.\h95.10.
\h95.11.\h95.12.\h95.13.\h95.14.\h95.15.\h95.16.\h95.17.\h95.18.
\h95.19.\h95.20.\h95.21.\h95.22.\h95.23.\h95.24.\h95.25.\h95.26.
\h95.27.\h95.28.\h95.29.\h95.30.\h95.31.\h95.32.\h95.33.\h95.34.
\h95.35.\h95.36.\h95.37.\h95.38.\h95.39.\h95.40.\h95.41.\h95.42.
\h95.43.\h95.44.\h95.45.\h95.46.\h95.47.\h95.48.\h95.49.\h95.50.
\h95.51.\h95.52.\h95.53.\h95.54.\h95.55.\h95.56.\h95.57.\h95.58.
\h95.59.\h95.60.\h95.61.\h95.62.\h95.63.\h95.64.\h95.65.\h95.66.
\h95.67.\h95.68.\h95.69.\h95.70.\h95.71.\h95.72.\h95.73.\h95.74.
\h95.75.\h95.76.\h95.77.\h95.78.\h95.79.\h95.80.\h95.81.\h95.82.
\h95.83.\h95.84.\h95.85.\h95.86.\h95.87.\h95.88.\h95.89.\h95.90.
\h95.91.\h95.92.\h95.93.\h95.94.\h95.95.\h96.4.\h96.6.\h96.7.
\h96.8.\h96.9.\h96.10.\h96.11.\h96.12.\h96.13.\h96.14.\h96.15.
\h96.16.\h96.17.\h96.18.\h96.19.\h96.20.\h96.21.\h96.22.\h96.23.
\h96.24.\h96.25.\h96.26.\h96.27.\h96.28.\h96.29.\h96.30.\h96.31.
\h96.32.\h96.33.\h96.34.\h96.35.\h96.36.\h96.37.\h96.38.\h96.39.
\h96.40.\h96.41.\h96.42.\h96.43.\h96.44.\h96.45.\h96.46.\h96.47.
\h96.48.\h96.49.\h96.50.\h96.51.\h96.52.\h96.53.\h96.54.\h96.55.
\h96.56.\h96.57.\h96.58.\h96.59.\h96.60.\h96.61.\h96.62.\h96.63.
\h96.64.\h96.65.\h96.66.\h96.67.\h96.68.\h96.69.\h96.70.\h96.71.
\h96.72.\h96.73.\h96.74.\h96.75.\h96.76.\h96.77.\h96.78.\h96.79.
\h96.80.\h96.81.\h96.82.\h96.83.\h96.84.\h96.86.\h96.87.\h96.88.
\h96.89.\h96.90.\h96.91.\h96.92.\h96.93.\h96.94.\h96.96.\h97.4.
\h97.5.\h97.6.\h97.7.\h97.8.\h97.9.\h97.10.\h97.11.\h97.12.
\h97.13.\h97.14.\h97.15.\h97.16.\h97.17.\h97.18.\h97.19.\h97.20.
\h97.21.\h97.22.\h97.23.\h97.24.\h97.25.\h97.26.\h97.27.\h97.28.
\h97.29.\h97.30.\h97.31.\h97.32.\h97.33.\h97.34.\h97.35.\h97.36.
\h97.37.\h97.38.\h97.39.\h97.40.\h97.41.\h97.42.\h97.43.\h97.44.
\h97.45.\h97.46.\h97.47.\h97.48.\h97.49.\h97.50.\h97.51.\h97.52.
\h97.53.\h97.54.\h97.55.\h97.56.\h97.57.\h97.58.\h97.59.\h97.60.
\h97.61.\h97.62.\h97.63.\h97.64.\h97.65.\h97.66.\h97.67.\h97.68.
\h97.69.\h97.70.\h97.71.\h97.72.\h97.73.\h97.74.\h97.75.\h97.76.
\h97.77.\h97.78.\h97.79.\h97.80.\h97.81.\h97.82.\h97.83.\h97.84.
\h97.85.\h97.86.\h97.87.\h97.88.\h97.89.\h97.90.\h97.91.\h97.92.
\h97.93.\h97.94.\h97.95.\h97.97.\h98.4.\h98.5.\h98.6.\h98.7.
\h98.8.\h98.9.\h98.10.\h98.11.\h98.12.\h98.13.\h98.14.\h98.15.
\h98.16.\h98.17.\h98.18.\h98.19.\h98.20.\h98.21.\h98.22.\h98.23.
\h98.24.\h98.25.\h98.26.\h98.27.\h98.28.\h98.29.\h98.30.\h98.31.
\h98.32.\h98.33.\h98.34.\h98.35.\h98.36.\h98.37.\h98.38.\h98.39.
\h98.40.\h98.41.\h98.42.\h98.43.\h98.44.\h98.45.\h98.46.\h98.47.
\h98.48.\h98.49.\h98.50.\h98.51.\h98.52.\h98.53.\h98.54.\h98.55.
\h98.56.\h98.57.\h98.58.\h98.59.\h98.60.\h98.61.\h98.62.\h98.63.
\h98.64.\h98.65.\h98.66.\h98.67.\h98.68.\h98.69.\h98.70.\h98.71.
\h98.72.\h98.73.\h98.74.\h98.75.\h98.76.\h98.77.\h98.78.\h98.79.
\h98.80.\h98.81.\h98.82.\h98.83.\h98.84.\h98.85.\h98.86.\h98.87.
\h98.88.\h98.89.\h98.90.\h98.92.\h98.93.\h98.94.\h98.95.\h98.96.
\h98.97.\h98.98.\h99.3.\h99.5.\h99.6.\h99.7.\h99.8.\h99.9.
\h99.10.\h99.11.\h99.12.\h99.13.\h99.14.\h99.15.\h99.16.\h99.17.
\h99.18.\h99.19.\h99.20.\h99.21.\h99.22.\h99.23.\h99.24.\h99.25.
\h99.26.\h99.27.\h99.28.\h99.29.\h99.30.\h99.31.\h99.32.\h99.33.
\h99.34.\h99.35.\h99.36.\h99.37.\h99.38.\h99.39.\h99.40.\h99.41.
\h99.42.\h99.43.\h99.44.\h99.45.\h99.46.\h99.47.\h99.48.\h99.49.
\h99.50.\h99.51.\h99.52.\h99.53.\h99.54.\h99.55.\h99.56.\h99.57.
\h99.58.\h99.59.\h99.60.\h99.61.\h99.62.\h99.63.\h99.64.\h99.65.
\h99.66.\h99.67.\h99.68.\h99.69.\h99.70.\h99.71.\h99.72.\h99.73.
\h99.74.\h99.75.\h99.76.\h99.77.\h99.78.\h99.79.\h99.80.\h99.81.
\h99.82.\h99.83.\h99.84.\h99.85.\h99.86.\h99.87.\h99.88.\h99.89.
\h99.90.\h99.91.\h99.92.\h99.93.\h99.95.\h99.96.\h99.97.\h99.98.
\h99.99.\h100.4.\h100.6.\h100.7.\h100.8.\h100.9.\h100.10.\h100.11.
\h100.12.\h100.13.\h100.14.\h100.15.\h100.16.\h100.17.\h100.18.\h100.19.
\h100.20.\h100.21.\h100.22.\h100.23.\h100.24.\h100.25.\h100.26.\h100.27.
\h100.28.\h100.29.\h100.30.\h100.31.\h100.32.\h100.33.\h100.34.\h100.35.
\h100.36.\h100.37.\h100.38.\h100.39.\h100.40.\h100.41.\h100.42.\h100.43.
\h100.44.\h100.45.\h100.46.\h100.47.\h100.48.\h100.49.\h100.50.\h100.51.
\h100.52.\h100.53.\h100.54.\h100.55.\h100.56.\h100.57.\h100.58.\h100.59.
\h100.60.\h100.61.\h100.62.\h100.63.\h100.64.\h100.65.\h100.66.\h100.67.
\h100.68.\h100.69.\h100.70.\h100.71.\h100.72.\h100.73.\h100.74.\h100.75.
\h100.76.\h100.77.\h100.78.\h100.79.\h100.80.\h100.81.\h100.82.\h100.83.
\h100.84.\h100.85.\h100.86.\h100.87.\h100.88.\h100.90.\h100.91.\h100.92.
\h100.94.\h100.95.\h100.96.\h100.97.\h100.98.\h100.99.\h100.100.\h101.1.
\h101.4.\h101.5.\h101.6.\h101.7.\h101.8.\h101.9.\h101.10.\h101.11.
\h101.12.\h101.13.\h101.14.\h101.15.\h101.16.\h101.17.\h101.18.\h101.19.
\h101.20.\h101.21.\h101.22.\h101.23.\h101.24.\h101.25.\h101.26.\h101.27.
\h101.28.\h101.29.\h101.30.\h101.31.\h101.32.\h101.33.\h101.34.\h101.35.
\h101.36.\h101.37.\h101.38.\h101.39.\h101.40.\h101.41.\h101.42.\h101.43.
\h101.44.\h101.45.\h101.46.\h101.47.\h101.48.\h101.49.\h101.50.\h101.51.
\h101.52.\h101.53.\h101.54.\h101.55.\h101.56.\h101.57.\h101.58.\h101.59.
\h101.60.\h101.61.\h101.62.\h101.63.\h101.64.\h101.65.\h101.66.\h101.67.
\h101.68.\h101.69.\h101.70.\h101.71.\h101.72.\h101.73.\h101.74.\h101.75.
\h101.76.\h101.77.\h101.78.\h101.79.\h101.80.\h101.81.\h101.82.\h101.83.
\h101.84.\h101.85.\h101.86.\h101.87.\h101.88.\h101.89.\h101.90.\h101.91.
\h101.92.\h101.93.\h101.94.\h101.95.\h101.96.\h101.97.\h101.98.\h101.99.
\h101.101.\h102.2.\h102.4.\h102.6.\h102.7.\h102.8.\h102.9.\h102.10.
\h102.11.\h102.12.\h102.13.\h102.14.\h102.15.\h102.16.\h102.17.\h102.18.
\h102.19.\h102.20.\h102.21.\h102.22.\h102.23.\h102.24.\h102.25.\h102.26.
\h102.27.\h102.28.\h102.29.\h102.30.\h102.31.\h102.32.\h102.33.\h102.34.
\h102.35.\h102.36.\h102.37.\h102.38.\h102.39.\h102.40.\h102.41.\h102.42.
\h102.43.\h102.44.\h102.45.\h102.46.\h102.47.\h102.48.\h102.49.\h102.50.
\h102.51.\h102.52.\h102.53.\h102.54.\h102.55.\h102.56.\h102.57.\h102.58.
\h102.59.\h102.60.\h102.61.\h102.62.\h102.63.\h102.64.\h102.65.\h102.66.
\h102.67.\h102.68.\h102.69.\h102.70.\h102.71.\h102.72.\h102.73.\h102.74.
\h102.75.\h102.76.\h102.77.\h102.78.\h102.79.\h102.80.\h102.81.\h102.82.
\h102.83.\h102.84.\h102.85.\h102.86.\h102.87.\h102.88.\h102.89.\h102.90.
\h102.91.\h102.92.\h102.93.\h102.94.\h102.96.\h102.97.\h102.98.\h102.99.
\h102.100.\h102.102.\h103.1.\h103.3.\h103.5.\h103.6.\h103.7.\h103.8.
\h103.9.\h103.10.\h103.11.\h103.12.\h103.13.\h103.14.\h103.15.\h103.16.
\h103.17.\h103.18.\h103.19.\h103.20.\h103.21.\h103.22.\h103.23.\h103.24.
\h103.25.\h103.26.\h103.27.\h103.28.\h103.29.\h103.30.\h103.31.\h103.32.
\h103.33.\h103.34.\h103.35.\h103.36.\h103.37.\h103.38.\h103.39.\h103.40.
\h103.41.\h103.42.\h103.43.\h103.44.\h103.45.\h103.46.\h103.47.\h103.48.
\h103.49.\h103.50.\h103.51.\h103.52.\h103.53.\h103.54.\h103.55.\h103.56.
\h103.57.\h103.58.\h103.59.\h103.60.\h103.61.\h103.62.\h103.63.\h103.64.
\h103.65.\h103.66.\h103.67.\h103.68.\h103.69.\h103.70.\h103.71.\h103.72.
\h103.73.\h103.74.\h103.75.\h103.76.\h103.77.\h103.78.\h103.79.\h103.80.
\h103.81.\h103.82.\h103.83.\h103.84.\h103.85.\h103.86.\h103.87.\h103.88.
\h103.89.\h103.90.\h103.91.\h103.92.\h103.93.\h103.94.\h103.95.\h103.96.
\h103.97.\h103.98.\h103.99.\h103.100.\h103.101.\h103.103.\h104.4.\h104.5.
\h104.6.\h104.7.\h104.8.\h104.9.\h104.10.\h104.11.\h104.12.\h104.13.
\h104.14.\h104.15.\h104.16.\h104.17.\h104.18.\h104.19.\h104.20.\h104.21.
\h104.22.\h104.23.\h104.24.\h104.25.\h104.26.\h104.27.\h104.28.\h104.29.
\h104.30.\h104.31.\h104.32.\h104.33.\h104.34.\h104.35.\h104.36.\h104.37.
\h104.38.\h104.39.\h104.40.\h104.41.\h104.42.\h104.43.\h104.44.\h104.45.
\h104.46.\h104.47.\h104.48.\h104.49.\h104.50.\h104.51.\h104.52.\h104.53.
\h104.54.\h104.55.\h104.56.\h104.57.\h104.58.\h104.59.\h104.60.\h104.61.
\h104.62.\h104.63.\h104.64.\h104.65.\h104.66.\h104.67.\h104.68.\h104.69.
\h104.70.\h104.71.\h104.72.\h104.73.\h104.74.\h104.75.\h104.76.\h104.77.
\h104.78.\h104.79.\h104.80.\h104.81.\h104.82.\h104.83.\h104.84.\h104.85.
\h104.86.\h104.87.\h104.88.\h104.89.\h104.90.\h104.91.\h104.92.\h104.94.
\h104.95.\h104.96.\h104.97.\h104.98.\h104.99.\h104.100.\h104.101.\h104.102.
\h104.104.\h105.3.\h105.5.\h105.6.\h105.7.\h105.8.\h105.9.\h105.10.
\h105.11.\h105.12.\h105.13.\h105.14.\h105.15.\h105.16.\h105.17.\h105.18.
\h105.19.\h105.20.\h105.21.\h105.22.\h105.23.\h105.24.\h105.25.\h105.26.
\h105.27.\h105.28.\h105.29.\h105.30.\h105.31.\h105.32.\h105.33.\h105.34.
\h105.35.\h105.36.\h105.37.\h105.38.\h105.39.\h105.40.\h105.41.\h105.42.
\h105.43.\h105.44.\h105.45.\h105.46.\h105.47.\h105.48.\h105.49.\h105.50.
\h105.51.\h105.52.\h105.53.\h105.54.\h105.55.\h105.56.\h105.57.\h105.58.
\h105.59.\h105.60.\h105.61.\h105.62.\h105.63.\h105.64.\h105.65.\h105.66.
\h105.67.\h105.68.\h105.69.\h105.70.\h105.71.\h105.72.\h105.73.\h105.74.
\h105.75.\h105.76.\h105.77.\h105.78.\h105.79.\h105.80.\h105.81.\h105.82.
\h105.83.\h105.84.\h105.85.\h105.87.\h105.88.\h105.89.\h105.90.\h105.91.
\h105.93.\h105.95.\h105.96.\h105.97.\h105.98.\h105.99.\h105.101.\h105.103.
\h105.105.\h106.2.\h106.4.\h106.6.\h106.7.\h106.8.\h106.9.\h106.10.
\h106.11.\h106.12.\h106.13.\h106.14.\h106.15.\h106.16.\h106.17.\h106.18.
\h106.19.\h106.20.\h106.21.\h106.22.\h106.23.\h106.24.\h106.25.\h106.26.
\h106.27.\h106.28.\h106.29.\h106.30.\h106.31.\h106.32.\h106.33.\h106.34.
\h106.35.\h106.36.\h106.37.\h106.38.\h106.39.\h106.40.\h106.41.\h106.42.
\h106.43.\h106.44.\h106.45.\h106.46.\h106.47.\h106.48.\h106.49.\h106.50.
\h106.51.\h106.52.\h106.53.\h106.54.\h106.55.\h106.56.\h106.57.\h106.58.
\h106.59.\h106.60.\h106.61.\h106.62.\h106.63.\h106.64.\h106.65.\h106.66.
\h106.67.\h106.68.\h106.69.\h106.70.\h106.71.\h106.72.\h106.73.\h106.74.
\h106.75.\h106.76.\h106.77.\h106.78.\h106.79.\h106.80.\h106.81.\h106.82.
\h106.83.\h106.84.\h106.86.\h106.88.\h106.89.\h106.90.\h106.91.\h106.92.
\h106.93.\h106.94.\h106.96.\h106.97.\h106.98.\h106.100.\h106.102.\h106.104.
\h106.106.\h107.3.\h107.5.\h107.6.\h107.7.\h107.8.\h107.9.\h107.10.
\h107.11.\h107.12.\h107.13.\h107.14.\h107.15.\h107.16.\h107.17.\h107.18.
\h107.19.\h107.20.\h107.21.\h107.22.\h107.23.\h107.24.\h107.25.\h107.26.
\h107.27.\h107.28.\h107.29.\h107.30.\h107.31.\h107.32.\h107.33.\h107.34.
\h107.35.\h107.36.\h107.37.\h107.38.\h107.39.\h107.40.\h107.41.\h107.42.
\h107.43.\h107.44.\h107.45.\h107.46.\h107.47.\h107.48.\h107.49.\h107.50.
\h107.51.\h107.52.\h107.53.\h107.54.\h107.55.\h107.56.\h107.57.\h107.58.
\h107.59.\h107.60.\h107.61.\h107.62.\h107.63.\h107.64.\h107.65.\h107.66.
\h107.67.\h107.68.\h107.69.\h107.70.\h107.71.\h107.72.\h107.73.\h107.74.
\h107.75.\h107.76.\h107.77.\h107.79.\h107.80.\h107.81.\h107.82.\h107.83.
\h107.84.\h107.85.\h107.86.\h107.87.\h107.88.\h107.89.\h107.91.\h107.92.
\h107.93.\h107.95.\h107.97.\h107.98.\h107.99.\h107.100.\h107.101.\h107.102.
\h107.103.\h107.104.\h107.105.\h107.107.\h108.4.\h108.5.\h108.6.\h108.7.
\h108.8.\h108.9.\h108.10.\h108.11.\h108.12.\h108.13.\h108.14.\h108.15.
\h108.16.\h108.17.\h108.18.\h108.19.\h108.20.\h108.21.\h108.22.\h108.23.
\h108.24.\h108.25.\h108.26.\h108.27.\h108.28.\h108.29.\h108.30.\h108.31.
\h108.32.\h108.33.\h108.34.\h108.35.\h108.36.\h108.37.\h108.38.\h108.39.
\h108.40.\h108.41.\h108.42.\h108.43.\h108.44.\h108.45.\h108.46.\h108.47.
\h108.48.\h108.49.\h108.50.\h108.51.\h108.52.\h108.53.\h108.54.\h108.55.
\h108.56.\h108.57.\h108.58.\h108.59.\h108.60.\h108.61.\h108.62.\h108.63.
\h108.64.\h108.65.\h108.66.\h108.67.\h108.68.\h108.69.\h108.70.\h108.71.
\h108.72.\h108.73.\h108.74.\h108.75.\h108.76.\h108.78.\h108.79.\h108.80.
\h108.81.\h108.82.\h108.83.\h108.84.\h108.86.\h108.87.\h108.88.\h108.89.
\h108.90.\h108.92.\h108.93.\h108.94.\h108.96.\h108.98.\h108.99.\h108.100.
\h108.101.\h108.102.\h108.104.\h108.105.\h108.106.\h108.107.\h108.108.\h109.4.
\h109.5.\h109.6.\h109.7.\h109.8.\h109.9.\h109.10.\h109.11.\h109.12.
\h109.13.\h109.14.\h109.15.\h109.16.\h109.17.\h109.18.\h109.19.\h109.20.
\h109.21.\h109.22.\h109.23.\h109.24.\h109.25.\h109.26.\h109.27.\h109.28.
\h109.29.\h109.30.\h109.31.\h109.32.\h109.33.\h109.34.\h109.35.\h109.36.
\h109.37.\h109.38.\h109.39.\h109.40.\h109.41.\h109.42.\h109.43.\h109.44.
\h109.45.\h109.46.\h109.47.\h109.48.\h109.49.\h109.50.\h109.51.\h109.52.
\h109.53.\h109.54.\h109.55.\h109.56.\h109.57.\h109.58.\h109.59.\h109.60.
\h109.61.\h109.62.\h109.63.\h109.64.\h109.65.\h109.66.\h109.67.\h109.68.
\h109.69.\h109.70.\h109.71.\h109.72.\h109.73.\h109.74.\h109.75.\h109.76.
\h109.77.\h109.78.\h109.79.\h109.80.\h109.81.\h109.82.\h109.83.\h109.85.
\h109.86.\h109.87.\h109.88.\h109.89.\h109.91.\h109.92.\h109.93.\h109.94.
\h109.95.\h109.96.\h109.97.\h109.98.\h109.99.\h109.100.\h109.101.\h109.103.
\h109.104.\h109.105.\h109.106.\h109.107.\h109.109.\h110.4.\h110.5.\h110.6.
\h110.7.\h110.8.\h110.9.\h110.10.\h110.11.\h110.12.\h110.13.\h110.14.
\h110.15.\h110.16.\h110.17.\h110.18.\h110.19.\h110.20.\h110.21.\h110.22.
\h110.23.\h110.24.\h110.25.\h110.26.\h110.27.\h110.28.\h110.29.\h110.30.
\h110.31.\h110.32.\h110.33.\h110.34.\h110.35.\h110.36.\h110.37.\h110.38.
\h110.39.\h110.40.\h110.41.\h110.42.\h110.43.\h110.44.\h110.45.\h110.46.
\h110.47.\h110.48.\h110.49.\h110.50.\h110.51.\h110.52.\h110.53.\h110.54.
\h110.55.\h110.56.\h110.57.\h110.58.\h110.59.\h110.60.\h110.61.\h110.62.
\h110.63.\h110.64.\h110.65.\h110.66.\h110.67.\h110.68.\h110.69.\h110.70.
\h110.71.\h110.72.\h110.73.\h110.74.\h110.75.\h110.76.\h110.77.\h110.78.
\h110.79.\h110.80.\h110.81.\h110.82.\h110.83.\h110.84.\h110.85.\h110.86.
\h110.87.\h110.88.\h110.89.\h110.90.\h110.91.\h110.92.\h110.93.\h110.94.
\h110.95.\h110.96.\h110.97.\h110.98.\h110.100.\h110.102.\h110.103.\h110.104.
\h110.105.\h110.106.\h110.107.\h110.108.\h110.110.\h111.3.\h111.5.\h111.6.
\h111.7.\h111.8.\h111.9.\h111.10.\h111.11.\h111.12.\h111.13.\h111.14.
\h111.15.\h111.16.\h111.17.\h111.18.\h111.19.\h111.20.\h111.21.\h111.22.
\h111.23.\h111.24.\h111.25.\h111.26.\h111.27.\h111.28.\h111.29.\h111.30.
\h111.31.\h111.32.\h111.33.\h111.34.\h111.35.\h111.36.\h111.37.\h111.38.
\h111.39.\h111.40.\h111.41.\h111.42.\h111.43.\h111.44.\h111.45.\h111.46.
\h111.47.\h111.48.\h111.49.\h111.50.\h111.51.\h111.52.\h111.53.\h111.54.
\h111.55.\h111.56.\h111.57.\h111.58.\h111.59.\h111.60.\h111.61.\h111.62.
\h111.63.\h111.64.\h111.65.\h111.66.\h111.67.\h111.68.\h111.69.\h111.70.
\h111.71.\h111.72.\h111.73.\h111.74.\h111.75.\h111.76.\h111.77.\h111.78.
\h111.79.\h111.80.\h111.81.\h111.82.\h111.83.\h111.84.\h111.85.\h111.86.
\h111.87.\h111.89.\h111.90.\h111.91.\h111.92.\h111.93.\h111.94.\h111.95.
\h111.96.\h111.97.\h111.99.\h111.101.\h111.102.\h111.103.\h111.104.\h111.105.
\h111.107.\h111.108.\h111.109.\h111.111.\h112.4.\h112.6.\h112.7.\h112.8.
\h112.9.\h112.10.\h112.11.\h112.12.\h112.13.\h112.14.\h112.15.\h112.16.
\h112.17.\h112.18.\h112.19.\h112.20.\h112.21.\h112.22.\h112.23.\h112.24.
\h112.25.\h112.26.\h112.27.\h112.28.\h112.29.\h112.30.\h112.31.\h112.32.
\h112.33.\h112.34.\h112.35.\h112.36.\h112.37.\h112.38.\h112.39.\h112.40.
\h112.41.\h112.42.\h112.43.\h112.44.\h112.45.\h112.46.\h112.47.\h112.48.
\h112.49.\h112.50.\h112.51.\h112.52.\h112.53.\h112.54.\h112.55.\h112.56.
\h112.57.\h112.58.\h112.59.\h112.60.\h112.61.\h112.62.\h112.63.\h112.64.
\h112.65.\h112.66.\h112.67.\h112.68.\h112.69.\h112.70.\h112.71.\h112.72.
\h112.73.\h112.74.\h112.75.\h112.76.\h112.77.\h112.78.\h112.79.\h112.80.
\h112.82.\h112.83.\h112.84.\h112.85.\h112.86.\h112.88.\h112.89.\h112.90.
\h112.91.\h112.92.\h112.94.\h112.96.\h112.97.\h112.98.\h112.100.\h112.102.
\h112.103.\h112.104.\h112.106.\h112.108.\h112.109.\h112.110.\h112.112.\h113.5.
\h113.7.\h113.8.\h113.9.\h113.10.\h113.11.\h113.12.\h113.13.\h113.14.
\h113.15.\h113.16.\h113.17.\h113.18.\h113.19.\h113.20.\h113.21.\h113.22.
\h113.23.\h113.24.\h113.25.\h113.26.\h113.27.\h113.28.\h113.29.\h113.30.
\h113.31.\h113.32.\h113.33.\h113.34.\h113.35.\h113.36.\h113.37.\h113.38.
\h113.39.\h113.40.\h113.41.\h113.42.\h113.43.\h113.44.\h113.45.\h113.46.
\h113.47.\h113.48.\h113.49.\h113.50.\h113.51.\h113.52.\h113.53.\h113.54.
\h113.55.\h113.56.\h113.57.\h113.58.\h113.59.\h113.60.\h113.61.\h113.62.
\h113.63.\h113.64.\h113.65.\h113.66.\h113.67.\h113.68.\h113.69.\h113.70.
\h113.71.\h113.72.\h113.73.\h113.74.\h113.75.\h113.77.\h113.78.\h113.79.
\h113.80.\h113.81.\h113.82.\h113.83.\h113.84.\h113.85.\h113.86.\h113.87.
\h113.88.\h113.89.\h113.90.\h113.91.\h113.93.\h113.95.\h113.97.\h113.98.
\h113.99.\h113.101.\h113.103.\h113.104.\h113.105.\h113.107.\h113.109.\h113.110.
\h113.111.\h113.113.\h114.4.\h114.6.\h114.7.\h114.8.\h114.9.\h114.10.
\h114.11.\h114.12.\h114.13.\h114.14.\h114.15.\h114.16.\h114.17.\h114.18.
\h114.19.\h114.20.\h114.21.\h114.22.\h114.23.\h114.24.\h114.25.\h114.26.
\h114.27.\h114.28.\h114.29.\h114.30.\h114.31.\h114.32.\h114.33.\h114.34.
\h114.35.\h114.36.\h114.37.\h114.38.\h114.39.\h114.40.\h114.41.\h114.42.
\h114.43.\h114.44.\h114.45.\h114.46.\h114.47.\h114.48.\h114.49.\h114.50.
\h114.51.\h114.52.\h114.53.\h114.54.\h114.55.\h114.56.\h114.57.\h114.58.
\h114.59.\h114.60.\h114.61.\h114.62.\h114.63.\h114.64.\h114.65.\h114.66.
\h114.67.\h114.68.\h114.69.\h114.70.\h114.71.\h114.72.\h114.73.\h114.74.
\h114.75.\h114.76.\h114.77.\h114.78.\h114.80.\h114.81.\h114.82.\h114.84.
\h114.85.\h114.86.\h114.87.\h114.88.\h114.89.\h114.90.\h114.92.\h114.93.
\h114.94.\h114.96.\h114.97.\h114.98.\h114.99.\h114.100.\h114.102.\h114.104.
\h114.106.\h114.108.\h114.109.\h114.110.\h114.111.\h114.112.\h114.114.\h115.3.
\h115.5.\h115.6.\h115.7.\h115.8.\h115.9.\h115.10.\h115.11.\h115.12.
\h115.13.\h115.14.\h115.15.\h115.16.\h115.17.\h115.18.\h115.19.\h115.20.
\h115.21.\h115.22.\h115.23.\h115.24.\h115.25.\h115.26.\h115.27.\h115.28.
\h115.29.\h115.30.\h115.31.\h115.32.\h115.33.\h115.34.\h115.35.\h115.36.
\h115.37.\h115.38.\h115.39.\h115.40.\h115.41.\h115.42.\h115.43.\h115.44.
\h115.45.\h115.46.\h115.47.\h115.48.\h115.49.\h115.50.\h115.51.\h115.52.
\h115.53.\h115.54.\h115.55.\h115.56.\h115.57.\h115.58.\h115.59.\h115.60.
\h115.61.\h115.62.\h115.63.\h115.64.\h115.65.\h115.66.\h115.67.\h115.68.
\h115.69.\h115.70.\h115.71.\h115.72.\h115.73.\h115.74.\h115.75.\h115.76.
\h115.77.\h115.79.\h115.80.\h115.81.\h115.82.\h115.83.\h115.85.\h115.86.
\h115.87.\h115.88.\h115.89.\h115.91.\h115.93.\h115.95.\h115.96.\h115.97.
\h115.99.\h115.100.\h115.101.\h115.103.\h115.105.\h115.106.\h115.107.\h115.108.
\h115.109.\h115.111.\h115.112.\h115.113.\h115.115.\h116.2.\h116.4.\h116.5.
\h116.6.\h116.7.\h116.8.\h116.9.\h116.10.\h116.11.\h116.12.\h116.13.
\h116.14.\h116.15.\h116.16.\h116.17.\h116.18.\h116.19.\h116.20.\h116.21.
\h116.22.\h116.23.\h116.24.\h116.25.\h116.26.\h116.27.\h116.28.\h116.29.
\h116.30.\h116.31.\h116.32.\h116.33.\h116.34.\h116.35.\h116.36.\h116.37.
\h116.38.\h116.39.\h116.40.\h116.41.\h116.42.\h116.43.\h116.44.\h116.45.
\h116.46.\h116.47.\h116.48.\h116.49.\h116.50.\h116.51.\h116.52.\h116.53.
\h116.54.\h116.55.\h116.56.\h116.57.\h116.58.\h116.59.\h116.60.\h116.61.
\h116.62.\h116.63.\h116.64.\h116.65.\h116.66.\h116.67.\h116.68.\h116.69.
\h116.70.\h116.71.\h116.72.\h116.74.\h116.75.\h116.76.\h116.77.\h116.78.
\h116.79.\h116.80.\h116.82.\h116.83.\h116.84.\h116.86.\h116.88.\h116.89.
\h116.90.\h116.92.\h116.94.\h116.95.\h116.96.\h116.98.\h116.99.\h116.100.
\h116.101.\h116.102.\h116.104.\h116.106.\h116.107.\h116.108.\h116.109.\h116.110.
\h116.112.\h116.113.\h116.114.\h116.115.\h116.116.\h117.5.\h117.6.\h117.7.
\h117.9.\h117.10.\h117.11.\h117.12.\h117.13.\h117.14.\h117.15.\h117.16.
\h117.17.\h117.18.\h117.19.\h117.20.\h117.21.\h117.22.\h117.23.\h117.24.
\h117.25.\h117.26.\h117.27.\h117.28.\h117.29.\h117.30.\h117.31.\h117.32.
\h117.33.\h117.34.\h117.35.\h117.36.\h117.37.\h117.38.\h117.39.\h117.40.
\h117.41.\h117.42.\h117.43.\h117.44.\h117.45.\h117.46.\h117.47.\h117.48.
\h117.49.\h117.50.\h117.51.\h117.52.\h117.53.\h117.54.\h117.55.\h117.56.
\h117.57.\h117.58.\h117.59.\h117.60.\h117.61.\h117.62.\h117.63.\h117.64.
\h117.65.\h117.66.\h117.67.\h117.68.\h117.69.\h117.70.\h117.71.\h117.72.
\h117.73.\h117.74.\h117.75.\h117.76.\h117.77.\h117.78.\h117.79.\h117.80.
\h117.81.\h117.82.\h117.83.\h117.84.\h117.85.\h117.86.\h117.87.\h117.89.
\h117.90.\h117.91.\h117.93.\h117.95.\h117.96.\h117.97.\h117.98.\h117.99.
\h117.101.\h117.102.\h117.103.\h117.105.\h117.106.\h117.107.\h117.108.\h117.109.
\h117.111.\h117.112.\h117.113.\h117.114.\h117.117.\h118.4.\h118.6.\h118.8.
\h118.9.\h118.10.\h118.11.\h118.12.\h118.13.\h118.14.\h118.15.\h118.16.
\h118.17.\h118.18.\h118.19.\h118.20.\h118.21.\h118.22.\h118.23.\h118.24.
\h118.25.\h118.26.\h118.27.\h118.28.\h118.29.\h118.30.\h118.31.\h118.32.
\h118.33.\h118.34.\h118.35.\h118.36.\h118.37.\h118.38.\h118.39.\h118.40.
\h118.41.\h118.42.\h118.43.\h118.44.\h118.45.\h118.46.\h118.47.\h118.48.
\h118.49.\h118.50.\h118.51.\h118.52.\h118.53.\h118.54.\h118.55.\h118.56.
\h118.57.\h118.58.\h118.59.\h118.60.\h118.61.\h118.62.\h118.63.\h118.64.
\h118.65.\h118.66.\h118.67.\h118.68.\h118.69.\h118.70.\h118.71.\h118.72.
\h118.73.\h118.74.\h118.75.\h118.76.\h118.77.\h118.78.\h118.79.\h118.80.
\h118.81.\h118.82.\h118.83.\h118.84.\h118.85.\h118.86.\h118.87.\h118.88.
\h118.89.\h118.90.\h118.91.\h118.92.\h118.94.\h118.95.\h118.96.\h118.97.
\h118.98.\h118.100.\h118.101.\h118.102.\h118.103.\h118.104.\h118.105.\h118.106.
\h118.107.\h118.108.\h118.110.\h118.111.\h118.112.\h118.113.\h118.114.\h118.116.
\h119.3.\h119.5.\h119.7.\h119.8.\h119.9.\h119.10.\h119.11.\h119.12.
\h119.13.\h119.14.\h119.15.\h119.16.\h119.17.\h119.18.\h119.19.\h119.20.
\h119.21.\h119.22.\h119.23.\h119.24.\h119.25.\h119.26.\h119.27.\h119.28.
\h119.29.\h119.30.\h119.31.\h119.32.\h119.33.\h119.34.\h119.35.\h119.36.
\h119.37.\h119.38.\h119.39.\h119.40.\h119.41.\h119.42.\h119.43.\h119.44.
\h119.45.\h119.46.\h119.47.\h119.48.\h119.49.\h119.50.\h119.51.\h119.52.
\h119.53.\h119.54.\h119.55.\h119.56.\h119.57.\h119.58.\h119.59.\h119.60.
\h119.61.\h119.62.\h119.63.\h119.64.\h119.65.\h119.66.\h119.67.\h119.68.
\h119.69.\h119.70.\h119.71.\h119.72.\h119.73.\h119.74.\h119.75.\h119.76.
\h119.77.\h119.78.\h119.79.\h119.80.\h119.81.\h119.82.\h119.83.\h119.84.
\h119.85.\h119.86.\h119.87.\h119.89.\h119.91.\h119.92.\h119.93.\h119.94.
\h119.95.\h119.96.\h119.97.\h119.99.\h119.100.\h119.101.\h119.103.\h119.104.
\h119.105.\h119.107.\h119.109.\h119.110.\h119.111.\h119.113.\h119.115.\h119.116.
\h119.117.\h119.119.\h120.2.\h120.4.\h120.6.\h120.7.\h120.8.\h120.9.
\h120.10.\h120.11.\h120.12.\h120.13.\h120.14.\h120.15.\h120.16.\h120.17.
\h120.18.\h120.19.\h120.20.\h120.21.\h120.22.\h120.23.\h120.24.\h120.25.
\h120.26.\h120.27.\h120.28.\h120.29.\h120.30.\h120.31.\h120.32.\h120.33.
\h120.34.\h120.35.\h120.36.\h120.37.\h120.38.\h120.39.\h120.40.\h120.41.
\h120.42.\h120.43.\h120.44.\h120.45.\h120.46.\h120.47.\h120.48.\h120.49.
\h120.50.\h120.51.\h120.52.\h120.53.\h120.54.\h120.55.\h120.56.\h120.57.
\h120.58.\h120.59.\h120.60.\h120.61.\h120.62.\h120.63.\h120.64.\h120.65.
\h120.66.\h120.67.\h120.68.\h120.69.\h120.70.\h120.71.\h120.72.\h120.74.
\h120.75.\h120.76.\h120.78.\h120.79.\h120.80.\h120.81.\h120.82.\h120.83.
\h120.84.\h120.85.\h120.86.\h120.87.\h120.88.\h120.90.\h120.91.\h120.92.
\h120.93.\h120.94.\h120.95.\h120.96.\h120.98.\h120.99.\h120.100.\h120.102.
\h120.104.\h120.105.\h120.106.\h120.108.\h120.110.\h120.111.\h120.114.\h120.116.
\h120.118.\h120.120.\h121.4.\h121.5.\h121.7.\h121.8.\h121.9.\h121.10.
\h121.11.\h121.12.\h121.13.\h121.14.\h121.15.\h121.16.\h121.17.\h121.18.
\h121.19.\h121.20.\h121.21.\h121.22.\h121.23.\h121.24.\h121.25.\h121.26.
\h121.27.\h121.28.\h121.29.\h121.30.\h121.31.\h121.32.\h121.33.\h121.34.
\h121.35.\h121.36.\h121.37.\h121.38.\h121.39.\h121.40.\h121.41.\h121.42.
\h121.43.\h121.44.\h121.45.\h121.46.\h121.47.\h121.48.\h121.49.\h121.50.
\h121.51.\h121.52.\h121.53.\h121.54.\h121.55.\h121.56.\h121.57.\h121.58.
\h121.59.\h121.60.\h121.61.\h121.62.\h121.63.\h121.64.\h121.65.\h121.66.
\h121.67.\h121.68.\h121.69.\h121.70.\h121.71.\h121.73.\h121.74.\h121.75.
\h121.76.\h121.77.\h121.78.\h121.79.\h121.81.\h121.82.\h121.83.\h121.85.
\h121.87.\h121.89.\h121.90.\h121.91.\h121.93.\h121.94.\h121.95.\h121.97.
\h121.98.\h121.99.\h121.100.\h121.101.\h121.103.\h121.105.\h121.107.\h121.109.
\h121.113.\h121.115.\h121.117.\h121.118.\h121.121.\h122.2.\h122.4.\h122.5.
\h122.6.\h122.7.\h122.8.\h122.9.\h122.10.\h122.11.\h122.12.\h122.13.
\h122.14.\h122.15.\h122.16.\h122.17.\h122.18.\h122.19.\h122.20.\h122.21.
\h122.22.\h122.23.\h122.24.\h122.25.\h122.26.\h122.27.\h122.28.\h122.29.
\h122.30.\h122.31.\h122.32.\h122.33.\h122.34.\h122.35.\h122.36.\h122.37.
\h122.38.\h122.39.\h122.40.\h122.41.\h122.42.\h122.43.\h122.44.\h122.45.
\h122.46.\h122.47.\h122.48.\h122.49.\h122.50.\h122.51.\h122.52.\h122.53.
\h122.54.\h122.55.\h122.56.\h122.57.\h122.58.\h122.59.\h122.60.\h122.61.
\h122.62.\h122.63.\h122.64.\h122.65.\h122.66.\h122.67.\h122.68.\h122.69.
\h122.70.\h122.71.\h122.72.\h122.74.\h122.76.\h122.77.\h122.78.\h122.80.
\h122.82.\h122.83.\h122.84.\h122.85.\h122.86.\h122.87.\h122.88.\h122.89.
\h122.90.\h122.92.\h122.93.\h122.94.\h122.96.\h122.97.\h122.98.\h122.100.
\h122.102.\h122.104.\h122.106.\h122.108.\h122.110.\h122.112.\h122.116.\h122.118.
\h122.120.\h122.122.\h123.3.\h123.5.\h123.7.\h123.8.\h123.9.\h123.10.
\h123.11.\h123.12.\h123.13.\h123.14.\h123.15.\h123.16.\h123.17.\h123.18.
\h123.19.\h123.20.\h123.21.\h123.22.\h123.23.\h123.24.\h123.25.\h123.26.
\h123.27.\h123.28.\h123.29.\h123.30.\h123.31.\h123.32.\h123.33.\h123.34.
\h123.35.\h123.36.\h123.37.\h123.38.\h123.39.\h123.40.\h123.41.\h123.42.
\h123.43.\h123.44.\h123.45.\h123.46.\h123.47.\h123.48.\h123.49.\h123.50.
\h123.51.\h123.52.\h123.53.\h123.54.\h123.55.\h123.56.\h123.57.\h123.58.
\h123.59.\h123.60.\h123.61.\h123.62.\h123.63.\h123.64.\h123.65.\h123.66.
\h123.67.\h123.68.\h123.69.\h123.70.\h123.71.\h123.72.\h123.73.\h123.75.
\h123.76.\h123.77.\h123.78.\h123.79.\h123.81.\h123.83.\h123.84.\h123.85.
\h123.86.\h123.87.\h123.89.\h123.90.\h123.91.\h123.93.\h123.95.\h123.96.
\h123.97.\h123.99.\h123.101.\h123.103.\h123.105.\h123.107.\h123.109.\h123.111.
\h123.115.\h123.117.\h123.119.\h123.123.\h124.4.\h124.6.\h124.7.\h124.8.
\h124.9.\h124.10.\h124.11.\h124.12.\h124.13.\h124.14.\h124.15.\h124.16.
\h124.17.\h124.18.\h124.19.\h124.20.\h124.21.\h124.22.\h124.23.\h124.24.
\h124.25.\h124.26.\h124.27.\h124.28.\h124.29.\h124.30.\h124.31.\h124.32.
\h124.33.\h124.34.\h124.35.\h124.36.\h124.37.\h124.38.\h124.39.\h124.40.
\h124.41.\h124.42.\h124.43.\h124.44.\h124.45.\h124.46.\h124.47.\h124.48.
\h124.49.\h124.50.\h124.51.\h124.52.\h124.53.\h124.54.\h124.55.\h124.56.
\h124.57.\h124.58.\h124.59.\h124.60.\h124.61.\h124.62.\h124.63.\h124.64.
\h124.65.\h124.66.\h124.67.\h124.68.\h124.69.\h124.70.\h124.71.\h124.72.
\h124.73.\h124.74.\h124.75.\h124.76.\h124.77.\h124.78.\h124.79.\h124.80.
\h124.81.\h124.82.\h124.83.\h124.84.\h124.85.\h124.86.\h124.88.\h124.90.
\h124.91.\h124.92.\h124.94.\h124.96.\h124.98.\h124.100.\h124.102.\h124.104.
\h124.106.\h124.108.\h124.109.\h124.110.\h124.112.\h124.114.\h124.115.\h124.116.
\h124.117.\h124.118.\h124.122.\h125.5.\h125.7.\h125.8.\h125.9.\h125.10.
\h125.11.\h125.12.\h125.13.\h125.14.\h125.15.\h125.16.\h125.17.\h125.18.
\h125.19.\h125.20.\h125.21.\h125.22.\h125.23.\h125.24.\h125.25.\h125.26.
\h125.27.\h125.28.\h125.29.\h125.30.\h125.31.\h125.32.\h125.33.\h125.34.
\h125.35.\h125.36.\h125.37.\h125.38.\h125.39.\h125.40.\h125.41.\h125.42.
\h125.43.\h125.44.\h125.45.\h125.46.\h125.47.\h125.48.\h125.49.\h125.50.
\h125.51.\h125.52.\h125.53.\h125.54.\h125.55.\h125.56.\h125.57.\h125.58.
\h125.59.\h125.60.\h125.61.\h125.62.\h125.63.\h125.64.\h125.65.\h125.66.
\h125.67.\h125.68.\h125.69.\h125.70.\h125.71.\h125.73.\h125.74.\h125.75.
\h125.76.\h125.77.\h125.79.\h125.80.\h125.81.\h125.82.\h125.83.\h125.85.
\h125.86.\h125.87.\h125.88.\h125.89.\h125.91.\h125.93.\h125.95.\h125.97.
\h125.98.\h125.99.\h125.101.\h125.103.\h125.104.\h125.105.\h125.107.\h125.109.
\h125.113.\h125.114.\h125.115.\h125.116.\h125.117.\h125.119.\h125.121.\h126.4.
\h126.6.\h126.8.\h126.9.\h126.10.\h126.11.\h126.12.\h126.13.\h126.14.
\h126.15.\h126.16.\h126.17.\h126.18.\h126.19.\h126.20.\h126.21.\h126.22.
\h126.23.\h126.24.\h126.25.\h126.26.\h126.27.\h126.28.\h126.29.\h126.30.
\h126.31.\h126.32.\h126.33.\h126.34.\h126.35.\h126.36.\h126.37.\h126.38.
\h126.39.\h126.40.\h126.41.\h126.42.\h126.43.\h126.44.\h126.45.\h126.46.
\h126.47.\h126.48.\h126.49.\h126.50.\h126.51.\h126.52.\h126.53.\h126.54.
\h126.55.\h126.56.\h126.57.\h126.58.\h126.59.\h126.60.\h126.61.\h126.62.
\h126.63.\h126.64.\h126.65.\h126.66.\h126.68.\h126.69.\h126.70.\h126.71.
\h126.72.\h126.73.\h126.74.\h126.75.\h126.76.\h126.78.\h126.80.\h126.81.
\h126.82.\h126.84.\h126.85.\h126.86.\h126.87.\h126.88.\h126.89.\h126.90.
\h126.92.\h126.93.\h126.94.\h126.96.\h126.98.\h126.100.\h126.102.\h126.103.
\h126.104.\h126.106.\h126.108.\h126.112.\h126.113.\h126.114.\h126.116.\h126.117.
\h126.118.\h126.120.\h126.122.\h126.126.\h127.3.\h127.5.\h127.6.\h127.7.
\h127.9.\h127.10.\h127.11.\h127.12.\h127.13.\h127.14.\h127.15.\h127.16.
\h127.17.\h127.18.\h127.19.\h127.20.\h127.21.\h127.22.\h127.23.\h127.24.
\h127.25.\h127.26.\h127.27.\h127.28.\h127.29.\h127.30.\h127.31.\h127.32.
\h127.33.\h127.34.\h127.35.\h127.36.\h127.37.\h127.38.\h127.39.\h127.40.
\h127.41.\h127.42.\h127.43.\h127.44.\h127.45.\h127.46.\h127.47.\h127.48.
\h127.49.\h127.50.\h127.51.\h127.52.\h127.53.\h127.54.\h127.55.\h127.56.
\h127.57.\h127.58.\h127.59.\h127.60.\h127.61.\h127.62.\h127.63.\h127.64.
\h127.65.\h127.66.\h127.67.\h127.68.\h127.69.\h127.70.\h127.71.\h127.73.
\h127.74.\h127.75.\h127.76.\h127.77.\h127.78.\h127.79.\h127.81.\h127.83.
\h127.84.\h127.85.\h127.87.\h127.88.\h127.89.\h127.91.\h127.92.\h127.93.
\h127.95.\h127.97.\h127.99.\h127.101.\h127.102.\h127.103.\h127.105.\h127.109.
\h127.111.\h127.112.\h127.115.\h127.117.\h127.119.\h127.121.\h127.127.\h128.2.
\h128.4.\h128.5.\h128.6.\h128.8.\h128.9.\h128.10.\h128.11.\h128.12.
\h128.13.\h128.14.\h128.15.\h128.16.\h128.17.\h128.18.\h128.19.\h128.20.
\h128.21.\h128.22.\h128.23.\h128.24.\h128.25.\h128.26.\h128.27.\h128.28.
\h128.29.\h128.30.\h128.31.\h128.32.\h128.33.\h128.34.\h128.35.\h128.36.
\h128.37.\h128.38.\h128.39.\h128.40.\h128.41.\h128.42.\h128.43.\h128.44.
\h128.45.\h128.46.\h128.47.\h128.48.\h128.49.\h128.50.\h128.51.\h128.52.
\h128.53.\h128.54.\h128.55.\h128.56.\h128.57.\h128.58.\h128.59.\h128.60.
\h128.61.\h128.62.\h128.63.\h128.64.\h128.65.\h128.66.\h128.67.\h128.68.
\h128.69.\h128.70.\h128.71.\h128.72.\h128.73.\h128.74.\h128.75.\h128.76.
\h128.77.\h128.78.\h128.79.\h128.80.\h128.82.\h128.83.\h128.84.\h128.86.
\h128.88.\h128.90.\h128.91.\h128.92.\h128.94.\h128.95.\h128.96.\h128.98.
\h128.100.\h128.101.\h128.102.\h128.103.\h128.104.\h128.106.\h128.107.\h128.108.
\h128.110.\h128.114.\h128.116.\h128.118.\h128.120.\h128.122.\h128.126.\h128.128.
\h129.5.\h129.9.\h129.10.\h129.11.\h129.12.\h129.13.\h129.14.\h129.15.
\h129.16.\h129.17.\h129.18.\h129.19.\h129.20.\h129.21.\h129.22.\h129.23.
\h129.24.\h129.25.\h129.26.\h129.27.\h129.28.\h129.29.\h129.30.\h129.31.
\h129.32.\h129.33.\h129.34.\h129.35.\h129.36.\h129.37.\h129.38.\h129.39.
\h129.40.\h129.41.\h129.42.\h129.43.\h129.44.\h129.45.\h129.46.\h129.47.
\h129.48.\h129.49.\h129.50.\h129.51.\h129.52.\h129.53.\h129.54.\h129.55.
\h129.56.\h129.57.\h129.58.\h129.59.\h129.60.\h129.61.\h129.62.\h129.63.
\h129.64.\h129.65.\h129.66.\h129.67.\h129.69.\h129.70.\h129.71.\h129.72.
\h129.73.\h129.74.\h129.75.\h129.76.\h129.77.\h129.78.\h129.79.\h129.81.
\h129.83.\h129.84.\h129.85.\h129.87.\h129.89.\h129.90.\h129.91.\h129.93.
\h129.96.\h129.97.\h129.99.\h129.101.\h129.102.\h129.103.\h129.105.\h129.107.
\h129.111.\h129.113.\h129.114.\h129.117.\h129.120.\h129.125.\h129.127.\h129.129.
\h130.4.\h130.7.\h130.8.\h130.9.\h130.10.\h130.11.\h130.12.\h130.13.
\h130.14.\h130.15.\h130.16.\h130.17.\h130.18.\h130.19.\h130.20.\h130.21.
\h130.22.\h130.23.\h130.24.\h130.25.\h130.26.\h130.27.\h130.28.\h130.29.
\h130.30.\h130.31.\h130.32.\h130.33.\h130.34.\h130.35.\h130.36.\h130.37.
\h130.38.\h130.39.\h130.40.\h130.41.\h130.42.\h130.43.\h130.44.\h130.45.
\h130.46.\h130.47.\h130.48.\h130.49.\h130.50.\h130.51.\h130.52.\h130.53.
\h130.54.\h130.55.\h130.56.\h130.57.\h130.58.\h130.59.\h130.60.\h130.61.
\h130.62.\h130.64.\h130.65.\h130.66.\h130.67.\h130.68.\h130.69.\h130.70.
\h130.71.\h130.72.\h130.73.\h130.74.\h130.76.\h130.78.\h130.79.\h130.80.
\h130.82.\h130.84.\h130.85.\h130.86.\h130.88.\h130.90.\h130.91.\h130.92.
\h130.94.\h130.96.\h130.98.\h130.100.\h130.102.\h130.103.\h130.108.\h130.110.
\h130.112.\h130.115.\h130.118.\h130.120.\h130.122.\h130.124.\h130.127.\h130.128.
\h130.130.\h131.3.\h131.5.\h131.7.\h131.8.\h131.9.\h131.10.\h131.11.
\h131.12.\h131.13.\h131.14.\h131.15.\h131.16.\h131.17.\h131.18.\h131.19.
\h131.20.\h131.21.\h131.22.\h131.23.\h131.24.\h131.25.\h131.26.\h131.27.
\h131.28.\h131.29.\h131.30.\h131.31.\h131.32.\h131.33.\h131.34.\h131.35.
\h131.36.\h131.37.\h131.38.\h131.39.\h131.40.\h131.41.\h131.42.\h131.43.
\h131.44.\h131.45.\h131.46.\h131.47.\h131.48.\h131.49.\h131.50.\h131.51.
\h131.52.\h131.53.\h131.54.\h131.55.\h131.56.\h131.57.\h131.58.\h131.59.
\h131.60.\h131.61.\h131.62.\h131.63.\h131.64.\h131.65.\h131.66.\h131.67.
\h131.68.\h131.69.\h131.70.\h131.71.\h131.73.\h131.74.\h131.75.\h131.77.
\h131.79.\h131.80.\h131.81.\h131.83.\h131.85.\h131.86.\h131.87.\h131.89.
\h131.90.\h131.91.\h131.93.\h131.95.\h131.97.\h131.98.\h131.99.\h131.101.
\h131.107.\h131.109.\h131.110.\h131.111.\h131.113.\h131.115.\h131.116.\h131.117.
\h131.119.\h131.125.\h131.127.\h131.128.\h131.131.\h132.2.\h132.6.\h132.8.
\h132.9.\h132.10.\h132.11.\h132.12.\h132.13.\h132.14.\h132.15.\h132.16.
\h132.17.\h132.18.\h132.19.\h132.20.\h132.21.\h132.22.\h132.23.\h132.24.
\h132.25.\h132.26.\h132.27.\h132.28.\h132.29.\h132.30.\h132.31.\h132.32.
\h132.33.\h132.34.\h132.35.\h132.36.\h132.37.\h132.38.\h132.39.\h132.40.
\h132.41.\h132.42.\h132.43.\h132.44.\h132.45.\h132.46.\h132.47.\h132.48.
\h132.49.\h132.50.\h132.51.\h132.52.\h132.53.\h132.54.\h132.55.\h132.56.
\h132.57.\h132.58.\h132.59.\h132.60.\h132.61.\h132.62.\h132.63.\h132.64.
\h132.65.\h132.66.\h132.67.\h132.68.\h132.69.\h132.70.\h132.72.\h132.73.
\h132.74.\h132.76.\h132.78.\h132.80.\h132.81.\h132.82.\h132.84.\h132.85.
\h132.86.\h132.87.\h132.88.\h132.90.\h132.92.\h132.94.\h132.96.\h132.97.
\h132.98.\h132.99.\h132.100.\h132.102.\h132.106.\h132.108.\h132.109.\h132.110.
\h132.111.\h132.112.\h132.114.\h132.117.\h132.118.\h132.120.\h132.122.\h132.124.
\h132.126.\h132.132.\h133.5.\h133.7.\h133.9.\h133.10.\h133.11.\h133.12.
\h133.13.\h133.14.\h133.15.\h133.16.\h133.17.\h133.18.\h133.19.\h133.20.
\h133.21.\h133.22.\h133.23.\h133.24.\h133.25.\h133.26.\h133.27.\h133.28.
\h133.29.\h133.30.\h133.31.\h133.32.\h133.33.\h133.34.\h133.35.\h133.36.
\h133.37.\h133.38.\h133.39.\h133.40.\h133.41.\h133.42.\h133.43.\h133.44.
\h133.45.\h133.46.\h133.47.\h133.48.\h133.49.\h133.50.\h133.51.\h133.52.
\h133.53.\h133.54.\h133.55.\h133.56.\h133.57.\h133.58.\h133.59.\h133.60.
\h133.61.\h133.63.\h133.64.\h133.65.\h133.66.\h133.67.\h133.68.\h133.69.
\h133.70.\h133.71.\h133.72.\h133.73.\h133.75.\h133.76.\h133.77.\h133.79.
\h133.80.\h133.81.\h133.82.\h133.83.\h133.84.\h133.85.\h133.87.\h133.88.
\h133.89.\h133.91.\h133.93.\h133.95.\h133.97.\h133.98.\h133.99.\h133.101.
\h133.103.\h133.105.\h133.106.\h133.107.\h133.109.\h133.111.\h133.113.\h133.115.
\h133.117.\h133.121.\h133.123.\h133.124.\h133.125.\h133.127.\h133.133.\h134.8.
\h134.9.\h134.10.\h134.11.\h134.12.\h134.13.\h134.14.\h134.15.\h134.16.
\h134.17.\h134.18.\h134.19.\h134.20.\h134.21.\h134.22.\h134.23.\h134.24.
\h134.25.\h134.26.\h134.27.\h134.28.\h134.29.\h134.30.\h134.31.\h134.32.
\h134.33.\h134.34.\h134.35.\h134.36.\h134.37.\h134.38.\h134.39.\h134.40.
\h134.41.\h134.42.\h134.43.\h134.44.\h134.45.\h134.46.\h134.47.\h134.48.
\h134.49.\h134.50.\h134.51.\h134.52.\h134.53.\h134.54.\h134.55.\h134.56.
\h134.57.\h134.58.\h134.59.\h134.60.\h134.61.\h134.62.\h134.63.\h134.64.
\h134.65.\h134.66.\h134.67.\h134.68.\h134.69.\h134.70.\h134.71.\h134.72.
\h134.74.\h134.76.\h134.77.\h134.78.\h134.80.\h134.81.\h134.82.\h134.83.
\h134.84.\h134.86.\h134.88.\h134.89.\h134.90.\h134.92.\h134.94.\h134.95.
\h134.96.\h134.98.\h134.100.\h134.104.\h134.106.\h134.108.\h134.110.\h134.112.
\h134.116.\h134.120.\h134.122.\h134.124.\h134.128.\h134.130.\h134.132.\h135.5.
\h135.7.\h135.9.\h135.10.\h135.11.\h135.12.\h135.13.\h135.14.\h135.15.
\h135.16.\h135.17.\h135.18.\h135.19.\h135.20.\h135.21.\h135.22.\h135.23.
\h135.24.\h135.25.\h135.26.\h135.27.\h135.28.\h135.29.\h135.30.\h135.31.
\h135.32.\h135.33.\h135.34.\h135.35.\h135.36.\h135.37.\h135.38.\h135.39.
\h135.40.\h135.41.\h135.42.\h135.43.\h135.44.\h135.45.\h135.46.\h135.47.
\h135.48.\h135.49.\h135.50.\h135.51.\h135.52.\h135.53.\h135.54.\h135.55.
\h135.56.\h135.57.\h135.58.\h135.59.\h135.60.\h135.61.\h135.63.\h135.65.
\h135.66.\h135.67.\h135.68.\h135.69.\h135.71.\h135.73.\h135.75.\h135.76.
\h135.77.\h135.78.\h135.79.\h135.81.\h135.83.\h135.85.\h135.87.\h135.89.
\h135.91.\h135.93.\h135.95.\h135.99.\h135.101.\h135.103.\h135.105.\h135.107.
\h135.109.\h135.111.\h135.114.\h135.115.\h135.117.\h135.119.\h135.121.\h135.123.
\h135.125.\h135.129.\h135.131.\h135.135.\h136.4.\h136.6.\h136.8.\h136.9.
\h136.10.\h136.11.\h136.12.\h136.13.\h136.14.\h136.15.\h136.16.\h136.17.
\h136.18.\h136.19.\h136.20.\h136.21.\h136.22.\h136.23.\h136.24.\h136.25.
\h136.26.\h136.27.\h136.28.\h136.29.\h136.30.\h136.31.\h136.32.\h136.33.
\h136.34.\h136.35.\h136.36.\h136.37.\h136.38.\h136.39.\h136.40.\h136.41.
\h136.42.\h136.43.\h136.44.\h136.45.\h136.46.\h136.47.\h136.48.\h136.49.
\h136.50.\h136.51.\h136.52.\h136.53.\h136.54.\h136.55.\h136.56.\h136.57.
\h136.58.\h136.59.\h136.60.\h136.61.\h136.62.\h136.64.\h136.66.\h136.67.
\h136.68.\h136.70.\h136.72.\h136.74.\h136.75.\h136.76.\h136.78.\h136.80.
\h136.82.\h136.85.\h136.86.\h136.88.\h136.90.\h136.92.\h136.94.\h136.96.
\h136.98.\h136.100.\h136.102.\h136.103.\h136.104.\h136.106.\h136.112.\h136.114.
\h136.118.\h136.120.\h136.122.\h136.124.\h136.130.\h136.132.\h136.134.\h137.5.
\h137.7.\h137.8.\h137.9.\h137.10.\h137.11.\h137.12.\h137.13.\h137.14.
\h137.15.\h137.16.\h137.17.\h137.18.\h137.19.\h137.20.\h137.21.\h137.22.
\h137.23.\h137.24.\h137.25.\h137.26.\h137.27.\h137.28.\h137.29.\h137.30.
\h137.31.\h137.32.\h137.33.\h137.34.\h137.35.\h137.36.\h137.37.\h137.38.
\h137.39.\h137.40.\h137.41.\h137.42.\h137.43.\h137.44.\h137.45.\h137.46.
\h137.47.\h137.48.\h137.49.\h137.50.\h137.51.\h137.52.\h137.53.\h137.54.
\h137.55.\h137.56.\h137.57.\h137.58.\h137.59.\h137.60.\h137.61.\h137.62.
\h137.63.\h137.64.\h137.65.\h137.66.\h137.67.\h137.68.\h137.69.\h137.71.
\h137.73.\h137.74.\h137.75.\h137.77.\h137.79.\h137.81.\h137.83.\h137.85.
\h137.87.\h137.89.\h137.91.\h137.92.\h137.93.\h137.95.\h137.97.\h137.98.
\h137.99.\h137.101.\h137.103.\h137.105.\h137.109.\h137.111.\h137.113.\h137.115.
\h137.117.\h137.119.\h137.121.\h137.125.\h137.127.\h137.128.\h137.131.\h138.6.
\h138.8.\h138.9.\h138.10.\h138.11.\h138.12.\h138.13.\h138.14.\h138.15.
\h138.16.\h138.17.\h138.18.\h138.19.\h138.20.\h138.21.\h138.22.\h138.23.
\h138.24.\h138.25.\h138.26.\h138.27.\h138.28.\h138.29.\h138.30.\h138.31.
\h138.32.\h138.33.\h138.34.\h138.35.\h138.36.\h138.37.\h138.38.\h138.39.
\h138.40.\h138.41.\h138.42.\h138.43.\h138.44.\h138.45.\h138.46.\h138.47.
\h138.48.\h138.49.\h138.50.\h138.51.\h138.52.\h138.53.\h138.54.\h138.56.
\h138.57.\h138.58.\h138.59.\h138.60.\h138.62.\h138.63.\h138.64.\h138.65.
\h138.66.\h138.67.\h138.68.\h138.69.\h138.70.\h138.72.\h138.74.\h138.75.
\h138.76.\h138.77.\h138.78.\h138.80.\h138.81.\h138.82.\h138.84.\h138.86.
\h138.87.\h138.88.\h138.90.\h138.92.\h138.93.\h138.94.\h138.96.\h138.98.
\h138.99.\h138.100.\h138.102.\h138.104.\h138.106.\h138.108.\h138.110.\h138.113.
\h138.114.\h138.116.\h138.120.\h138.124.\h138.126.\h138.132.\h139.5.\h139.7.
\h139.9.\h139.10.\h139.11.\h139.12.\h139.13.\h139.14.\h139.15.\h139.16.
\h139.17.\h139.18.\h139.19.\h139.20.\h139.21.\h139.22.\h139.23.\h139.24.
\h139.25.\h139.26.\h139.27.\h139.28.\h139.29.\h139.30.\h139.31.\h139.32.
\h139.33.\h139.34.\h139.35.\h139.36.\h139.37.\h139.38.\h139.39.\h139.40.
\h139.41.\h139.42.\h139.43.\h139.44.\h139.45.\h139.46.\h139.47.\h139.48.
\h139.49.\h139.50.\h139.51.\h139.52.\h139.53.\h139.54.\h139.55.\h139.56.
\h139.57.\h139.58.\h139.59.\h139.61.\h139.63.\h139.64.\h139.65.\h139.66.
\h139.67.\h139.69.\h139.70.\h139.71.\h139.73.\h139.75.\h139.76.\h139.77.
\h139.79.\h139.81.\h139.83.\h139.84.\h139.85.\h139.87.\h139.89.\h139.91.
\h139.93.\h139.95.\h139.97.\h139.99.\h139.101.\h139.102.\h139.103.\h139.105.
\h139.107.\h139.109.\h139.111.\h139.112.\h139.115.\h139.117.\h139.119.\h139.121.
\h139.123.\h139.125.\h139.127.\h139.131.\h139.139.\h140.6.\h140.7.\h140.8.
\h140.10.\h140.11.\h140.12.\h140.13.\h140.14.\h140.15.\h140.16.\h140.17.
\h140.18.\h140.19.\h140.20.\h140.21.\h140.22.\h140.23.\h140.24.\h140.25.
\h140.26.\h140.27.\h140.28.\h140.29.\h140.30.\h140.31.\h140.32.\h140.33.
\h140.34.\h140.35.\h140.36.\h140.37.\h140.38.\h140.39.\h140.40.\h140.41.
\h140.42.\h140.43.\h140.44.\h140.45.\h140.46.\h140.47.\h140.48.\h140.49.
\h140.50.\h140.51.\h140.52.\h140.53.\h140.54.\h140.55.\h140.56.\h140.58.
\h140.59.\h140.60.\h140.61.\h140.62.\h140.63.\h140.64.\h140.65.\h140.66.
\h140.68.\h140.70.\h140.72.\h140.73.\h140.74.\h140.76.\h140.77.\h140.78.
\h140.80.\h140.82.\h140.83.\h140.84.\h140.86.\h140.88.\h140.90.\h140.91.
\h140.92.\h140.94.\h140.96.\h140.98.\h140.100.\h140.101.\h140.102.\h140.104.
\h140.110.\h140.114.\h140.116.\h140.118.\h140.119.\h140.122.\h140.126.\h140.130.
\h140.134.\h141.3.\h141.5.\h141.7.\h141.9.\h141.10.\h141.11.\h141.12.
\h141.13.\h141.14.\h141.15.\h141.16.\h141.17.\h141.18.\h141.19.\h141.20.
\h141.21.\h141.22.\h141.23.\h141.24.\h141.25.\h141.26.\h141.27.\h141.28.
\h141.29.\h141.30.\h141.31.\h141.32.\h141.33.\h141.34.\h141.35.\h141.36.
\h141.37.\h141.38.\h141.39.\h141.40.\h141.41.\h141.42.\h141.43.\h141.44.
\h141.45.\h141.46.\h141.47.\h141.48.\h141.49.\h141.50.\h141.51.\h141.53.
\h141.54.\h141.55.\h141.56.\h141.57.\h141.59.\h141.60.\h141.61.\h141.62.
\h141.63.\h141.64.\h141.65.\h141.66.\h141.67.\h141.68.\h141.69.\h141.71.
\h141.72.\h141.73.\h141.74.\h141.75.\h141.77.\h141.78.\h141.79.\h141.81.
\h141.83.\h141.85.\h141.87.\h141.89.\h141.90.\h141.91.\h141.93.\h141.97.
\h141.99.\h141.101.\h141.102.\h141.103.\h141.105.\h141.108.\h141.109.\h141.111.
\h141.113.\h141.116.\h141.117.\h141.121.\h141.123.\h141.129.\h141.133.\h141.141.
\h142.4.\h142.6.\h142.7.\h142.8.\h142.10.\h142.11.\h142.12.\h142.13.
\h142.14.\h142.15.\h142.16.\h142.17.\h142.18.\h142.19.\h142.20.\h142.21.
\h142.22.\h142.23.\h142.24.\h142.25.\h142.26.\h142.27.\h142.28.\h142.29.
\h142.30.\h142.31.\h142.32.\h142.33.\h142.34.\h142.35.\h142.36.\h142.37.
\h142.38.\h142.39.\h142.40.\h142.41.\h142.42.\h142.43.\h142.44.\h142.45.
\h142.46.\h142.47.\h142.48.\h142.49.\h142.50.\h142.52.\h142.53.\h142.54.
\h142.55.\h142.56.\h142.58.\h142.59.\h142.60.\h142.61.\h142.62.\h142.64.
\h142.65.\h142.66.\h142.68.\h142.69.\h142.70.\h142.72.\h142.73.\h142.74.
\h142.76.\h142.78.\h142.80.\h142.82.\h142.84.\h142.86.\h142.87.\h142.88.
\h142.90.\h142.92.\h142.94.\h142.96.\h142.97.\h142.100.\h142.104.\h142.106.
\h142.108.\h142.112.\h142.114.\h142.115.\h142.116.\h142.120.\h142.124.\h142.128.
\h142.130.\h142.132.\h142.136.\h142.142.\h143.5.\h143.7.\h143.8.\h143.9.
\h143.10.\h143.11.\h143.13.\h143.14.\h143.15.\h143.16.\h143.17.\h143.18.
\h143.19.\h143.20.\h143.21.\h143.22.\h143.23.\h143.24.\h143.25.\h143.26.
\h143.27.\h143.28.\h143.29.\h143.30.\h143.31.\h143.32.\h143.33.\h143.34.
\h143.35.\h143.36.\h143.37.\h143.38.\h143.39.\h143.40.\h143.41.\h143.42.
\h143.43.\h143.44.\h143.45.\h143.46.\h143.47.\h143.48.\h143.49.\h143.50.
\h143.51.\h143.53.\h143.54.\h143.55.\h143.56.\h143.57.\h143.58.\h143.59.
\h143.60.\h143.61.\h143.63.\h143.65.\h143.67.\h143.69.\h143.70.\h143.71.
\h143.73.\h143.75.\h143.77.\h143.79.\h143.81.\h143.83.\h143.85.\h143.86.
\h143.87.\h143.88.\h143.89.\h143.91.\h143.93.\h143.95.\h143.97.\h143.99.
\h143.101.\h143.103.\h143.107.\h143.111.\h143.113.\h143.115.\h143.116.\h143.117.
\h143.119.\h143.123.\h143.125.\h143.127.\h143.128.\h143.131.\h143.135.\h143.141.
\h143.143.\h144.2.\h144.4.\h144.6.\h144.8.\h144.9.\h144.10.\h144.11.
\h144.12.\h144.13.\h144.14.\h144.15.\h144.16.\h144.17.\h144.18.\h144.19.
\h144.20.\h144.21.\h144.22.\h144.23.\h144.24.\h144.25.\h144.26.\h144.27.
\h144.28.\h144.29.\h144.30.\h144.31.\h144.32.\h144.33.\h144.34.\h144.35.
\h144.36.\h144.37.\h144.38.\h144.39.\h144.40.\h144.41.\h144.42.\h144.43.
\h144.44.\h144.45.\h144.46.\h144.47.\h144.48.\h144.49.\h144.50.\h144.51.
\h144.52.\h144.53.\h144.54.\h144.56.\h144.57.\h144.58.\h144.59.\h144.60.
\h144.62.\h144.63.\h144.64.\h144.66.\h144.67.\h144.68.\h144.69.\h144.70.
\h144.72.\h144.74.\h144.75.\h144.76.\h144.78.\h144.80.\h144.81.\h144.82.
\h144.84.\h144.86.\h144.87.\h144.88.\h144.90.\h144.92.\h144.94.\h144.96.
\h144.98.\h144.102.\h144.104.\h144.106.\h144.108.\h144.110.\h144.111.\h144.112.
\h144.114.\h144.116.\h144.120.\h144.122.\h144.124.\h144.126.\h144.129.\h144.130.
\h144.132.\h144.142.\h145.1.\h145.5.\h145.7.\h145.9.\h145.10.\h145.11.
\h145.12.\h145.13.\h145.14.\h145.15.\h145.16.\h145.17.\h145.18.\h145.19.
\h145.20.\h145.21.\h145.22.\h145.23.\h145.24.\h145.25.\h145.26.\h145.27.
\h145.28.\h145.29.\h145.30.\h145.31.\h145.32.\h145.33.\h145.34.\h145.35.
\h145.36.\h145.37.\h145.38.\h145.39.\h145.40.\h145.41.\h145.42.\h145.43.
\h145.44.\h145.45.\h145.46.\h145.47.\h145.48.\h145.49.\h145.50.\h145.51.
\h145.52.\h145.53.\h145.55.\h145.56.\h145.57.\h145.58.\h145.59.\h145.60.
\h145.61.\h145.63.\h145.64.\h145.65.\h145.66.\h145.67.\h145.69.\h145.70.
\h145.71.\h145.73.\h145.74.\h145.75.\h145.76.\h145.77.\h145.79.\h145.81.
\h145.82.\h145.83.\h145.85.\h145.87.\h145.89.\h145.91.\h145.93.\h145.94.
\h145.95.\h145.97.\h145.99.\h145.100.\h145.101.\h145.103.\h145.105.\h145.106.
\h145.107.\h145.109.\h145.111.\h145.112.\h145.113.\h145.115.\h145.121.\h145.125.
\h145.127.\h145.130.\h145.141.\h146.8.\h146.9.\h146.10.\h146.11.\h146.12.
\h146.13.\h146.14.\h146.15.\h146.16.\h146.17.\h146.18.\h146.19.\h146.20.
\h146.21.\h146.22.\h146.23.\h146.24.\h146.25.\h146.26.\h146.27.\h146.28.
\h146.29.\h146.30.\h146.31.\h146.32.\h146.33.\h146.34.\h146.35.\h146.36.
\h146.37.\h146.38.\h146.39.\h146.40.\h146.41.\h146.42.\h146.43.\h146.44.
\h146.46.\h146.47.\h146.48.\h146.49.\h146.50.\h146.51.\h146.52.\h146.53.
\h146.54.\h146.56.\h146.58.\h146.59.\h146.60.\h146.62.\h146.63.\h146.64.
\h146.65.\h146.66.\h146.68.\h146.69.\h146.70.\h146.71.\h146.72.\h146.74.
\h146.75.\h146.76.\h146.77.\h146.78.\h146.80.\h146.81.\h146.82.\h146.83.
\h146.84.\h146.86.\h146.88.\h146.89.\h146.90.\h146.92.\h146.94.\h146.95.
\h146.96.\h146.98.\h146.100.\h146.101.\h146.102.\h146.104.\h146.105.\h146.106.
\h146.108.\h146.110.\h146.112.\h146.113.\h146.116.\h146.124.\h146.126.\h146.128.
\h146.140.\h147.7.\h147.9.\h147.11.\h147.12.\h147.13.\h147.14.\h147.15.
\h147.16.\h147.17.\h147.18.\h147.19.\h147.20.\h147.21.\h147.22.\h147.23.
\h147.24.\h147.25.\h147.26.\h147.27.\h147.28.\h147.29.\h147.30.\h147.31.
\h147.32.\h147.33.\h147.34.\h147.35.\h147.36.\h147.37.\h147.38.\h147.39.
\h147.40.\h147.41.\h147.42.\h147.43.\h147.44.\h147.45.\h147.46.\h147.47.
\h147.48.\h147.49.\h147.50.\h147.51.\h147.52.\h147.53.\h147.54.\h147.55.
\h147.56.\h147.57.\h147.58.\h147.59.\h147.60.\h147.61.\h147.62.\h147.63.
\h147.64.\h147.65.\h147.66.\h147.67.\h147.69.\h147.70.\h147.71.\h147.72.
\h147.73.\h147.74.\h147.75.\h147.76.\h147.77.\h147.79.\h147.81.\h147.82.
\h147.83.\h147.84.\h147.85.\h147.87.\h147.89.\h147.91.\h147.93.\h147.94.
\h147.95.\h147.96.\h147.97.\h147.99.\h147.102.\h147.103.\h147.105.\h147.107.
\h147.109.\h147.111.\h147.119.\h147.123.\h147.125.\h147.127.\h147.139.\h147.147.
\h148.4.\h148.6.\h148.8.\h148.9.\h148.10.\h148.11.\h148.12.\h148.13.
\h148.14.\h148.15.\h148.16.\h148.17.\h148.18.\h148.19.\h148.20.\h148.21.
\h148.22.\h148.23.\h148.24.\h148.25.\h148.26.\h148.27.\h148.28.\h148.29.
\h148.30.\h148.31.\h148.32.\h148.33.\h148.34.\h148.35.\h148.36.\h148.37.
\h148.38.\h148.39.\h148.40.\h148.41.\h148.42.\h148.43.\h148.44.\h148.45.
\h148.46.\h148.47.\h148.48.\h148.49.\h148.50.\h148.51.\h148.52.\h148.53.
\h148.54.\h148.55.\h148.56.\h148.57.\h148.58.\h148.60.\h148.61.\h148.62.
\h148.63.\h148.64.\h148.65.\h148.66.\h148.67.\h148.68.\h148.69.\h148.70.
\h148.72.\h148.73.\h148.74.\h148.75.\h148.76.\h148.78.\h148.80.\h148.82.
\h148.83.\h148.84.\h148.86.\h148.88.\h148.90.\h148.91.\h148.92.\h148.94.
\h148.96.\h148.98.\h148.100.\h148.102.\h148.104.\h148.106.\h148.108.\h148.109.
\h148.110.\h148.118.\h148.122.\h148.124.\h148.138.\h148.142.\h148.148.\h149.1.
\h149.5.\h149.7.\h149.8.\h149.9.\h149.10.\h149.11.\h149.12.\h149.13.
\h149.14.\h149.15.\h149.16.\h149.17.\h149.18.\h149.19.\h149.20.\h149.21.
\h149.22.\h149.23.\h149.24.\h149.25.\h149.26.\h149.27.\h149.28.\h149.29.
\h149.30.\h149.31.\h149.32.\h149.33.\h149.34.\h149.35.\h149.36.\h149.37.
\h149.38.\h149.39.\h149.40.\h149.41.\h149.42.\h149.43.\h149.44.\h149.45.
\h149.46.\h149.47.\h149.48.\h149.49.\h149.50.\h149.51.\h149.52.\h149.53.
\h149.54.\h149.55.\h149.56.\h149.57.\h149.59.\h149.61.\h149.62.\h149.63.
\h149.65.\h149.66.\h149.67.\h149.68.\h149.69.\h149.71.\h149.73.\h149.74.
\h149.75.\h149.77.\h149.79.\h149.80.\h149.81.\h149.83.\h149.85.\h149.86.
\h149.87.\h149.89.\h149.91.\h149.92.\h149.93.\h149.95.\h149.97.\h149.101.
\h149.105.\h149.107.\h149.109.\h149.113.\h149.117.\h149.119.\h149.121.\h149.125.
\h149.137.\h149.149.\h150.6.\h150.8.\h150.9.\h150.10.\h150.11.\h150.12.
\h150.13.\h150.14.\h150.15.\h150.16.\h150.17.\h150.18.\h150.19.\h150.20.
\h150.21.\h150.22.\h150.23.\h150.24.\h150.25.\h150.26.\h150.27.\h150.28.
\h150.29.\h150.30.\h150.32.\h150.33.\h150.34.\h150.35.\h150.36.\h150.37.
\h150.38.\h150.39.\h150.40.\h150.41.\h150.42.\h150.43.\h150.44.\h150.45.
\h150.46.\h150.48.\h150.50.\h150.51.\h150.52.\h150.53.\h150.54.\h150.55.
\h150.56.\h150.57.\h150.58.\h150.60.\h150.61.\h150.62.\h150.63.\h150.64.
\h150.65.\h150.66.\h150.68.\h150.69.\h150.70.\h150.72.\h150.74.\h150.76.
\h150.78.\h150.80.\h150.81.\h150.82.\h150.84.\h150.86.\h150.88.\h150.90.
\h150.92.\h150.96.\h150.100.\h150.102.\h150.104.\h150.108.\h150.114.\h150.116.
\h150.120.\h150.124.\h150.126.\h150.132.\h150.140.\h150.150.\h151.7.\h151.8.
\h151.9.\h151.10.\h151.11.\h151.12.\h151.13.\h151.14.\h151.15.\h151.16.
\h151.17.\h151.18.\h151.19.\h151.20.\h151.21.\h151.22.\h151.23.\h151.24.
\h151.25.\h151.26.\h151.27.\h151.28.\h151.29.\h151.30.\h151.31.\h151.32.
\h151.33.\h151.34.\h151.35.\h151.36.\h151.37.\h151.38.\h151.39.\h151.40.
\h151.41.\h151.42.\h151.43.\h151.44.\h151.45.\h151.46.\h151.47.\h151.49.
\h151.50.\h151.51.\h151.52.\h151.53.\h151.54.\h151.55.\h151.56.\h151.57.
\h151.59.\h151.60.\h151.61.\h151.63.\h151.64.\h151.65.\h151.67.\h151.69.
\h151.70.\h151.71.\h151.73.\h151.75.\h151.77.\h151.79.\h151.83.\h151.85.
\h151.87.\h151.88.\h151.91.\h151.97.\h151.103.\h151.107.\h151.109.\h151.111.
\h151.115.\h151.119.\h151.121.\h151.123.\h151.127.\h151.135.\h151.139.\h151.151.
\h152.6.\h152.8.\h152.10.\h152.11.\h152.12.\h152.13.\h152.14.\h152.15.
\h152.16.\h152.17.\h152.18.\h152.19.\h152.20.\h152.21.\h152.22.\h152.23.
\h152.24.\h152.25.\h152.26.\h152.27.\h152.28.\h152.29.\h152.30.\h152.31.
\h152.32.\h152.33.\h152.34.\h152.35.\h152.36.\h152.37.\h152.38.\h152.39.
\h152.40.\h152.41.\h152.42.\h152.44.\h152.45.\h152.46.\h152.47.\h152.48.
\h152.49.\h152.50.\h152.52.\h152.53.\h152.54.\h152.55.\h152.56.\h152.58.
\h152.59.\h152.60.\h152.62.\h152.64.\h152.65.\h152.66.\h152.67.\h152.68.
\h152.70.\h152.72.\h152.74.\h152.76.\h152.77.\h152.78.\h152.80.\h152.82.
\h152.84.\h152.86.\h152.90.\h152.92.\h152.94.\h152.98.\h152.100.\h152.102.
\h152.106.\h152.110.\h152.116.\h152.118.\h152.122.\h152.134.\h152.142.\h152.146.
\h152.150.\h153.5.\h153.6.\h153.7.\h153.9.\h153.10.\h153.11.\h153.12.
\h153.13.\h153.14.\h153.15.\h153.16.\h153.17.\h153.18.\h153.19.\h153.20.
\h153.21.\h153.22.\h153.23.\h153.24.\h153.25.\h153.26.\h153.27.\h153.28.
\h153.29.\h153.30.\h153.31.\h153.32.\h153.33.\h153.34.\h153.35.\h153.36.
\h153.37.\h153.38.\h153.39.\h153.40.\h153.41.\h153.42.\h153.43.\h153.44.
\h153.45.\h153.46.\h153.47.\h153.48.\h153.49.\h153.50.\h153.51.\h153.53.
\h153.54.\h153.55.\h153.56.\h153.57.\h153.59.\h153.60.\h153.61.\h153.62.
\h153.63.\h153.65.\h153.66.\h153.67.\h153.69.\h153.71.\h153.73.\h153.75.
\h153.77.\h153.79.\h153.81.\h153.83.\h153.85.\h153.87.\h153.89.\h153.93.
\h153.95.\h153.99.\h153.101.\h153.105.\h153.113.\h153.117.\h153.121.\h153.131.
\h153.133.\h153.141.\h153.145.\h153.147.\h154.4.\h154.7.\h154.8.\h154.9.
\h154.10.\h154.12.\h154.13.\h154.14.\h154.15.\h154.16.\h154.17.\h154.18.
\h154.19.\h154.20.\h154.21.\h154.22.\h154.23.\h154.24.\h154.25.\h154.26.
\h154.27.\h154.28.\h154.29.\h154.30.\h154.31.\h154.32.\h154.33.\h154.34.
\h154.35.\h154.36.\h154.37.\h154.38.\h154.39.\h154.40.\h154.41.\h154.42.
\h154.43.\h154.44.\h154.45.\h154.46.\h154.47.\h154.48.\h154.49.\h154.50.
\h154.51.\h154.52.\h154.54.\h154.55.\h154.56.\h154.58.\h154.60.\h154.61.
\h154.62.\h154.63.\h154.64.\h154.66.\h154.68.\h154.70.\h154.72.\h154.73.
\h154.74.\h154.76.\h154.78.\h154.80.\h154.82.\h154.84.\h154.86.\h154.88.
\h154.92.\h154.94.\h154.96.\h154.100.\h154.102.\h154.104.\h154.112.\h154.115.
\h154.116.\h154.118.\h154.120.\h154.130.\h154.142.\h154.144.\h155.7.\h155.8.
\h155.9.\h155.10.\h155.11.\h155.12.\h155.13.\h155.14.\h155.15.\h155.16.
\h155.17.\h155.18.\h155.19.\h155.20.\h155.21.\h155.22.\h155.23.\h155.24.
\h155.25.\h155.26.\h155.27.\h155.28.\h155.29.\h155.30.\h155.31.\h155.32.
\h155.33.\h155.34.\h155.35.\h155.36.\h155.37.\h155.38.\h155.39.\h155.40.
\h155.41.\h155.42.\h155.43.\h155.44.\h155.45.\h155.47.\h155.48.\h155.49.
\h155.50.\h155.51.\h155.52.\h155.53.\h155.55.\h155.56.\h155.57.\h155.58.
\h155.59.\h155.61.\h155.62.\h155.63.\h155.65.\h155.67.\h155.68.\h155.69.
\h155.71.\h155.73.\h155.75.\h155.77.\h155.79.\h155.80.\h155.81.\h155.83.
\h155.85.\h155.86.\h155.87.\h155.89.\h155.91.\h155.95.\h155.97.\h155.98.
\h155.99.\h155.101.\h155.103.\h155.105.\h155.111.\h155.112.\h155.113.\h155.115.
\h155.117.\h155.119.\h155.121.\h155.137.\h155.149.\h155.155.\h156.6.\h156.8.
\h156.9.\h156.10.\h156.11.\h156.12.\h156.13.\h156.14.\h156.15.\h156.16.
\h156.17.\h156.18.\h156.20.\h156.21.\h156.22.\h156.23.\h156.24.\h156.25.
\h156.26.\h156.27.\h156.28.\h156.29.\h156.30.\h156.31.\h156.32.\h156.33.
\h156.34.\h156.35.\h156.36.\h156.37.\h156.38.\h156.39.\h156.40.\h156.41.
\h156.42.\h156.43.\h156.44.\h156.45.\h156.46.\h156.48.\h156.49.\h156.50.
\h156.51.\h156.52.\h156.53.\h156.54.\h156.55.\h156.56.\h156.57.\h156.58.
\h156.59.\h156.60.\h156.62.\h156.63.\h156.64.\h156.66.\h156.67.\h156.68.
\h156.69.\h156.70.\h156.72.\h156.74.\h156.75.\h156.76.\h156.78.\h156.80.
\h156.82.\h156.83.\h156.84.\h156.86.\h156.87.\h156.90.\h156.94.\h156.96.
\h156.98.\h156.100.\h156.101.\h156.102.\h156.104.\h156.108.\h156.110.\h156.111.
\h156.114.\h156.120.\h156.126.\h156.144.\h156.146.\h157.7.\h157.9.\h157.10.
\h157.11.\h157.12.\h157.13.\h157.14.\h157.15.\h157.16.\h157.17.\h157.18.
\h157.19.\h157.20.\h157.21.\h157.22.\h157.23.\h157.24.\h157.25.\h157.26.
\h157.27.\h157.28.\h157.29.\h157.30.\h157.31.\h157.32.\h157.33.\h157.34.
\h157.35.\h157.37.\h157.38.\h157.39.\h157.40.\h157.41.\h157.42.\h157.43.
\h157.44.\h157.45.\h157.46.\h157.47.\h157.48.\h157.49.\h157.51.\h157.52.
\h157.53.\h157.54.\h157.55.\h157.56.\h157.57.\h157.58.\h157.59.\h157.61.
\h157.63.\h157.64.\h157.65.\h157.67.\h157.69.\h157.71.\h157.72.\h157.73.
\h157.75.\h157.76.\h157.77.\h157.79.\h157.81.\h157.82.\h157.83.\h157.85.
\h157.87.\h157.90.\h157.91.\h157.93.\h157.95.\h157.97.\h157.99.\h157.100.
\h157.103.\h157.109.\h157.113.\h157.115.\h157.117.\h157.119.\h157.141.\h157.145.
\h158.8.\h158.9.\h158.10.\h158.11.\h158.12.\h158.13.\h158.14.\h158.16.
\h158.17.\h158.18.\h158.19.\h158.20.\h158.21.\h158.22.\h158.23.\h158.24.
\h158.25.\h158.26.\h158.27.\h158.28.\h158.29.\h158.30.\h158.31.\h158.32.
\h158.33.\h158.34.\h158.35.\h158.36.\h158.37.\h158.38.\h158.39.\h158.40.
\h158.41.\h158.42.\h158.43.\h158.44.\h158.45.\h158.46.\h158.47.\h158.48.
\h158.50.\h158.52.\h158.53.\h158.54.\h158.56.\h158.57.\h158.58.\h158.59.
\h158.60.\h158.61.\h158.62.\h158.63.\h158.64.\h158.65.\h158.66.\h158.68.
\h158.70.\h158.71.\h158.72.\h158.73.\h158.74.\h158.76.\h158.77.\h158.80.
\h158.82.\h158.84.\h158.86.\h158.88.\h158.89.\h158.92.\h158.96.\h158.98.
\h158.100.\h158.101.\h158.102.\h158.110.\h158.112.\h158.116.\h158.119.\h158.140.
\h158.143.\h159.7.\h159.9.\h159.10.\h159.11.\h159.12.\h159.13.\h159.15.
\h159.16.\h159.17.\h159.18.\h159.19.\h159.20.\h159.21.\h159.22.\h159.23.
\h159.24.\h159.25.\h159.26.\h159.27.\h159.28.\h159.29.\h159.30.\h159.31.
\h159.32.\h159.33.\h159.34.\h159.35.\h159.36.\h159.37.\h159.38.\h159.39.
\h159.41.\h159.42.\h159.43.\h159.44.\h159.45.\h159.46.\h159.47.\h159.48.
\h159.49.\h159.51.\h159.53.\h159.54.\h159.55.\h159.57.\h159.59.\h159.60.
\h159.61.\h159.63.\h159.64.\h159.65.\h159.67.\h159.69.\h159.71.\h159.72.
\h159.73.\h159.75.\h159.79.\h159.81.\h159.83.\h159.85.\h159.87.\h159.90.
\h159.91.\h159.95.\h159.99.\h159.102.\h159.111.\h159.114.\h159.115.\h159.117.
\h159.123.\h159.141.\h159.144.\h160.6.\h160.7.\h160.8.\h160.10.\h160.12.
\h160.13.\h160.14.\h160.15.\h160.16.\h160.17.\h160.18.\h160.19.\h160.20.
\h160.21.\h160.22.\h160.23.\h160.24.\h160.25.\h160.26.\h160.27.\h160.28.
\h160.29.\h160.30.\h160.31.\h160.32.\h160.33.\h160.34.\h160.35.\h160.36.
\h160.37.\h160.38.\h160.39.\h160.40.\h160.41.\h160.42.\h160.43.\h160.44.
\h160.45.\h160.46.\h160.47.\h160.48.\h160.50.\h160.52.\h160.54.\h160.55.
\h160.56.\h160.57.\h160.58.\h160.60.\h160.62.\h160.64.\h160.66.\h160.68.
\h160.70.\h160.72.\h160.73.\h160.74.\h160.76.\h160.78.\h160.79.\h160.80.
\h160.82.\h160.85.\h160.86.\h160.88.\h160.90.\h160.94.\h160.97.\h160.98.
\h160.100.\h160.106.\h160.109.\h160.112.\h160.114.\h160.115.\h160.116.\h160.118.
\h160.136.\h160.142.\h160.160.\h161.5.\h161.7.\h161.8.\h161.9.\h161.11.
\h161.12.\h161.13.\h161.14.\h161.15.\h161.16.\h161.17.\h161.18.\h161.19.
\h161.20.\h161.21.\h161.22.\h161.23.\h161.24.\h161.25.\h161.26.\h161.27.
\h161.28.\h161.29.\h161.30.\h161.31.\h161.32.\h161.33.\h161.34.\h161.35.
\h161.36.\h161.37.\h161.38.\h161.39.\h161.40.\h161.41.\h161.42.\h161.43.
\h161.44.\h161.45.\h161.46.\h161.47.\h161.49.\h161.51.\h161.52.\h161.53.
\h161.54.\h161.55.\h161.56.\h161.57.\h161.59.\h161.61.\h161.62.\h161.63.
\h161.65.\h161.66.\h161.67.\h161.68.\h161.69.\h161.71.\h161.72.\h161.73.
\h161.75.\h161.77.\h161.80.\h161.81.\h161.83.\h161.84.\h161.85.\h161.86.
\h161.87.\h161.89.\h161.93.\h161.95.\h161.97.\h161.98.\h161.99.\h161.101.
\h161.107.\h161.113.\h161.117.\h161.131.\h161.141.\h161.149.\h162.4.\h162.6.
\h162.8.\h162.9.\h162.10.\h162.11.\h162.12.\h162.14.\h162.15.\h162.16.
\h162.17.\h162.18.\h162.19.\h162.20.\h162.21.\h162.22.\h162.24.\h162.25.
\h162.26.\h162.27.\h162.28.\h162.29.\h162.30.\h162.31.\h162.32.\h162.33.
\h162.34.\h162.35.\h162.36.\h162.37.\h162.38.\h162.39.\h162.40.\h162.41.
\h162.42.\h162.44.\h162.45.\h162.46.\h162.48.\h162.49.\h162.50.\h162.51.
\h162.52.\h162.53.\h162.54.\h162.55.\h162.56.\h162.57.\h162.58.\h162.60.
\h162.62.\h162.63.\h162.64.\h162.65.\h162.66.\h162.67.\h162.68.\h162.69.
\h162.70.\h162.72.\h162.74.\h162.76.\h162.78.\h162.80.\h162.81.\h162.82.
\h162.84.\h162.86.\h162.87.\h162.90.\h162.94.\h162.96.\h162.100.\h162.102.
\h162.108.\h162.112.\h162.116.\h162.120.\h162.140.\h163.7.\h163.8.\h163.9.
\h163.10.\h163.11.\h163.13.\h163.14.\h163.15.\h163.16.\h163.17.\h163.18.
\h163.19.\h163.20.\h163.21.\h163.22.\h163.23.\h163.24.\h163.25.\h163.26.
\h163.27.\h163.28.\h163.29.\h163.30.\h163.31.\h163.32.\h163.33.\h163.34.
\h163.35.\h163.36.\h163.37.\h163.38.\h163.39.\h163.40.\h163.41.\h163.43.
\h163.45.\h163.46.\h163.47.\h163.49.\h163.50.\h163.51.\h163.52.\h163.53.
\h163.54.\h163.55.\h163.57.\h163.58.\h163.59.\h163.61.\h163.63.\h163.65.
\h163.67.\h163.68.\h163.69.\h163.71.\h163.73.\h163.75.\h163.76.\h163.77.
\h163.79.\h163.83.\h163.85.\h163.91.\h163.95.\h163.97.\h163.99.\h163.107.
\h163.111.\h163.112.\h163.115.\h163.135.\h163.139.\h163.151.\h163.163.\h164.6.
\h164.8.\h164.10.\h164.11.\h164.12.\h164.13.\h164.14.\h164.15.\h164.16.
\h164.17.\h164.18.\h164.19.\h164.20.\h164.21.\h164.22.\h164.23.\h164.24.
\h164.25.\h164.26.\h164.27.\h164.28.\h164.29.\h164.30.\h164.31.\h164.32.
\h164.34.\h164.35.\h164.36.\h164.37.\h164.38.\h164.39.\h164.40.\h164.41.
\h164.42.\h164.44.\h164.45.\h164.46.\h164.47.\h164.48.\h164.49.\h164.50.
\h164.51.\h164.52.\h164.53.\h164.54.\h164.56.\h164.58.\h164.59.\h164.60.
\h164.62.\h164.64.\h164.65.\h164.66.\h164.68.\h164.70.\h164.72.\h164.74.
\h164.76.\h164.78.\h164.82.\h164.83.\h164.86.\h164.88.\h164.94.\h164.98.
\h164.106.\h164.110.\h164.114.\h164.122.\h164.134.\h164.138.\h164.150.\h164.162.
\h165.3.\h165.5.\h165.9.\h165.11.\h165.12.\h165.13.\h165.14.\h165.15.
\h165.16.\h165.17.\h165.18.\h165.19.\h165.20.\h165.21.\h165.22.\h165.23.
\h165.24.\h165.25.\h165.26.\h165.27.\h165.28.\h165.29.\h165.30.\h165.31.
\h165.32.\h165.33.\h165.34.\h165.35.\h165.36.\h165.37.\h165.38.\h165.39.
\h165.40.\h165.41.\h165.42.\h165.43.\h165.44.\h165.45.\h165.47.\h165.48.
\h165.49.\h165.51.\h165.53.\h165.54.\h165.55.\h165.57.\h165.59.\h165.61.
\h165.63.\h165.65.\h165.66.\h165.67.\h165.69.\h165.71.\h165.73.\h165.74.
\h165.75.\h165.77.\h165.81.\h165.85.\h165.87.\h165.89.\h165.93.\h165.97.
\h165.105.\h165.109.\h165.113.\h165.117.\h165.121.\h165.133.\h165.137.\h165.141.
\h165.159.\h165.165.\h166.4.\h166.8.\h166.10.\h166.11.\h166.12.\h166.13.
\h166.14.\h166.16.\h166.17.\h166.18.\h166.19.\h166.20.\h166.21.\h166.22.
\h166.24.\h166.25.\h166.26.\h166.28.\h166.29.\h166.30.\h166.31.\h166.32.
\h166.33.\h166.34.\h166.35.\h166.36.\h166.37.\h166.38.\h166.39.\h166.40.
\h166.41.\h166.42.\h166.43.\h166.44.\h166.45.\h166.46.\h166.47.\h166.48.
\h166.49.\h166.50.\h166.52.\h166.54.\h166.55.\h166.56.\h166.58.\h166.60.
\h166.64.\h166.66.\h166.68.\h166.70.\h166.72.\h166.73.\h166.76.\h166.80.
\h166.82.\h166.84.\h166.88.\h166.92.\h166.96.\h166.100.\h166.104.\h166.108.
\h166.112.\h166.116.\h166.120.\h166.130.\h166.136.\h166.148.\h166.154.\h166.164.
\h167.7.\h167.8.\h167.9.\h167.10.\h167.11.\h167.13.\h167.15.\h167.17.
\h167.18.\h167.19.\h167.20.\h167.21.\h167.22.\h167.23.\h167.24.\h167.25.
\h167.26.\h167.27.\h167.28.\h167.29.\h167.30.\h167.31.\h167.32.\h167.33.
\h167.34.\h167.35.\h167.36.\h167.37.\h167.38.\h167.39.\h167.40.\h167.41.
\h167.42.\h167.43.\h167.44.\h167.45.\h167.47.\h167.49.\h167.50.\h167.51.
\h167.53.\h167.55.\h167.59.\h167.61.\h167.63.\h167.64.\h167.65.\h167.67.
\h167.69.\h167.71.\h167.73.\h167.75.\h167.79.\h167.83.\h167.86.\h167.87.
\h167.89.\h167.91.\h167.95.\h167.101.\h167.103.\h167.107.\h167.109.\h167.115.
\h167.119.\h167.125.\h167.131.\h167.147.\h167.167.\h168.6.\h168.8.\h168.9.
\h168.10.\h168.11.\h168.12.\h168.13.\h168.14.\h168.15.\h168.16.\h168.18.
\h168.19.\h168.20.\h168.21.\h168.22.\h168.23.\h168.24.\h168.25.\h168.26.
\h168.27.\h168.28.\h168.30.\h168.31.\h168.32.\h168.33.\h168.34.\h168.35.
\h168.36.\h168.37.\h168.38.\h168.39.\h168.40.\h168.42.\h168.44.\h168.45.
\h168.46.\h168.47.\h168.48.\h168.50.\h168.52.\h168.54.\h168.58.\h168.60.
\h168.62.\h168.63.\h168.64.\h168.66.\h168.68.\h168.70.\h168.72.\h168.74.
\h168.75.\h168.78.\h168.84.\h168.86.\h168.90.\h168.92.\h168.96.\h168.102.
\h168.104.\h168.106.\h168.108.\h168.111.\h168.114.\h168.118.\h168.130.\h168.138.
\h168.150.\h168.162.\h169.7.\h169.9.\h169.10.\h169.11.\h169.12.\h169.13.
\h169.15.\h169.16.\h169.17.\h169.18.\h169.19.\h169.20.\h169.21.\h169.22.
\h169.23.\h169.24.\h169.25.\h169.26.\h169.27.\h169.28.\h169.29.\h169.30.
\h169.31.\h169.33.\h169.34.\h169.35.\h169.36.\h169.37.\h169.39.\h169.41.
\h169.43.\h169.44.\h169.45.\h169.47.\h169.49.\h169.51.\h169.52.\h169.53.
\h169.55.\h169.57.\h169.59.\h169.61.\h169.63.\h169.64.\h169.65.\h169.67.
\h169.69.\h169.73.\h169.75.\h169.77.\h169.79.\h169.81.\h169.82.\h169.83.
\h169.85.\h169.89.\h169.91.\h169.97.\h169.100.\h169.101.\h169.103.\h169.109.
\h169.113.\h169.115.\h169.117.\h169.121.\h169.133.\h169.145.\h169.149.\h169.157.
\h169.161.\h170.8.\h170.10.\h170.11.\h170.12.\h170.14.\h170.15.\h170.16.
\h170.17.\h170.18.\h170.19.\h170.20.\h170.21.\h170.22.\h170.23.\h170.24.
\h170.25.\h170.26.\h170.27.\h170.28.\h170.29.\h170.30.\h170.31.\h170.32.
\h170.34.\h170.35.\h170.36.\h170.37.\h170.38.\h170.40.\h170.41.\h170.42.
\h170.43.\h170.44.\h170.46.\h170.47.\h170.48.\h170.50.\h170.52.\h170.53.
\h170.54.\h170.56.\h170.58.\h170.60.\h170.61.\h170.62.\h170.64.\h170.66.
\h170.68.\h170.70.\h170.71.\h170.72.\h170.74.\h170.76.\h170.78.\h170.80.
\h170.84.\h170.86.\h170.88.\h170.89.\h170.92.\h170.96.\h170.98.\h170.100.
\h170.102.\h170.104.\h170.108.\h170.112.\h170.116.\h170.120.\h170.128.\h170.132.
\h170.140.\h170.144.\h170.158.\h170.160.\h171.6.\h171.7.\h171.9.\h171.10.
\h171.11.\h171.12.\h171.13.\h171.15.\h171.16.\h171.17.\h171.18.\h171.19.
\h171.20.\h171.21.\h171.23.\h171.24.\h171.25.\h171.26.\h171.27.\h171.29.
\h171.30.\h171.31.\h171.32.\h171.33.\h171.35.\h171.36.\h171.37.\h171.39.
\h171.41.\h171.42.\h171.43.\h171.44.\h171.45.\h171.47.\h171.48.\h171.49.
\h171.51.\h171.52.\h171.53.\h171.55.\h171.57.\h171.59.\h171.60.\h171.61.
\h171.63.\h171.65.\h171.67.\h171.69.\h171.71.\h171.72.\h171.73.\h171.75.
\h171.77.\h171.79.\h171.81.\h171.83.\h171.87.\h171.91.\h171.93.\h171.95.
\h171.99.\h171.103.\h171.105.\h171.111.\h171.114.\h171.115.\h171.119.\h171.127.
\h171.129.\h171.131.\h171.135.\h171.141.\h171.147.\h171.159.\h172.6.\h172.9.
\h172.10.\h172.12.\h172.13.\h172.14.\h172.15.\h172.16.\h172.17.\h172.18.
\h172.19.\h172.20.\h172.22.\h172.23.\h172.24.\h172.25.\h172.26.\h172.28.
\h172.29.\h172.30.\h172.31.\h172.32.\h172.33.\h172.34.\h172.35.\h172.36.
\h172.37.\h172.38.\h172.40.\h172.41.\h172.42.\h172.43.\h172.44.\h172.46.
\h172.47.\h172.48.\h172.50.\h172.51.\h172.52.\h172.54.\h172.56.\h172.57.
\h172.58.\h172.60.\h172.62.\h172.64.\h172.66.\h172.68.\h172.70.\h172.72.
\h172.74.\h172.76.\h172.78.\h172.79.\h172.80.\h172.82.\h172.84.\h172.85.
\h172.86.\h172.88.\h172.90.\h172.92.\h172.94.\h172.98.\h172.100.\h172.102.
\h172.106.\h172.110.\h172.112.\h172.114.\h172.118.\h172.124.\h172.130.\h172.157.
\h173.5.\h173.7.\h173.8.\h173.9.\h173.11.\h173.13.\h173.14.\h173.15.
\h173.17.\h173.18.\h173.19.\h173.20.\h173.21.\h173.23.\h173.24.\h173.25.
\h173.26.\h173.27.\h173.28.\h173.29.\h173.30.\h173.31.\h173.32.\h173.33.
\h173.35.\h173.36.\h173.37.\h173.38.\h173.39.\h173.40.\h173.41.\h173.43.
\h173.44.\h173.45.\h173.46.\h173.47.\h173.49.\h173.50.\h173.51.\h173.52.
\h173.53.\h173.55.\h173.56.\h173.57.\h173.58.\h173.59.\h173.61.\h173.63.
\h173.65.\h173.67.\h173.68.\h173.69.\h173.71.\h173.73.\h173.74.\h173.75.
\h173.77.\h173.81.\h173.83.\h173.85.\h173.86.\h173.87.\h173.89.\h173.93.
\h173.95.\h173.97.\h173.101.\h173.103.\h173.110.\h173.113.\h173.128.\h173.153.
\h173.155.\h173.158.\h174.6.\h174.8.\h174.9.\h174.10.\h174.12.\h174.14.
\h174.15.\h174.16.\h174.17.\h174.18.\h174.19.\h174.20.\h174.21.\h174.22.
\h174.23.\h174.24.\h174.26.\h174.27.\h174.28.\h174.29.\h174.30.\h174.31.
\h174.32.\h174.33.\h174.34.\h174.35.\h174.36.\h174.37.\h174.38.\h174.39.
\h174.40.\h174.41.\h174.42.\h174.44.\h174.45.\h174.46.\h174.47.\h174.48.
\h174.50.\h174.51.\h174.52.\h174.54.\h174.56.\h174.57.\h174.58.\h174.60.
\h174.62.\h174.63.\h174.64.\h174.65.\h174.66.\h174.68.\h174.70.\h174.72.
\h174.74.\h174.75.\h174.76.\h174.78.\h174.80.\h174.81.\h174.84.\h174.93.
\h174.96.\h174.98.\h174.99.\h174.102.\h174.105.\h174.108.\h174.112.\h174.124.
\h174.126.\h174.129.\h174.140.\h174.152.\h174.156.\h174.168.\h175.7.\h175.9.
\h175.10.\h175.11.\h175.12.\h175.13.\h175.15.\h175.16.\h175.17.\h175.19.
\h175.20.\h175.21.\h175.22.\h175.23.\h175.25.\h175.27.\h175.28.\h175.29.
\h175.31.\h175.32.\h175.33.\h175.34.\h175.35.\h175.36.\h175.37.\h175.38.
\h175.39.\h175.40.\h175.41.\h175.43.\h175.45.\h175.46.\h175.47.\h175.49.
\h175.51.\h175.52.\h175.53.\h175.55.\h175.57.\h175.58.\h175.59.\h175.61.
\h175.62.\h175.63.\h175.64.\h175.65.\h175.67.\h175.69.\h175.70.\h175.71.
\h175.73.\h175.74.\h175.75.\h175.76.\h175.79.\h175.82.\h175.83.\h175.88.
\h175.91.\h175.94.\h175.95.\h175.97.\h175.100.\h175.101.\h175.103.\h175.107.
\h175.109.\h175.111.\h175.123.\h175.127.\h175.139.\h175.151.\h175.155.\h176.6.
\h176.8.\h176.10.\h176.11.\h176.12.\h176.14.\h176.16.\h176.17.\h176.18.
\h176.19.\h176.20.\h176.21.\h176.22.\h176.23.\h176.24.\h176.25.\h176.26.
\h176.28.\h176.29.\h176.30.\h176.32.\h176.33.\h176.34.\h176.35.\h176.36.
\h176.38.\h176.40.\h176.41.\h176.42.\h176.44.\h176.46.\h176.47.\h176.48.
\h176.50.\h176.51.\h176.52.\h176.53.\h176.54.\h176.56.\h176.57.\h176.58.
\h176.59.\h176.62.\h176.63.\h176.64.\h176.65.\h176.66.\h176.68.\h176.69.
\h176.70.\h176.71.\h176.72.\h176.74.\h176.76.\h176.78.\h176.80.\h176.82.
\h176.83.\h176.86.\h176.90.\h176.92.\h176.94.\h176.96.\h176.98.\h176.106.
\h176.108.\h176.110.\h176.122.\h176.126.\h176.146.\h176.152.\h177.5.\h177.9.
\h177.11.\h177.12.\h177.13.\h177.15.\h177.16.\h177.17.\h177.18.\h177.19.
\h177.20.\h177.21.\h177.22.\h177.23.\h177.24.\h177.25.\h177.27.\h177.28.
\h177.29.\h177.30.\h177.31.\h177.32.\h177.33.\h177.34.\h177.35.\h177.36.
\h177.37.\h177.38.\h177.39.\h177.40.\h177.41.\h177.42.\h177.43.\h177.44.
\h177.45.\h177.46.\h177.47.\h177.48.\h177.49.\h177.51.\h177.52.\h177.53.
\h177.54.\h177.55.\h177.57.\h177.61.\h177.62.\h177.63.\h177.65.\h177.66.
\h177.67.\h177.69.\h177.71.\h177.72.\h177.73.\h177.75.\h177.77.\h177.79.
\h177.81.\h177.85.\h177.89.\h177.91.\h177.93.\h177.97.\h177.101.\h177.105.
\h177.107.\h177.117.\h177.123.\h177.137.\h177.141.\h177.147.\h178.4.\h178.7.
\h178.8.\h178.10.\h178.12.\h178.13.\h178.16.\h178.17.\h178.18.\h178.19.
\h178.20.\h178.21.\h178.22.\h178.23.\h178.24.\h178.26.\h178.27.\h178.28.
\h178.29.\h178.30.\h178.31.\h178.32.\h178.33.\h178.34.\h178.35.\h178.36.
\h178.37.\h178.38.\h178.40.\h178.42.\h178.43.\h178.44.\h178.45.\h178.46.
\h178.48.\h178.49.\h178.50.\h178.51.\h178.52.\h178.53.\h178.56.\h178.58.
\h178.60.\h178.61.\h178.62.\h178.64.\h178.68.\h178.70.\h178.72.\h178.73.
\h178.74.\h178.76.\h178.80.\h178.82.\h178.84.\h178.88.\h178.90.\h178.92.
\h178.94.\h178.96.\h178.100.\h178.104.\h178.108.\h178.112.\h178.118.\h178.136.
\h178.154.\h179.5.\h179.7.\h179.8.\h179.9.\h179.11.\h179.13.\h179.14.
\h179.15.\h179.17.\h179.18.\h179.19.\h179.21.\h179.22.\h179.23.\h179.25.
\h179.26.\h179.27.\h179.28.\h179.29.\h179.31.\h179.32.\h179.33.\h179.34.
\h179.35.\h179.36.\h179.37.\h179.38.\h179.39.\h179.41.\h179.43.\h179.44.
\h179.45.\h179.47.\h179.49.\h179.50.\h179.51.\h179.53.\h179.55.\h179.57.
\h179.59.\h179.61.\h179.63.\h179.64.\h179.65.\h179.67.\h179.69.\h179.71.
\h179.73.\h179.75.\h179.77.\h179.79.\h179.80.\h179.83.\h179.87.\h179.89.
\h179.95.\h179.99.\h179.107.\h179.125.\h179.131.\h179.151.\h179.167.\h179.179.
\h180.6.\h180.8.\h180.10.\h180.12.\h180.14.\h180.15.\h180.16.\h180.17.
\h180.18.\h180.20.\h180.21.\h180.22.\h180.24.\h180.26.\h180.27.\h180.28.
\h180.29.\h180.30.\h180.31.\h180.32.\h180.33.\h180.34.\h180.35.\h180.36.
\h180.38.\h180.39.\h180.40.\h180.42.\h180.43.\h180.44.\h180.45.\h180.46.
\h180.48.\h180.49.\h180.50.\h180.52.\h180.54.\h180.56.\h180.57.\h180.58.
\h180.59.\h180.60.\h180.62.\h180.63.\h180.66.\h180.68.\h180.69.\h180.70.
\h180.72.\h180.74.\h180.78.\h180.82.\h180.84.\h180.90.\h180.96.\h180.102.
\h180.106.\h180.120.\h180.122.\h180.138.\h180.150.\h180.162.\h181.7.\h181.9.
\h181.11.\h181.13.\h181.14.\h181.15.\h181.16.\h181.17.\h181.18.\h181.19.
\h181.20.\h181.21.\h181.22.\h181.23.\h181.24.\h181.25.\h181.26.\h181.27.
\h181.28.\h181.29.\h181.30.\h181.31.\h181.32.\h181.33.\h181.34.\h181.35.
\h181.37.\h181.38.\h181.39.\h181.41.\h181.43.\h181.44.\h181.45.\h181.46.
\h181.47.\h181.49.\h181.51.\h181.53.\h181.55.\h181.57.\h181.58.\h181.61.
\h181.65.\h181.67.\h181.69.\h181.71.\h181.73.\h181.75.\h181.76.\h181.77.
\h181.79.\h181.85.\h181.89.\h181.91.\h181.93.\h181.97.\h181.105.\h181.109.
\h181.121.\h181.133.\h181.145.\h181.149.\h182.6.\h182.8.\h182.11.\h182.12.
\h182.14.\h182.16.\h182.17.\h182.18.\h182.19.\h182.20.\h182.21.\h182.22.
\h182.23.\h182.24.\h182.25.\h182.26.\h182.27.\h182.28.\h182.29.\h182.30.
\h182.31.\h182.32.\h182.33.\h182.34.\h182.36.\h182.38.\h182.39.\h182.40.
\h182.41.\h182.42.\h182.43.\h182.44.\h182.46.\h182.47.\h182.48.\h182.50.
\h182.52.\h182.53.\h182.54.\h182.56.\h182.58.\h182.60.\h182.62.\h182.64.
\h182.65.\h182.66.\h182.68.\h182.70.\h182.72.\h182.74.\h182.76.\h182.80.
\h182.84.\h182.88.\h182.92.\h182.104.\h182.116.\h182.120.\h182.132.\h182.140.
\h182.144.\h182.146.\h183.7.\h183.9.\h183.11.\h183.12.\h183.13.\h183.15.
\h183.17.\h183.18.\h183.19.\h183.20.\h183.21.\h183.23.\h183.24.\h183.25.
\h183.26.\h183.27.\h183.28.\h183.29.\h183.30.\h183.31.\h183.32.\h183.33.
\h183.35.\h183.36.\h183.37.\h183.38.\h183.39.\h183.41.\h183.42.\h183.43.
\h183.45.\h183.47.\h183.48.\h183.49.\h183.50.\h183.51.\h183.53.\h183.54.
\h183.55.\h183.57.\h183.59.\h183.60.\h183.61.\h183.63.\h183.65.\h183.67.
\h183.69.\h183.71.\h183.72.\h183.75.\h183.81.\h183.87.\h183.91.\h183.103.
\h183.111.\h183.115.\h183.117.\h183.119.\h183.129.\h183.143.\h184.6.\h184.8.
\h184.10.\h184.12.\h184.13.\h184.14.\h184.16.\h184.17.\h184.18.\h184.19.
\h184.20.\h184.22.\h184.23.\h184.24.\h184.25.\h184.26.\h184.27.\h184.28.
\h184.29.\h184.30.\h184.31.\h184.32.\h184.34.\h184.35.\h184.36.\h184.37.
\h184.38.\h184.39.\h184.40.\h184.41.\h184.42.\h184.43.\h184.44.\h184.46.
\h184.48.\h184.49.\h184.50.\h184.52.\h184.54.\h184.56.\h184.58.\h184.60.
\h184.62.\h184.64.\h184.67.\h184.68.\h184.70.\h184.72.\h184.74.\h184.80.
\h184.82.\h184.86.\h184.88.\h184.90.\h184.100.\h184.102.\h184.114.\h184.118.
\h184.136.\h184.148.\h185.5.\h185.7.\h185.8.\h185.9.\h185.11.\h185.13.
\h185.14.\h185.15.\h185.16.\h185.17.\h185.19.\h185.20.\h185.21.\h185.22.
\h185.23.\h185.25.\h185.26.\h185.27.\h185.28.\h185.29.\h185.30.\h185.31.
\h185.32.\h185.33.\h185.35.\h185.36.\h185.37.\h185.38.\h185.39.\h185.40.
\h185.41.\h185.43.\h185.45.\h185.47.\h185.48.\h185.49.\h185.51.\h185.53.
\h185.55.\h185.56.\h185.57.\h185.59.\h185.61.\h185.62.\h185.63.\h185.65.
\h185.67.\h185.69.\h185.71.\h185.73.\h185.79.\h185.83.\h185.85.\h185.89.
\h185.93.\h185.95.\h185.99.\h185.101.\h185.103.\h185.107.\h185.110.\h185.113.
\h185.119.\h185.125.\h185.143.\h185.145.\h186.6.\h186.9.\h186.12.\h186.14.
\h186.15.\h186.16.\h186.17.\h186.18.\h186.20.\h186.21.\h186.22.\h186.23.
\h186.24.\h186.25.\h186.26.\h186.27.\h186.28.\h186.29.\h186.30.\h186.31.
\h186.32.\h186.33.\h186.34.\h186.36.\h186.38.\h186.39.\h186.40.\h186.41.
\h186.42.\h186.43.\h186.44.\h186.46.\h186.48.\h186.50.\h186.51.\h186.52.
\h186.54.\h186.56.\h186.57.\h186.58.\h186.60.\h186.62.\h186.64.\h186.66.
\h186.68.\h186.70.\h186.72.\h186.74.\h186.76.\h186.78.\h186.81.\h186.82.
\h186.84.\h186.86.\h186.90.\h186.92.\h186.96.\h186.99.\h186.102.\h186.108.
\h186.114.\h186.116.\h186.144.\h187.7.\h187.9.\h187.11.\h187.13.\h187.15.
\h187.16.\h187.19.\h187.20.\h187.21.\h187.22.\h187.23.\h187.24.\h187.25.
\h187.26.\h187.27.\h187.28.\h187.29.\h187.31.\h187.32.\h187.33.\h187.34.
\h187.35.\h187.36.\h187.37.\h187.39.\h187.40.\h187.41.\h187.42.\h187.43.
\h187.45.\h187.46.\h187.47.\h187.49.\h187.51.\h187.52.\h187.53.\h187.55.
\h187.57.\h187.59.\h187.61.\h187.63.\h187.65.\h187.67.\h187.69.\h187.70.
\h187.73.\h187.75.\h187.79.\h187.85.\h187.87.\h187.88.\h187.95.\h187.97.
\h187.101.\h187.115.\h187.142.\h188.6.\h188.8.\h188.10.\h188.11.\h188.14.
\h188.15.\h188.16.\h188.17.\h188.18.\h188.19.\h188.20.\h188.21.\h188.22.
\h188.24.\h188.25.\h188.26.\h188.28.\h188.30.\h188.31.\h188.32.\h188.34.
\h188.35.\h188.36.\h188.38.\h188.40.\h188.41.\h188.42.\h188.43.\h188.44.
\h188.45.\h188.46.\h188.48.\h188.50.\h188.52.\h188.54.\h188.56.\h188.58.
\h188.59.\h188.60.\h188.62.\h188.66.\h188.68.\h188.70.\h188.71.\h188.72.
\h188.80.\h188.84.\h188.86.\h188.94.\h188.98.\h188.101.\h188.113.\h188.140.
\h188.143.\h189.5.\h189.9.\h189.13.\h189.14.\h189.15.\h189.17.\h189.18.
\h189.19.\h189.20.\h189.21.\h189.23.\h189.24.\h189.25.\h189.27.\h189.29.
\h189.30.\h189.31.\h189.32.\h189.33.\h189.35.\h189.37.\h189.39.\h189.41.
\h189.42.\h189.43.\h189.44.\h189.45.\h189.47.\h189.49.\h189.51.\h189.53.
\h189.54.\h189.55.\h189.57.\h189.59.\h189.65.\h189.67.\h189.69.\h189.73.
\h189.78.\h189.79.\h189.84.\h189.85.\h189.93.\h189.99.\h189.111.\h189.114.
\h189.135.\h189.141.\h190.4.\h190.8.\h190.10.\h190.12.\h190.13.\h190.14.
\h190.15.\h190.16.\h190.17.\h190.18.\h190.19.\h190.20.\h190.22.\h190.24.
\h190.25.\h190.26.\h190.28.\h190.29.\h190.30.\h190.31.\h190.32.\h190.34.
\h190.36.\h190.37.\h190.38.\h190.39.\h190.40.\h190.42.\h190.43.\h190.44.
\h190.46.\h190.48.\h190.50.\h190.52.\h190.54.\h190.55.\h190.56.\h190.58.
\h190.64.\h190.66.\h190.67.\h190.68.\h190.70.\h190.72.\h190.76.\h190.82.
\h190.84.\h190.85.\h190.88.\h190.94.\h190.98.\h190.100.\h190.106.\h190.112.
\h190.130.\h190.140.\h191.9.\h191.11.\h191.13.\h191.14.\h191.15.\h191.17.
\h191.18.\h191.19.\h191.21.\h191.23.\h191.25.\h191.26.\h191.27.\h191.28.
\h191.29.\h191.30.\h191.31.\h191.32.\h191.33.\h191.35.\h191.36.\h191.37.
\h191.38.\h191.39.\h191.40.\h191.41.\h191.43.\h191.45.\h191.47.\h191.48.
\h191.50.\h191.53.\h191.54.\h191.55.\h191.56.\h191.57.\h191.59.\h191.64.
\h191.65.\h191.67.\h191.68.\h191.69.\h191.71.\h191.77.\h191.79.\h191.81.
\h191.83.\h191.95.\h191.99.\h191.101.\h191.111.\h191.119.\h191.139.\h191.167.
\h192.6.\h192.7.\h192.8.\h192.10.\h192.12.\h192.13.\h192.14.\h192.16.
\h192.17.\h192.18.\h192.20.\h192.21.\h192.22.\h192.24.\h192.26.\h192.27.
\h192.28.\h192.29.\h192.30.\h192.32.\h192.33.\h192.34.\h192.35.\h192.36.
\h192.37.\h192.38.\h192.39.\h192.40.\h192.42.\h192.44.\h192.45.\h192.46.
\h192.48.\h192.50.\h192.51.\h192.52.\h192.54.\h192.56.\h192.60.\h192.64.
\h192.66.\h192.68.\h192.72.\h192.75.\h192.78.\h192.82.\h192.90.\h192.98.
\h192.110.\h192.114.\h192.134.\h192.138.\h193.7.\h193.9.\h193.11.\h193.13.
\h193.15.\h193.16.\h193.17.\h193.19.\h193.20.\h193.21.\h193.22.\h193.23.
\h193.25.\h193.27.\h193.28.\h193.29.\h193.31.\h193.33.\h193.34.\h193.35.
\h193.36.\h193.37.\h193.38.\h193.39.\h193.40.\h193.43.\h193.45.\h193.46.
\h193.47.\h193.49.\h193.51.\h193.53.\h193.55.\h193.57.\h193.59.\h193.61.
\h193.64.\h193.65.\h193.67.\h193.69.\h193.71.\h193.73.\h193.75.\h193.77.
\h193.81.\h193.85.\h193.93.\h193.97.\h193.105.\h193.109.\h193.121.\h193.133.
\h194.4.\h194.8.\h194.10.\h194.11.\h194.12.\h194.14.\h194.16.\h194.17.
\h194.18.\h194.19.\h194.20.\h194.21.\h194.22.\h194.23.\h194.24.\h194.26.
\h194.28.\h194.29.\h194.30.\h194.31.\h194.32.\h194.33.\h194.34.\h194.36.
\h194.37.\h194.38.\h194.40.\h194.41.\h194.42.\h194.44.\h194.46.\h194.47.
\h194.48.\h194.52.\h194.53.\h194.56.\h194.58.\h194.62.\h194.64.\h194.66.
\h194.68.\h194.70.\h194.74.\h194.76.\h194.80.\h194.84.\h194.92.\h194.96.
\h194.104.\h194.132.\h194.164.\h195.3.\h195.11.\h195.13.\h195.15.\h195.17.
\h195.18.\h195.19.\h195.20.\h195.21.\h195.23.\h195.25.\h195.27.\h195.28.
\h195.29.\h195.30.\h195.31.\h195.33.\h195.35.\h195.36.\h195.37.\h195.39.
\h195.40.\h195.41.\h195.42.\h195.43.\h195.45.\h195.46.\h195.47.\h195.51.
\h195.53.\h195.55.\h195.59.\h195.63.\h195.65.\h195.67.\h195.69.\h195.72.
\h195.75.\h195.79.\h195.91.\h195.95.\h195.103.\h195.129.\h195.147.\h195.195.
\h196.6.\h196.7.\h196.10.\h196.12.\h196.14.\h196.16.\h196.18.\h196.20.
\h196.22.\h196.25.\h196.26.\h196.27.\h196.28.\h196.29.\h196.30.\h196.31.
\h196.32.\h196.33.\h196.34.\h196.36.\h196.37.\h196.38.\h196.39.\h196.40.
\h196.42.\h196.43.\h196.44.\h196.46.\h196.50.\h196.52.\h196.54.\h196.58.
\h196.62.\h196.64.\h196.66.\h196.68.\h196.70.\h196.74.\h196.78.\h196.82.
\h196.94.\h196.100.\h196.118.\h196.124.\h196.130.\h196.166.\h196.178.\h197.5.
\h197.7.\h197.9.\h197.13.\h197.14.\h197.17.\h197.19.\h197.21.\h197.23.
\h197.25.\h197.26.\h197.27.\h197.29.\h197.30.\h197.31.\h197.32.\h197.33.
\h197.35.\h197.37.\h197.39.\h197.41.\h197.43.\h197.45.\h197.46.\h197.49.
\h197.51.\h197.53.\h197.57.\h197.59.\h197.61.\h197.62.\h197.63.\h197.65.
\h197.67.\h197.71.\h197.73.\h197.74.\h197.77.\h197.89.\h197.91.\h197.95.
\h197.101.\h197.113.\h197.117.\h197.137.\h197.149.\h197.161.\h198.6.\h198.8.
\h198.12.\h198.14.\h198.15.\h198.16.\h198.17.\h198.18.\h198.20.\h198.21.
\h198.22.\h198.23.\h198.24.\h198.25.\h198.26.\h198.27.\h198.28.\h198.30.
\h198.32.\h198.33.\h198.34.\h198.36.\h198.38.\h198.40.\h198.42.\h198.44.
\h198.48.\h198.52.\h198.54.\h198.56.\h198.58.\h198.60.\h198.62.\h198.63.
\h198.64.\h198.66.\h198.72.\h198.74.\h198.84.\h198.88.\h198.96.\h198.100.
\h198.108.\h198.120.\h198.132.\h198.148.\h199.7.\h199.11.\h199.13.\h199.15.
\h199.16.\h199.17.\h199.19.\h199.21.\h199.23.\h199.24.\h199.25.\h199.27.
\h199.29.\h199.31.\h199.33.\h199.35.\h199.37.\h199.39.\h199.41.\h199.43.
\h199.45.\h199.47.\h199.49.\h199.51.\h199.52.\h199.53.\h199.55.\h199.57.
\h199.59.\h199.61.\h199.63.\h199.65.\h199.67.\h199.71.\h199.79.\h199.83.
\h199.91.\h199.95.\h199.99.\h199.103.\h199.119.\h199.131.\h199.139.\h200.12.
\h200.14.\h200.15.\h200.18.\h200.20.\h200.22.\h200.23.\h200.24.\h200.25.
\h200.26.\h200.28.\h200.29.\h200.30.\h200.32.\h200.34.\h200.36.\h200.38.
\h200.40.\h200.41.\h200.42.\h200.43.\h200.44.\h200.46.\h200.48.\h200.50.
\h200.52.\h200.54.\h200.56.\h200.58.\h200.59.\h200.60.\h200.62.\h200.64.
\h200.66.\h200.68.\h200.70.\h200.74.\h200.78.\h200.82.\h200.86.\h200.90.
\h200.92.\h200.94.\h200.102.\h200.110.\h200.118.\h200.128.\h200.130.\h201.5.
\h201.9.\h201.11.\h201.13.\h201.14.\h201.15.\h201.19.\h201.21.\h201.23.
\h201.24.\h201.25.\h201.27.\h201.28.\h201.29.\h201.30.\h201.31.\h201.33.
\h201.35.\h201.36.\h201.37.\h201.39.\h201.41.\h201.42.\h201.43.\h201.45.
\h201.49.\h201.51.\h201.53.\h201.55.\h201.57.\h201.59.\h201.61.\h201.63.
\h201.65.\h201.69.\h201.73.\h201.75.\h201.77.\h201.81.\h201.87.\h201.89.
\h201.93.\h201.99.\h201.101.\h201.117.\h201.129.\h202.4.\h202.8.\h202.12.
\h202.13.\h202.14.\h202.16.\h202.20.\h202.22.\h202.24.\h202.25.\h202.26.
\h202.27.\h202.28.\h202.29.\h202.30.\h202.31.\h202.32.\h202.34.\h202.35.
\h202.36.\h202.38.\h202.40.\h202.42.\h202.44.\h202.45.\h202.46.\h202.48.
\h202.50.\h202.52.\h202.55.\h202.56.\h202.58.\h202.60.\h202.62.\h202.64.
\h202.70.\h202.72.\h202.73.\h202.76.\h202.80.\h202.88.\h202.92.\h202.100.
\h202.102.\h202.109.\h202.112.\h202.127.\h203.11.\h203.13.\h203.15.\h203.17.
\h203.19.\h203.20.\h203.21.\h203.23.\h203.25.\h203.26.\h203.27.\h203.28.
\h203.29.\h203.31.\h203.33.\h203.35.\h203.37.\h203.38.\h203.39.\h203.40.
\h203.41.\h203.43.\h203.44.\h203.45.\h203.47.\h203.50.\h203.51.\h203.53.
\h203.55.\h203.56.\h203.57.\h203.59.\h203.62.\h203.63.\h203.69.\h203.71.
\h203.73.\h203.74.\h203.75.\h203.80.\h203.83.\h203.92.\h203.98.\h203.101.
\h203.107.\h203.123.\h203.125.\h203.128.\h204.6.\h204.10.\h204.12.\h204.14.
\h204.18.\h204.19.\h204.20.\h204.21.\h204.22.\h204.24.\h204.26.\h204.27.
\h204.28.\h204.30.\h204.32.\h204.34.\h204.35.\h204.36.\h204.38.\h204.39.
\h204.40.\h204.42.\h204.44.\h204.45.\h204.46.\h204.48.\h204.50.\h204.51.
\h204.52.\h204.54.\h204.56.\h204.57.\h204.58.\h204.60.\h204.62.\h204.63.
\h204.68.\h204.69.\h204.72.\h204.75.\h204.78.\h204.87.\h204.90.\h204.94.
\h204.96.\h204.99.\h204.100.\h204.122.\h204.126.\h205.5.\h205.7.\h205.9.
\h205.10.\h205.11.\h205.13.\h205.16.\h205.17.\h205.18.\h205.19.\h205.21.
\h205.23.\h205.25.\h205.27.\h205.28.\h205.29.\h205.31.\h205.33.\h205.34.
\h205.35.\h205.37.\h205.39.\h205.40.\h205.41.\h205.43.\h205.45.\h205.46.
\h205.47.\h205.49.\h205.51.\h205.52.\h205.53.\h205.55.\h205.58.\h205.61.
\h205.62.\h205.63.\h205.65.\h205.67.\h205.70.\h205.71.\h205.73.\h205.79.
\h205.85.\h205.89.\h205.91.\h205.93.\h205.97.\h205.121.\h205.125.\h206.8.
\h206.10.\h206.12.\h206.14.\h206.16.\h206.17.\h206.18.\h206.19.\h206.20.
\h206.22.\h206.23.\h206.24.\h206.26.\h206.27.\h206.28.\h206.29.\h206.30.
\h206.32.\h206.34.\h206.35.\h206.36.\h206.38.\h206.41.\h206.42.\h206.44.
\h206.45.\h206.48.\h206.50.\h206.52.\h206.53.\h206.56.\h206.60.\h206.62.
\h206.64.\h206.68.\h206.71.\h206.72.\h206.74.\h206.78.\h206.84.\h206.88.
\h206.90.\h206.92.\h206.96.\h206.116.\h206.122.\h206.140.\h207.11.\h207.12.
\h207.15.\h207.19.\h207.21.\h207.22.\h207.23.\h207.24.\h207.25.\h207.27.
\h207.28.\h207.29.\h207.30.\h207.31.\h207.33.\h207.35.\h207.36.\h207.39.
\h207.42.\h207.43.\h207.45.\h207.47.\h207.49.\h207.51.\h207.55.\h207.57.
\h207.59.\h207.60.\h207.61.\h207.63.\h207.67.\h207.71.\h207.73.\h207.83.
\h207.87.\h207.89.\h207.93.\h207.95.\h207.107.\h207.111.\h207.117.\h208.4.
\h208.10.\h208.13.\h208.16.\h208.18.\h208.20.\h208.22.\h208.24.\h208.26.
\h208.28.\h208.30.\h208.32.\h208.34.\h208.38.\h208.40.\h208.42.\h208.44.
\h208.46.\h208.49.\h208.50.\h208.54.\h208.56.\h208.58.\h208.60.\h208.63.
\h208.64.\h208.66.\h208.70.\h208.72.\h208.78.\h208.82.\h208.86.\h208.88.
\h208.94.\h208.106.\h208.124.\h209.5.\h209.9.\h209.11.\h209.13.\h209.14.
\h209.17.\h209.19.\h209.21.\h209.23.\h209.25.\h209.26.\h209.27.\h209.29.
\h209.30.\h209.32.\h209.33.\h209.34.\h209.35.\h209.37.\h209.38.\h209.39.
\h209.41.\h209.43.\h209.44.\h209.45.\h209.47.\h209.49.\h209.51.\h209.53.
\h209.55.\h209.57.\h209.59.\h209.61.\h209.65.\h209.67.\h209.68.\h209.69.
\h209.77.\h209.81.\h209.83.\h209.89.\h209.95.\h209.101.\h209.119.\h209.121.
\h210.6.\h210.8.\h210.10.\h210.12.\h210.15.\h210.16.\h210.18.\h210.20.
\h210.21.\h210.24.\h210.26.\h210.27.\h210.28.\h210.30.\h210.32.\h210.33.
\h210.34.\h210.36.\h210.38.\h210.40.\h210.42.\h210.43.\h210.44.\h210.48.
\h210.50.\h210.52.\h210.54.\h210.56.\h210.57.\h210.60.\h210.64.\h210.66.
\h210.70.\h210.72.\h210.74.\h210.76.\h210.78.\h210.88.\h210.90.\h210.92.
\h210.96.\h210.120.\h211.7.\h211.10.\h211.11.\h211.13.\h211.15.\h211.17.
\h211.19.\h211.20.\h211.22.\h211.23.\h211.25.\h211.26.\h211.27.\h211.28.
\h211.29.\h211.31.\h211.33.\h211.35.\h211.37.\h211.38.\h211.39.\h211.41.
\h211.43.\h211.45.\h211.46.\h211.47.\h211.49.\h211.51.\h211.53.\h211.55.
\h211.59.\h211.61.\h211.63.\h211.64.\h211.67.\h211.69.\h211.71.\h211.73.
\h211.79.\h211.87.\h211.91.\h211.115.\h211.119.\h212.12.\h212.14.\h212.17.
\h212.18.\h212.20.\h212.22.\h212.26.\h212.27.\h212.28.\h212.29.\h212.30.
\h212.32.\h212.34.\h212.35.\h212.36.\h212.37.\h212.38.\h212.40.\h212.41.
\h212.42.\h212.44.\h212.46.\h212.48.\h212.50.\h212.52.\h212.53.\h212.54.
\h212.56.\h212.60.\h212.62.\h212.66.\h212.68.\h212.70.\h212.71.\h212.86.
\h212.90.\h212.114.\h212.116.\h213.5.\h213.9.\h213.11.\h213.15.\h213.17.
\h213.18.\h213.21.\h213.25.\h213.26.\h213.27.\h213.29.\h213.30.\h213.31.
\h213.33.\h213.35.\h213.37.\h213.39.\h213.40.\h213.41.\h213.42.\h213.43.
\h213.45.\h213.47.\h213.48.\h213.49.\h213.51.\h213.53.\h213.57.\h213.59.
\h213.61.\h213.67.\h213.69.\h213.85.\h213.87.\h213.113.\h214.4.\h214.10.
\h214.12.\h214.16.\h214.24.\h214.25.\h214.26.\h214.27.\h214.28.\h214.29.
\h214.30.\h214.31.\h214.32.\h214.34.\h214.36.\h214.38.\h214.39.\h214.40.
\h214.42.\h214.44.\h214.46.\h214.48.\h214.52.\h214.54.\h214.56.\h214.58.
\h214.61.\h214.64.\h214.66.\h214.68.\h214.82.\h214.84.\h214.94.\h214.106.
\h214.118.\h215.5.\h215.11.\h215.14.\h215.15.\h215.19.\h215.20.\h215.22.
\h215.23.\h215.24.\h215.25.\h215.26.\h215.27.\h215.29.\h215.31.\h215.32.
\h215.33.\h215.35.\h215.37.\h215.38.\h215.39.\h215.41.\h215.42.\h215.43.
\h215.45.\h215.47.\h215.49.\h215.51.\h215.55.\h215.57.\h215.59.\h215.63.
\h215.65.\h215.67.\h215.75.\h215.77.\h215.81.\h215.85.\h215.89.\h215.95.
\h215.113.\h215.115.\h216.8.\h216.10.\h216.12.\h216.14.\h216.18.\h216.21.
\h216.22.\h216.23.\h216.24.\h216.25.\h216.26.\h216.27.\h216.28.\h216.30.
\h216.32.\h216.33.\h216.34.\h216.36.\h216.38.\h216.40.\h216.42.\h216.44.
\h216.46.\h216.48.\h216.50.\h216.54.\h216.56.\h216.58.\h216.60.\h216.62.
\h216.64.\h216.66.\h216.68.\h216.74.\h216.84.\h216.86.\h216.114.\h217.7.
\h217.13.\h217.17.\h217.19.\h217.21.\h217.22.\h217.25.\h217.27.\h217.29.
\h217.31.\h217.33.\h217.34.\h217.35.\h217.37.\h217.39.\h217.40.\h217.41.
\h217.42.\h217.43.\h217.45.\h217.47.\h217.51.\h217.53.\h217.55.\h217.57.
\h217.58.\h217.59.\h217.61.\h217.67.\h217.83.\h217.85.\h217.112.\h217.139.
\h218.8.\h218.12.\h218.13.\h218.14.\h218.18.\h218.20.\h218.22.\h218.24.
\h218.26.\h218.28.\h218.29.\h218.30.\h218.32.\h218.34.\h218.36.\h218.38.
\h218.40.\h218.41.\h218.44.\h218.46.\h218.48.\h218.50.\h218.56.\h218.58.
\h218.66.\h218.83.\h218.110.\h218.113.\h219.7.\h219.11.\h219.12.\h219.17.
\h219.19.\h219.21.\h219.23.\h219.25.\h219.27.\h219.29.\h219.30.\h219.31.
\h219.33.\h219.35.\h219.36.\h219.37.\h219.39.\h219.41.\h219.43.\h219.45.
\h219.47.\h219.49.\h219.54.\h219.57.\h219.66.\h219.81.\h219.84.\h219.93.
\h219.105.\h219.111.\h220.10.\h220.12.\h220.16.\h220.18.\h220.20.\h220.22.
\h220.25.\h220.26.\h220.27.\h220.28.\h220.30.\h220.32.\h220.34.\h220.35.
\h220.36.\h220.37.\h220.38.\h220.40.\h220.49.\h220.52.\h220.54.\h220.55.
\h220.56.\h220.64.\h220.76.\h220.80.\h220.82.\h220.100.\h220.110.\h220.118.
\h220.138.\h220.166.\h221.9.\h221.11.\h221.13.\h221.17.\h221.19.\h221.20.
\h221.21.\h221.23.\h221.25.\h221.26.\h221.27.\h221.29.\h221.31.\h221.33.
\h221.35.\h221.36.\h221.38.\h221.39.\h221.41.\h221.44.\h221.47.\h221.53.
\h221.55.\h221.63.\h221.65.\h221.71.\h221.77.\h221.81.\h221.89.\h221.109.
\h221.137.\h222.6.\h222.10.\h222.12.\h222.16.\h222.18.\h222.19.\h222.20.
\h222.21.\h222.22.\h222.24.\h222.26.\h222.27.\h222.28.\h222.30.\h222.32.
\h222.33.\h222.34.\h222.36.\h222.38.\h222.46.\h222.48.\h222.52.\h222.54.
\h222.56.\h222.60.\h222.63.\h222.64.\h222.72.\h222.80.\h222.92.\h222.104.
\h222.108.\h222.120.\h223.7.\h223.9.\h223.16.\h223.17.\h223.19.\h223.21.
\h223.22.\h223.23.\h223.25.\h223.27.\h223.28.\h223.31.\h223.33.\h223.35.
\h223.37.\h223.40.\h223.43.\h223.45.\h223.47.\h223.51.\h223.52.\h223.55.
\h223.61.\h223.63.\h223.75.\h223.79.\h223.91.\h223.103.\h223.223.\h224.14.
\h224.16.\h224.18.\h224.20.\h224.22.\h224.24.\h224.26.\h224.28.\h224.30.
\h224.32.\h224.34.\h224.35.\h224.38.\h224.39.\h224.40.\h224.41.\h224.44.
\h224.46.\h224.50.\h224.52.\h224.54.\h224.62.\h224.64.\h224.74.\h224.78.
\h224.90.\h224.102.\h224.194.\h225.9.\h225.15.\h225.17.\h225.19.\h225.21.
\h225.25.\h225.26.\h225.27.\h225.29.\h225.30.\h225.31.\h225.33.\h225.35.
\h225.37.\h225.38.\h225.39.\h225.43.\h225.45.\h225.51.\h225.53.\h225.57.
\h225.61.\h225.63.\h225.73.\h225.77.\h225.81.\h225.99.\h225.117.\h225.165.
\h225.177.\h226.4.\h226.6.\h226.8.\h226.13.\h226.16.\h226.18.\h226.20.
\h226.22.\h226.24.\h226.25.\h226.26.\h226.28.\h226.29.\h226.31.\h226.32.
\h226.34.\h226.36.\h226.40.\h226.42.\h226.44.\h226.50.\h226.52.\h226.60.
\h226.70.\h226.72.\h226.76.\h226.88.\h226.94.\h226.116.\h226.136.\h226.148.
\h226.160.\h227.7.\h227.11.\h227.15.\h227.19.\h227.21.\h227.23.\h227.24.
\h227.27.\h227.29.\h227.31.\h227.32.\h227.33.\h227.35.\h227.37.\h227.41.
\h227.43.\h227.53.\h227.59.\h227.62.\h227.65.\h227.71.\h227.73.\h227.87.
\h227.95.\h227.107.\h227.119.\h227.131.\h227.147.\h228.6.\h228.14.\h228.18.
\h228.21.\h228.24.\h228.26.\h228.28.\h228.30.\h228.36.\h228.42.\h228.44.
\h228.50.\h228.51.\h228.52.\h228.54.\h228.56.\h228.58.\h228.60.\h228.66.
\h228.70.\h228.78.\h228.90.\h228.102.\h228.118.\h228.130.\h228.138.\h229.13.
\h229.19.\h229.21.\h229.23.\h229.25.\h229.27.\h229.28.\h229.29.\h229.31.
\h229.33.\h229.35.\h229.39.\h229.40.\h229.43.\h229.49.\h229.53.\h229.59.
\h229.61.\h229.65.\h229.73.\h229.77.\h229.85.\h229.89.\h229.101.\h229.109.
\h229.117.\h229.127.\h230.10.\h230.20.\h230.22.\h230.24.\h230.26.\h230.28.
\h230.29.\h230.32.\h230.34.\h230.36.\h230.38.\h230.42.\h230.44.\h230.48.
\h230.52.\h230.56.\h230.60.\h230.62.\h230.72.\h230.76.\h230.80.\h230.88.
\h230.98.\h230.100.\h230.116.\h231.3.\h231.12.\h231.15.\h231.19.\h231.21.
\h231.24.\h231.25.\h231.26.\h231.27.\h231.31.\h231.33.\h231.34.\h231.35.
\h231.37.\h231.39.\h231.41.\h231.43.\h231.45.\h231.47.\h231.51.\h231.55.
\h231.57.\h231.59.\h231.63.\h231.69.\h231.71.\h231.75.\h231.87.\h231.99.
\h231.101.\h231.108.\h231.111.\h232.10.\h232.14.\h232.16.\h232.19.\h232.22.
\h232.24.\h232.25.\h232.26.\h232.30.\h232.32.\h232.34.\h232.37.\h232.38.
\h232.40.\h232.42.\h232.44.\h232.46.\h232.50.\h232.52.\h232.58.\h232.61.
\h232.70.\h232.72.\h232.74.\h232.79.\h232.82.\h232.97.\h232.100.\h232.106.
\h233.5.\h233.9.\h233.13.\h233.17.\h233.18.\h233.23.\h233.25.\h233.26.
\h233.27.\h233.29.\h233.31.\h233.33.\h233.35.\h233.37.\h233.38.\h233.41.
\h233.43.\h233.45.\h233.49.\h233.50.\h233.53.\h233.57.\h233.59.\h233.68.
\h233.71.\h233.74.\h233.77.\h233.86.\h233.93.\h233.95.\h233.98.\h233.99.
\h234.12.\h234.15.\h234.16.\h234.18.\h234.20.\h234.24.\h234.26.\h234.28.
\h234.32.\h234.33.\h234.34.\h234.36.\h234.39.\h234.40.\h234.42.\h234.48.
\h234.57.\h234.60.\h234.64.\h234.66.\h234.69.\h234.70.\h234.72.\h234.78.
\h234.84.\h234.92.\h234.96.\h235.11.\h235.13.\h235.15.\h235.19.\h235.22.
\h235.23.\h235.25.\h235.27.\h235.28.\h235.29.\h235.31.\h235.35.\h235.37.
\h235.40.\h235.41.\h235.43.\h235.49.\h235.51.\h235.55.\h235.59.\h235.63.
\h235.67.\h235.71.\h235.73.\h235.77.\h235.83.\h235.91.\h235.95.\h235.115.
\h235.139.\h236.10.\h236.11.\h236.14.\h236.18.\h236.20.\h236.22.\h236.24.
\h236.26.\h236.28.\h236.29.\h236.32.\h236.34.\h236.38.\h236.42.\h236.48.
\h236.50.\h236.54.\h236.56.\h236.59.\h236.62.\h236.66.\h236.70.\h236.72.
\h236.82.\h236.86.\h236.92.\h236.94.\h236.110.\h237.9.\h237.12.\h237.15.
\h237.17.\h237.19.\h237.21.\h237.25.\h237.27.\h237.31.\h237.33.\h237.37.
\h237.39.\h237.41.\h237.48.\h237.49.\h237.53.\h237.55.\h237.57.\h237.62.
\h237.63.\h237.65.\h237.69.\h237.71.\h237.77.\h237.81.\h237.87.\h237.93.
\h238.4.\h238.10.\h238.16.\h238.18.\h238.24.\h238.25.\h238.26.\h238.28.
\h238.33.\h238.34.\h238.37.\h238.38.\h238.40.\h238.48.\h238.50.\h238.52.
\h238.54.\h238.58.\h238.64.\h238.68.\h238.76.\h238.82.\h238.94.\h239.9.
\h239.15.\h239.19.\h239.20.\h239.23.\h239.25.\h239.26.\h239.27.\h239.29.
\h239.31.\h239.37.\h239.41.\h239.47.\h239.49.\h239.51.\h239.53.\h239.55.
\h239.56.\h239.63.\h239.65.\h239.69.\h239.71.\h239.73.\h239.75.\h239.91.
\h240.9.\h240.18.\h240.22.\h240.26.\h240.30.\h240.36.\h240.37.\h240.40.
\h240.44.\h240.45.\h240.46.\h240.48.\h240.54.\h240.60.\h240.62.\h240.68.
\h240.70.\h240.72.\h240.90.\h240.114.\h241.17.\h241.19.\h241.21.\h241.25.
\h241.27.\h241.29.\h241.33.\h241.34.\h241.36.\h241.37.\h241.39.\h241.41.
\h241.43.\h241.45.\h241.49.\h241.52.\h241.53.\h241.55.\h241.61.\h241.65.
\h241.67.\h241.69.\h241.70.\h241.85.\h241.89.\h241.113.\h242.8.\h242.14.
\h242.16.\h242.20.\h242.24.\h242.25.\h242.26.\h242.32.\h242.34.\h242.36.
\h242.38.\h242.40.\h242.41.\h242.44.\h242.48.\h242.50.\h242.52.\h242.56.
\h242.60.\h242.66.\h242.68.\h242.84.\h242.86.\h242.112.\h243.3.\h243.15.
\h243.23.\h243.24.\h243.27.\h243.28.\h243.31.\h243.35.\h243.37.\h243.39.
\h243.43.\h243.45.\h243.47.\h243.51.\h243.55.\h243.57.\h243.60.\h243.63.
\h243.65.\h243.67.\h243.83.\h243.93.\h244.10.\h244.14.\h244.18.\h244.21.
\h244.22.\h244.24.\h244.25.\h244.26.\h244.28.\h244.31.\h244.34.\h244.36.
\h244.40.\h244.46.\h244.48.\h244.50.\h244.54.\h244.56.\h244.62.\h244.64.
\h244.66.\h244.74.\h244.76.\h244.88.\h244.112.\h245.13.\h245.17.\h245.21.
\h245.23.\h245.25.\h245.29.\h245.33.\h245.35.\h245.37.\h245.39.\h245.45.
\h245.47.\h245.49.\h245.55.\h245.57.\h245.59.\h245.61.\h245.63.\h245.65.
\h245.67.\h245.83.\h245.85.\h246.12.\h246.18.\h246.20.\h246.24.\h246.26.
\h246.28.\h246.30.\h246.32.\h246.34.\h246.36.\h246.44.\h246.46.\h246.50.
\h246.54.\h246.56.\h246.57.\h246.58.\h246.60.\h246.66.\h246.84.\h246.111.
\h246.138.\h247.7.\h247.17.\h247.21.\h247.23.\h247.25.\h247.27.\h247.28.
\h247.29.\h247.31.\h247.33.\h247.35.\h247.39.\h247.43.\h247.45.\h247.47.
\h247.49.\h247.55.\h247.57.\h247.65.\h247.82.\h247.109.\h248.6.\h248.11.
\h248.20.\h248.22.\h248.24.\h248.26.\h248.28.\h248.29.\h248.32.\h248.35.
\h248.38.\h248.40.\h248.42.\h248.44.\h248.46.\h248.48.\h248.53.\h248.56.
\h248.65.\h248.80.\h248.83.\h248.92.\h248.110.\h249.9.\h249.19.\h249.21.
\h249.24.\h249.25.\h249.27.\h249.31.\h249.33.\h249.34.\h249.35.\h249.37.
\h249.48.\h249.51.\h249.53.\h249.54.\h249.63.\h249.75.\h249.81.\h249.109.
\h249.137.\h249.165.\h250.8.\h250.10.\h250.16.\h250.19.\h250.22.\h250.24.
\h250.25.\h250.28.\h250.30.\h250.32.\h250.43.\h250.46.\h250.52.\h250.62.
\h250.64.\h250.70.\h250.80.\h250.108.\h250.136.\h251.5.\h251.11.\h251.17.
\h251.21.\h251.23.\h251.25.\h251.27.\h251.29.\h251.32.\h251.35.\h251.45.
\h251.47.\h251.51.\h251.53.\h251.59.\h251.63.\h251.79.\h251.91.\h251.107.
\h251.119.\h251.251.\h252.16.\h252.18.\h252.20.\h252.22.\h252.24.\h252.26.
\h252.42.\h252.44.\h252.46.\h252.50.\h252.51.\h252.54.\h252.62.\h252.74.
\h252.78.\h252.90.\h252.102.\h252.222.\h253.13.\h253.17.\h253.19.\h253.21.
\h253.23.\h253.25.\h253.31.\h253.37.\h253.40.\h253.43.\h253.45.\h253.49.
\h253.61.\h253.73.\h253.89.\h253.193.\h254.16.\h254.18.\h254.20.\h254.30.
\h254.34.\h254.36.\h254.38.\h254.42.\h254.44.\h254.52.\h254.56.\h254.60.
\h254.72.\h254.80.\h254.164.\h254.176.\h255.15.\h255.21.\h255.23.\h255.24.
\h255.27.\h255.30.\h255.33.\h255.35.\h255.39.\h255.41.\h255.43.\h255.51.
\h255.59.\h255.69.\h255.71.\h255.135.\h255.147.\h255.159.\h256.14.\h256.22.
\h256.26.\h256.30.\h256.32.\h256.34.\h256.36.\h256.40.\h256.42.\h256.58.
\h256.64.\h256.106.\h256.118.\h256.130.\h256.146.\h257.13.\h257.25.\h257.29.
\h257.35.\h257.41.\h257.43.\h257.50.\h257.53.\h257.55.\h257.57.\h257.77.
\h257.89.\h257.101.\h257.117.\h257.129.\h257.137.\h258.12.\h258.20.\h258.24.
\h258.26.\h258.27.\h258.28.\h258.38.\h258.39.\h258.42.\h258.48.\h258.52.
\h258.60.\h258.72.\h258.84.\h258.88.\h258.100.\h258.108.\h258.116.\h258.126.
\h259.19.\h259.25.\h259.35.\h259.37.\h259.41.\h259.43.\h259.47.\h259.55.
\h259.59.\h259.71.\h259.79.\h259.87.\h259.97.\h259.99.\h259.115.\h260.14.
\h260.23.\h260.26.\h260.33.\h260.34.\h260.36.\h260.38.\h260.40.\h260.42.
\h260.44.\h260.50.\h260.54.\h260.58.\h260.62.\h260.68.\h260.70.\h260.86.
\h260.98.\h260.100.\h260.107.\h260.110.\h261.9.\h261.13.\h261.18.\h261.21.
\h261.23.\h261.24.\h261.29.\h261.33.\h261.39.\h261.41.\h261.45.\h261.51.
\h261.57.\h261.69.\h261.71.\h261.78.\h261.81.\h261.96.\h261.99.\h261.105.
\h262.12.\h262.16.\h262.22.\h262.28.\h262.32.\h262.36.\h262.40.\h262.42.
\h262.49.\h262.52.\h262.56.\h262.67.\h262.70.\h262.76.\h262.85.\h262.92.
\h262.94.\h262.98.\h263.11.\h263.15.\h263.19.\h263.23.\h263.25.\h263.27.
\h263.31.\h263.32.\h263.35.\h263.38.\h263.39.\h263.41.\h263.47.\h263.56.
\h263.63.\h263.65.\h263.68.\h263.69.\h263.77.\h263.83.\h263.91.\h263.95.
\h264.10.\h264.12.\h264.18.\h264.22.\h264.24.\h264.26.\h264.28.\h264.30.
\h264.34.\h264.36.\h264.39.\h264.40.\h264.48.\h264.54.\h264.62.\h264.66.
\h264.70.\h264.76.\h264.82.\h264.90.\h265.10.\h265.21.\h265.23.\h265.25.
\h265.27.\h265.28.\h265.33.\h265.37.\h265.41.\h265.47.\h265.53.\h265.55.
\h265.61.\h265.65.\h265.69.\h265.81.\h265.85.\h266.8.\h266.18.\h266.20.
\h266.24.\h266.26.\h266.32.\h266.36.\h266.38.\h266.40.\h266.47.\h266.52.
\h266.54.\h266.56.\h266.61.\h266.62.\h266.68.\h266.76.\h266.80.\h267.9.
\h267.15.\h267.23.\h267.24.\h267.27.\h267.32.\h267.33.\h267.36.\h267.37.
\h267.39.\h267.47.\h267.51.\h267.53.\h267.57.\h267.75.\h268.8.\h268.18.
\h268.22.\h268.28.\h268.36.\h268.40.\h268.46.\h268.50.\h268.64.\h268.68.
\h268.72.\h268.74.\h269.17.\h269.21.\h269.35.\h269.39.\h269.43.\h269.44.
\h269.45.\h269.53.\h269.61.\h269.67.\h269.71.\h270.16.\h270.20.\h270.24.
\h270.32.\h270.33.\h270.36.\h270.38.\h270.42.\h270.52.\h270.60.\h270.64.
\h270.66.\h270.69.\h271.7.\h271.15.\h271.19.\h271.31.\h271.35.\h271.37.
\h271.39.\h271.40.\h271.51.\h271.55.\h271.59.\h271.65.\h271.67.\h272.2.
\h272.14.\h272.22.\h272.23.\h272.26.\h272.27.\h272.30.\h272.34.\h272.36.
\h272.38.\h272.50.\h272.54.\h272.56.\h272.59.\h272.62.\h272.64.\h273.13.
\h273.21.\h273.25.\h273.27.\h273.30.\h273.33.\h273.35.\h273.45.\h273.49.
\h273.53.\h273.61.\h273.63.\h274.12.\h274.22.\h274.24.\h274.32.\h274.34.
\h274.38.\h274.44.\h274.46.\h274.48.\h274.58.\h274.60.\h275.11.\h275.17.
\h275.23.\h275.29.\h275.31.\h275.33.\h275.43.\h275.45.\h275.49.\h275.53.
\h275.55.\h275.56.\h275.57.\h275.59.\h276.24.\h276.26.\h276.27.\h276.28.
\h276.30.\h276.32.\h276.38.\h276.42.\h276.48.\h276.54.\h276.56.\h277.21.
\h277.25.\h277.27.\h277.31.\h277.37.\h277.43.\h277.47.\h277.52.\h277.55.
\h278.20.\h278.23.\h278.26.\h278.30.\h278.32.\h278.34.\h278.36.\h278.47.
\h278.50.\h278.52.\h278.53.\h279.21.\h279.23.\h279.24.\h279.27.\h279.29.
\h279.31.\h279.42.\h279.45.\h279.51.\h280.16.\h280.22.\h280.24.\h280.26.
\h280.28.\h280.31.\h280.44.\h280.46.\h280.50.\h281.21.\h281.23.\h281.25.
\h281.41.\h281.45.\h281.49.\h282.16.\h282.20.\h282.22.\h282.30.\h282.36.
\h282.39.\h282.44.\h282.48.\h283.15.\h283.19.\h283.29.\h283.35.\h283.37.
\h283.43.\h284.14.\h284.29.\h284.34.\h284.40.\h284.42.\h285.13.\h285.21.
\h285.29.\h285.33.\h285.35.\h285.39.\h285.41.\h286.28.\h286.34.\h286.40.
\h287.23.\h287.25.\h287.26.\h287.27.\h287.37.\h287.38.\h287.47.\h288.18.
\h288.24.\h288.34.\h288.36.\h288.42.\h289.13.\h289.22.\h289.33.\h289.35.
\h289.41.\h290.12.\h290.20.\h290.28.\h290.32.\h290.38.\h290.40.\h291.11.
\h291.15.\h291.21.\h291.27.\h291.31.\h291.39.\h292.10.\h292.22.\h292.26.
\h292.30.\h292.37.\h292.38.\h293.25.\h293.27.\h293.29.\h293.33.\h293.35.
\h293.38.\h294.9.\h294.20.\h294.24.\h294.26.\h294.32.\h294.36.\h295.7.
\h295.19.\h295.25.\h295.31.\h295.35.\h295.37.\h296.22.\h296.23.\h296.26.
\h296.32.\h296.35.\h296.36.\h297.21.\h297.27.\h297.35.\h298.20.\h298.34.
\h299.19.\h299.31.\h299.32.\h299.35.\h300.14.\h300.18.\h300.30.\h300.34.
\h301.13.\h301.25.\h301.26.\h301.29.\h301.33.\h302.12.\h302.20.\h302.24.
\h302.26.\h302.32.\h303.23.\h303.33.\h304.28.\h304.30.\h304.32.\h305.23.
\h305.25.\h305.31.\h306.24.\h306.30.\h307.19.\h307.22.\h307.29.\h307.31.
\h308.20.\h308.23.\h308.26.\h308.28.\h309.21.\h309.23.\h309.25.\h309.30.
\h310.20.\h310.22.\h310.24.\h311.19.\h311.21.\h311.29.\h312.18.\h312.28.
\h313.13.\h314.28.\h315.27.\h316.22.\h316.24.\h316.26.\h317.17.\h317.23.
\h318.12.\h318.21.\h319.11.\h319.19.\h320.10.\h320.26.\h321.9.\h321.21.
\h321.25.\h322.24.\h323.23.\h324.18.\h325.22.\h325.25.\h326.20.\h327.19.
\h328.18.\h329.17.\h330.12.\h330.24.\h331.23.\h332.22.\h334.22.\h335.23.
\h336.21.\h337.19.\h338.20.\h338.22.\h339.19.\h339.21.\h340.18.\h340.20.
\h341.17.\h345.21.\h346.16.\h347.11.\h347.20.\h348.10.\h348.18.\h350.20.
\h355.19.\h356.18.\h357.17.\h358.16.\h366.18.\h369.17.\h370.16.\h375.15.
\h376.10.\h377.17.\h386.16.\h387.15.\h399.15.\h404.14.\h416.14.\h433.13.
\h462.12.\h491.11.
}

\def\MaxCHC{24}                                         \def\Msum{150}
\vbox to 99pt   {\vss   \unitlength=2.7pt               
\begin{center}  \def\putlin#1,#2,#3,#4,#5){\put#1,#2){\line(#3,#4){#5}}}
                \def\putvec#1,#2,#3,#4,#5){\put#1,#2){\vector(#3,#4){#5}}}
        \def\Xlab#1 {\put(#1,#1){\drawline(-.5,.5)(.5,-.5)}}
        \def\Ylab#1 {\put(-#1,#1){\drawline(.5,.5)(-.5,-.5)}}

\unitlength=2.4pt               \def\CSH{.4} \def\CSC{.8}       
\BP(80,25)(35,2)        \putvec(0,0,1,0,153)    \putvec(0,0,0,1,24)
        \put(1,\MaxCHC){\put(0,1){$h_{12}$}}  \put(\Msum,-3.8){$\!\!\!h_{11}$}
        \put(1,9.3){\tiny10}\put(8.2,-3.5){\tiny10}\put(48.2,-3.5){\tiny50}
        \put(97.5,-3.5){\tiny100}
        \def\Xlab#1 {\put(#1,0){\drawline(0,-.5)(0,.5)}}
        \def\Ylab#1 {\put(0,#1){\drawline(.5,0)(-.5,0)}}
        \Xlab10 \Xlab20 \Xlab30 \Xlab40 \Xlab50 \Xlab60 \Xlab70 \Xlab80 
        \Xlab90 \Xlab100 \Xlab110 \Xlab120 \Xlab130 \Xlab140 \Xlab150   
        \Ylab10 \Ylab20
\def\h#1.#2.{ \ifnum \MaxCHC < #1 \else \put(#2,#1){\blue \circle{\CSC}} \fi }

\FanoCICYs
\def\h#1.#2.{ \ifnum \MaxCHC < #1 \else \ifnum #2 > \Msum \else
                        \put(#2,#1){\black\circle*{\CSH}} \fi \fi
              \ifnum \MaxCHC < #2 \else \ifnum #1 > \Msum \else
                        \put(#1,#2){\black\circle*{\CSH}} \fi \fi}

\HypersurfCY
\EP\end{center}
                                        \stepcounter{figure}    \makeatletter
\immediate\write\@auxout{\string\newlabel{Conifolds}{\thefigure}}\makeatother
\vspace{14pt}
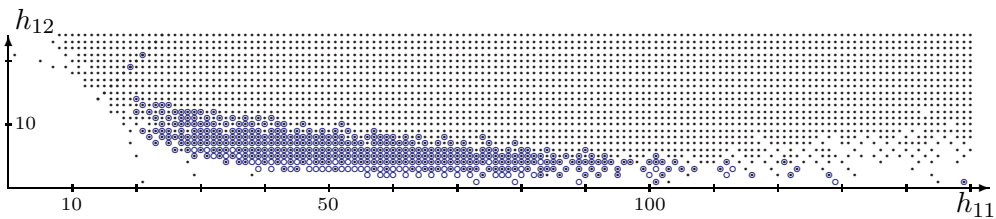
\centerline{Figure: CICYs from nef-partitions (circles) and toric 
        hypersurfaces (dots).}
\vspace{-5pt}}

\end{document}